\documentclass[12pt,oneside]{amsart}
\usepackage[margin=3cm]{geometry}
\usepackage[utf8]{inputenc}
\usepackage{amsmath}
\usepackage{tikz-cd}
\usepackage{bm}
\usepackage{amssymb}
\usepackage{mathtools}
\usepackage{amsthm}
\usepackage{graphicx}
\usepackage{comment}
\usepackage{bbold}
\usepackage{comment}
\usepackage{appendix}
\usepackage{caption}
\usepackage{subcaption}
\usepackage{hyperref}
\usepackage{upgreek}
\usepackage{enumitem}
\usepackage{comment}
\usepackage{tikz}
\usepackage{pict2e,mathrsfs}
\usepackage[T1]{fontenc}

\title{Immaculate complexes and Strong Generation by Injectives}
\author{Isambard Goodbody}
\date{\today}

\DeclareUnicodeCharacter{0327}{\c{}}

\DeclareMathOperator*{\Inj}{Inj}

\DeclareMathOperator{\coloc}{coloc}
\DeclareMathOperator*{\coker}{coker}
\DeclareMathOperator*{\hocolim}{hocolim}

\DeclareMathOperator*{\colim}{colim}

\DeclareMathOperator{\im}{im}

\DeclareMathOperator{\Hom}{Hom}

\DeclareMathOperator{\RHom}{RHom}
\DeclareMathOperator{\End}{End}
\DeclareMathOperator{\loc}{loc}

\DeclareMathOperator{\thick}{thick}

\DeclareMathOperator{\Mod}{Mod}

\DeclareMathOperator{\fdmod}{mod}
\DeclareMathOperator{\coh}{coh}
\DeclareMathOperator{\QCoh}{QCoh}
\DeclareMathOperator{\perf}{perf}

\DeclareMathOperator{\Ext}{Ext}

\DeclareMathOperator{\id}{id}

\DeclareMathOperator{\idim}{idim}
\DeclareMathOperator{\gldim}{gldim}

\DeclareMathOperator{\inj}{inj}
\DeclareMathOperator{\minj}{maxinj}

\DeclareMathOperator{\supp}{supp}

\DeclareMathOperator{\immac}{im}
\DeclareMathOperator{\mimmac}{mim}
\DeclareMathOperator{\ipscoh}{ipc}
\DeclareMathOperator{\ifp}{ifp}
\DeclareMathOperator{\pc}{pc}

\DeclareMathOperator{\codim}{codim}
\DeclareMathOperator{\Add}{Add}
\DeclareMathOperator{\Prod}{Prod}
\DeclareMathOperator{\hgt}{ht}
\DeclareMathOperator{\Spec}{Spec}
\DeclareMathOperator{\Frac}{Frac}

\DeclareMathOperator{\sg}{sg}
\DeclareMathOperator{\isg}{isg}
\DeclareMathOperator{\Ab}{Ab}
\DeclareMathOperator{\fsupp}{fsupp}

\begin{document}

\usetikzlibrary{matrix}
\usetikzlibrary{shapes}

\theoremstyle{plain}
\newtheorem{prop}{Proposition}[section]
\newtheorem{lemma}[prop]{Lemma}
\newtheorem{theorem}[prop]{Theorem}
\newtheorem{cor}[prop]{Corollary}
\newtheorem{quest}[prop]{Question}

\theoremstyle{definition}
\newtheorem{defn}[prop]{Definition}
\newtheorem{ass}[prop]{Assumption}
\newtheorem{conj}[prop]{Conjecture}
\newtheorem{ex}[prop]{Example}
\newtheorem{remark}[prop]{Remark}
\newtheorem{notation}[prop]{Notation}
\newtheorem*{ack}{Acknowledgements}
\newtheorem*{AID}{AI Declaration}

\begin{abstract}
We study the category of bounded complexes of finite sums of indecomposable injectives over a locally Noetherian scheme, which we call the immaculate complexes. We show that regularity and finiteness of dimension can be detected from this category using strong generation. Our methods involve studying when the unbounded derived category of a Grothendieck category is strongly generated by injectives. As applications we prove a conjecture of Hartshorne and recover a result of Keller and Krause. To any specialisation closed subset we associate a semi-orthogonal decomposition of the immaculates and use this to compute additive finitary invariants in terms of the completions of the local rings. The immaculates are equivalent to the perfect derived category of a left coherent preadditive category i.e.\ a ring with many objects and we compute this category for one-dimensional integral or Gorenstein schemes. 
\end{abstract}

\maketitle

\tikzset{
    vert/.style={anchor=south, rotate=90, inner sep=.5mm}
}

\tableofcontents

\section{Introduction}

The derived category of an algebraic or geometric object can be thought of as a highly structured invariant. The value of any invariant can be judged by two opposing factors:  how much information it contains about the original object and how easy it is to compute. Many well studied questions about derived categories can be framed in these terms. For example, a fundamental result in noncommutative algebraic geometry is that smoothness and properness of reasonable schemes can be detected from their derived categories \cite{Orl16}. Going further, reconstruction or categorical Torelli theorems, as in \cite{BO01}, occasionally allow us to reconstruct the entire scheme from its derived category. On the other hand, results that shed light on the internal structure of the derived category such as classification problems in tensor-triangular geometry \cite{Nee92} or moduli problems can be interpreted as computing the invariant.

In this article we introduce a new categorical invariant of schemes (or Grothendieck categories) and show that it contains significant geometric information and is amenable to computation. Morally dual to the perfect complexes, the immaculate complexes $\mathcal{D}^{\immac}(\mathcal{A})$ over a Grothendieck category $\mathcal{A}$ are the bounded complexes of finite sums of indecomposable injectives. The immaculates come with a set of well behaved generators: the indecomposable injectives. While a landmark result of \cite{VB02} shows that the perfect complexes over a reasonable scheme are generated by a single object, this generator is difficult to get one's hands on. Matlis's theory of injectives lets us control the generators of $\mathcal{D}^{\immac}(X)$ for locally Noetherian schemes. The cost is that instead of a single generator there is a (usually infinite) set of generators. It turns out that this can be dealt with using a certain finiteness property of injective resolutions (see Proposition \ref{finitenessprop}). 

Our first result is motivated by Conjecture 10 of \cite{Orl09} which states that one can extract the dimension of a smooth projective variety using the minimum number of steps needed to build any object of the perfect derived category from any generator. One difficulty with this conjecture is getting a handle on the possible generators of the derived category. We show that if one keeps track of the indecomposable injectives, then one can extract the dimension and regularity of reasonable schemes. For $n \geq 0$, and $S \subseteq \mathcal{T}$ a collection of objects in a triangulated category, we let $\thick_n(S) \subseteq \mathcal{T}$ denote the subcategory of objects built from $S$ is $n$ steps (see Definition \ref{buildingnotation} for a precise definition). For a scheme $X$, we let $\inj$ denote the indecomposable injective quasicoherent sheaves, which we view as complexes concentrated in degree zero. 

\begin{theorem}[Theorem \ref{regularscheme} b)] \label{introregularscheme1} Let $X$ be a locally Noetherian scheme. Then, $X$ is regular and of dimension at most $n$ if and only if $\mathcal{D}^{\immac}(X) = \thick_n(\inj)$.
\end{theorem}

Next, we consider regularity without reference to dimension and prove a version of a well-known result about the perfect complexes. We let $\mathcal{D}^b_{\ipscoh}(X)$ denote the bounded below complexes of finite sums of indecomposable injectives with bounded cohomology (Definition \ref{immacdefn}). The relationship between this category and the immaculate complexes is similar to the relationship between the bounded derived category of coherent sheaves and the perfect complexes.

\begin{theorem}[Theorem \ref{regularscheme} a)] \label{introregularscheme2}Let $X$ be a locally Noetherian scheme. Then $X$ is regular if and only if $\mathcal{D}^{\immac}(X) = \mathcal{D}^b_{\ipscoh}(X)$.
\end{theorem}

Our results follow from more general results about Grothendieck categories. It is useful to first consider generation in the big setting. If $S \subseteq \mathcal{T}$ is a collection of objects in a triangulated category with coproducts and products, we let $\loc_n(S)$ and $ \coloc_n(S)$ denote the subcategories of $\mathcal{T}$ built from $S$ in $n$ steps using all (including infinite) coproducts and products respectively (Definition \ref{bigbuildingnotation}). 

\begin{theorem}[Theorem \ref{bigproductglobdim}] \label{introbigglobdim} Let $\mathcal{A}$ denote a Grothendieck category. 
\begin{enumerate}
\item The global dimension of $\mathcal{A}$ is at most $n$ if and only if there is an injective $I$ such that $\mathcal{D}(\mathcal{A}) = \coloc_n(I)$.

\item Suppose $\mathcal{A}$ is locally Noetherian. Then the global dimension of $\mathcal{A}$ is at most $n$ if and only if there is an injective $I$ such that $\mathcal{D}(\mathcal{A}) = \loc_n(I)$, if and only if $\mathcal{D}(\mathcal{A}) = \loc_n(\inj)$. 
\end{enumerate}
\end{theorem}

We use Theorem \ref{introbigglobdim} to prove a conjecture (Theorem \ref{Bigstronggenregular}) of Hartshorne from \cite{Har66} about the difference between complexes of quasicoherent sheaves and complexes with quasicoherent cohomology for regular finite-dimensional schemes (not assuming quasicompactness). Another consequence of Theorem \ref{introbigglobdim} is on the behaviour of global dimension under derived equivalence for Grothendieck categories. 

\begin{theorem}[Theorem \ref{bigmoritaglobdim}]
Suppose $\mathcal{A}, \mathcal{B}$ are Grothendieck categories and there is an equivalence $G \colon \mathcal{D}(\mathcal{A}) \xrightarrow{\sim} \mathcal{D}(\mathcal{B})$ which preserves bounded complexes. Suppose that either 
\begin{enumerate}
\item $\mathcal{A}$ is locally Noetherian,

\item or products in $\mathcal{B}$ are exact.
\end{enumerate}
If $\mathcal{A}$ has finite global dimension, then so does $\mathcal{B}$. 
\end{theorem}

Specialising to reasonable schemes we show in Theorem \ref{bigmoritascheme} that being regular and finite-dimensional is derived invariant. We prove the converse to a result of \cite{BK20} about global dimension of endomorphism rings of tilting objects as well as giving a new proof of the original result.

\begin{theorem}[Theorem \ref{tiltingglobdim}] Let $\mathcal{A}$ denote a Grothendieck category and suppose $T \in \mathcal{D}(\mathcal{A})$ is a tilting object. Then the global dimension of $\mathcal{A}$ is finite if and only if the global dimension of $\End_{\mathcal{A}}(T)$ is finite.  
\end{theorem}

To pass to the small setting, that is to $\mathcal{D}^{\immac}(\mathcal{A})$, we need $\mathcal{A}$ to satisfy a condition dual to being coherent which we call $\inj$-coherence (Definition \ref{injcohdefn}). We call an object of $\mathcal{A}$ $\inj$-finitely presented if it admits a copresentation by finite sums of indecomposable injectives. 

\begin{theorem}[Theorem \ref{smallglobdim}] \label{introsmallglobdim}
Let $\mathcal{A}$ be an $\inj$-coherent Grothendieck category. Then 
\begin{enumerate}
\item $\mathcal{D}^{\immac}(\mathcal{A}) = \mathcal{D}^b_{\ipscoh}(\mathcal{A})$ if and only if every $\inj$-finitely presented object has finite injective dimension.

\item $\mathcal{D}^{\immac}(\mathcal{A}) = \thick_n(\inj)$ if and only if every $\inj$-finitely presented object has injective dimension at most $n$.
\end{enumerate}
\end{theorem}

So Theorem \ref{introsmallglobdim} reduces the strong generation question to the classical problem of determining which classes of objects in $\mathcal{A}$ can detect if $\mathcal{A}$ has finite global dimension. To apply Theorem \ref{introsmallglobdim} to locally Noetherian schemes we show that $\QCoh(X)$ is $\inj$-coherent. This is the key property underpinning most of the results in the paper and is proved in Section \ref{immaculateschemes}. Theorems \ref{introregularscheme1} and \ref{introregularscheme2} then follow from Theorem \ref{introsmallglobdim} as the residue fields are all $\inj$-finitely presented.

We note that $\inj$-coherence of $\mathcal{A}$ is equivalent to the standard $t$-structure restricting to $\mathcal{D}^b_{\ipscoh}(\mathcal{A})$ (Proposition \ref{injcohequiv}). This also provides a characterisation of regularity. 

\begin{theorem} [Theorem \ref{regularscheme} c)]
Let $X$ be a locally Noetherian scheme. Then $X$ is regular if and only if the standard $t$-structure restricts to $\mathcal{D}^{\immac}(X)$. 
\end{theorem}

As we have kept track of a set of  generators in Theorem \ref{introregularscheme1}, it is a priori not clear that this characterisation of regularity and dimension is inherent to $\mathcal{D}^{\immac}(X)$. To address this we prove the following.

\begin{theorem}[Theorem \ref{dimsame}] \label{introdimsame} Suppose $X$ and $Y$ are locally Noetherian schemes and $\mathcal{D}^{\immac}(X) \simeq \mathcal{D}^{\immac}(Y)$. Then $X$ is regular if and only if $Y$ is regular. If, in addition, $X$ has finite dimension, then $Y$ has finite dimension and $\dim(X) = \dim(Y)$.
\end{theorem}

The key ingredient for showing the dimensions coincide is the recent result of \cite{Let25} computing the Rouquier dimension of all commutative Noetherian rings.  We prove Theorem \ref{introdimsame} by showing that $\mathcal{D}^{\immac}(X)$ and $\mathcal{D}^b_{\ifp}(X)$ determine each other in the following sense. 

\begin{theorem}[Theorem \ref{determineschemes}]
Suppose $X$ and $Y$ are locally Noetherian schemes, then $\mathcal{D}^{\immac}(X) \simeq \mathcal{D}^{\immac}(Y)$ if and only if $\mathcal{D}^b_{\ifp}(X) \simeq \mathcal{D}^{b}_{\ifp}(Y)$.
\end{theorem}

This kind of problem for the perfects has been studied by many authors such as in \cite{Ric89}, \cite{Nee25}, \cite{KS22} and \cite{Kra20}. A key difference here is that the immaculates don't usually have a single generator.

The previously mentioned results show that the immaculates contain significant geometric information. We now move onto ways of understanding the internal structure of this category. One of these is via semi-orthogonal decompositions. We construct a semi-orthogonal decomposition for any specialisation closed subset of $X$ and in particular produce the following description. 

\begin{theorem}[Theorem \ref{SODfiltration}]
Let $X$ be a finite-dimensional locally Noetherian scheme. There is a semi-orthogonal decomposition
\[
\mathcal{D}^{\immac}(X) = \langle \mathcal{F}_0, \mathcal{F}_1,\dots \mathcal{F}_{\dim(X)} \rangle
\]
where $\mathcal{F}_i \simeq \bigoplus_{\dim \mathcal{O}_{X,x} = i }\mathcal{D}^{\perf}(\widehat{\mathcal{O}}_{X,x})$.
\end{theorem}

As a consequence, we compute any additive finitary invariant of $\mathcal{D}^{\immac}(X)$ in terms of the completions of the local rings. Examples include (non)connective $K$-theory and (topological) Hochschild homology. See \cite{BGT13} for more details on these invariants.

\begin{theorem}[Theorem \ref{additivethm}] \label{introKthy}
Let $X$ be a locally Noetherian scheme and $U$ any additive finitary invariant. Then there are equivalences
\[
U(\mathcal{D}^{\immac}(X)) \simeq \bigoplus_{x \in X} U(\mathcal{D}^{\perf}(\widehat{\mathcal{O}}_{X,x}))
\]

\end{theorem}

In particular, $K_0(\mathcal{D}^{\immac}(X))$ is free on the cardinality of $X$. Using negative $K$-theory, we consider when the immaculate singularity category $\mathcal{D}_{\isg}(X) \coloneq \mathcal{D}^b_{\ifp}(X)/ \mathcal{D}^{\immac}(X)$ is idempotent complete. Here $\mathcal{D}_{\sg}(\widehat{\mathcal{O}}_{X,x}) = \mathcal{D}^b_{\fdmod}(\widehat{\mathcal{O}}_{X,x})/\mathcal{D}^{\perf}(\widehat{\mathcal{O}}_{X,x})$ denotes the ordinary singularity category.

\begin{cor}[Corollary \ref{idempotentcomplete}]
Let $X$ be a locally Noetherian scheme. Then $\mathcal{D}_{\isg}(X)$ is idempotent complete if and only if $\mathcal{D}_{\sg}(\widehat{\mathcal{O}}_{X,x})$ is idempotent complete for every $x \in X$. 
\end{cor}

In Section \ref{preadditivesec}, we take a noncommutative viewpoint on $\mathcal{D}^{\immac}(\mathcal{A})$ as the indecomposable injectives provide a set of tilting generators. This allows us to describe $\mathcal{D}^{\immac}(\mathcal{A})$ using a small preadditive category. Recall that a preadditive category is a category enriched in Abelian groups and can be thought of as a ring with many objects. Many natural constructions and properties of rings generalise to preadditive categories. If $\mathcal{A}$ is a Grothendieck category, then we view the indecomposable injectives $\inj$ as a preadditive subcategory of $\mathcal{A}$.

\begin{theorem}[Theorems \ref{tiltingequiv1} and \ref{tiltingequiv2}]
If $\mathcal{A}$ is an $\inj$-coherent Grothendieck category, then $\inj$ is a left coherent preadditive category and the restricted Yoneda embedding $\mathcal{A}^{op} \to \Mod \inj^{op}$ induces equivalences
\[
\begin{tikzcd}
\mathcal{D}^b_{\ifp}(\mathcal{A})^{op} \arrow[r,"\sim"] & \mathcal{D}^b_{\fdmod}(\inj^{op}) \\
\mathcal{D}^{\immac}(\mathcal{A})^{op} \arrow[u,hook] \arrow[r,"\sim"] & \mathcal{D}^{\perf}(\inj^{op}) \arrow[u,hook]
\end{tikzcd}
\]
and the top equivalence is $t$-exact with respect to the standard $t$-structures.
\end{theorem}

Specialising to locally Noetherian schemes, we can combine our results to make a statement which does not mention derived categories. We recall that weak dimension is analogous to global dimension but uses flats instead of injectives or projectives and is left-right symmetric.

\begin{theorem}[Theorem \ref{weakdimscheme}] \label{introweakdimscheme} Let $X$ be a locally Noetherian scheme. Then 
\begin{enumerate}
\item $\inj(X)$ is a left coherent preadditive category.
\item $X$ is regular if and only if every finitely presented left $\inj(X)$ module has finite projective dimension.
\item $X$ is regular of dimension $n$ if and only if $\inj(X)$  has weak dimension $n$.
\end{enumerate}
\end{theorem}

Examples (see Section \ref{curvessec}) show that $\inj(X)$ often fails to be left Noetherian and that weak dimension cannot be replaced by left global dimension in the statement of Theorem \ref{introweakdimscheme}. 
It was shown by Rickard in \cite{Ric89} that a commutative ring is determined by its perfect derived category. In Section \ref{reconstructionsec}, we show that a Noetherian affine scheme can be recovered from its category of immaculates as its zero-th Hochschild cohomology. As a consequence we deduce the following. 

\begin{theorem}[Corollary \ref{recoveraffine}] Let $R,S$ be commutative Noetherian rings such that $\mathcal{D}^{\immac}(R) \simeq \mathcal{D}^{\immac}(S)$,  then $R \cong S$. 
\end{theorem}

In Section \ref{curvessec}, we explicitly describe $\inj(X)$ for one-dimensional schemes which are either Gorenstein or integral. We make use of a result of \cite{Sch15}, which describes the morphisms between indecomposable injectives for one-dimensional domains. We use these results to explicitly describe the category $\inj$ for many natural examples. For example (Example  \ref{injegslist}), if $C$ is a smooth projective curve over an algebraically closed field $k$, $\inj(X)$ has the form

\[
\begin{tikzcd}[column sep = huge]
\eta \arrow[r,"{k((t_x))}"] \arrow[loop left,"K(C)"] & x \arrow[loop right, "{k[[t_x]]}"]
\end{tikzcd}
\]
where $x$ varies over the closed points and $\eta$ is the generic point. For these small examples we also consider other properties of $\inj(X)$ such as right coherence and being Noetherian (Proposition \ref{notnoetheg}). 

In Section \ref{mimmacsec} we consider a variant called the maximal immaculates $\mathcal{D}^{\mimmac}(\mathcal{A})$ over a Grothendieck category $\mathcal{A}$. Here the indecomposable injectives are replaced with injective hulls of the simple objects. The analogous notion of $\inj$-coherence (or really $\inj$-Noetherianness) has been previously studied under the name co-Noetherian in noncommutative ring theory (see \cite{Jan69}, \cite{Vam68}). We relate strong generation to an invariant $\delta(R)$ of a ring $R$ introduced in \cite{Fac82} (Definition \ref{deltainvariant}). 

\begin{theorem}
Let $R$ be a co-Noetherian ring. Then $\delta(R) \leq n$ if and only if
\[
\mathcal{D}^{\mimmac}(R) = \thick_n(E(S) \mid S \text{ simple})
\]
\end{theorem}

In Section \ref{Matlissec}, we recall some facts about injective sheaves over locally Noetherian schemes. In Sections \ref{bigstronggensec}, \ref{derivedequivsec}, \ref{immacsec} and \ref{preadditivesec} we work in the abstract setting of Grothendieck categories. In Sections \ref{immaculateschemes}, \ref{perfimmacsec}, \ref{SODsec}, \ref{immacequivsec} and \ref{reconstructionsec} we work with locally Noetherian schemes. In Section \ref{curvessec}, we consider one-dimensional examples. Section \ref{mimmacsec} studies maximal immaculates in the context of noncommutative ring theory.   
\\

\begin{ack}
I'm extremely grateful to Ben Briggs, Timothy De Deyn, Martin Gallauer, Henning Krause, Kabeer Manali-Rahul, Greg Stevenson and Jan Šťovíček for many useful comments, suggestions, questions and their reading of drafts. In particular to Greg Stevenson for discussions about strong generation which started the project and to Henning Krause for asking questions which led to the results of Sections \ref{derivedequivsec} and \ref{reconstructionsec}. This project was funded by the London Mathematical Society Early Career Research Fellowship ECF-2026-15. This work was started as a PhD student supported by the Carnegie Trust for the
Universities of Scotland Grant Number PHD010651. 
\end{ack}

\begin{AID}
Claude Opus 5.5, Chat GPT Free and Google Gemini were used for literature review. None of the ideas, proofs or writing were produced by AI. 
\end{AID}

\section{Matlis's theory of Injectives}
\label{Matlissec}

For a scheme $X$, $\QCoh(X)$ will denote the quasicoherent sheaves on $X$ and $\Mod \mathcal{O}_X$, the $\mathcal{O}_X$-modules. By results of Gabber (\cite[Tag 077K]{stacks-project}) and Grothendieck (\cite[3.1.1]{Groth57}) both are Grothendieck categories. If $X$ is locally Noetherian, there is a good structural theory of the injectives in $\QCoh(X)$ due to Matlis \cite{Ma58}, Gabriel \cite{Gab62} and Hartshorne \cite{Har66}. See also \cite{Con00}. We summarise these below.

\begin{theorem} \label{matlisthy}
Let $X$ be a locally Noetherian scheme.
\begin{enumerate}
\item The injectives in $\QCoh(X)$ are exactly the injectives in $\Mod \mathcal{O}_X$ which are quasicoherent.

\item Direct sums of injectives in $\QCoh(X)$ are injective.

\item There is a bijection between isomorphism classes of indecomposable injectives in $\QCoh(X)$ and $X$, where $x \in X$ corresponds to $E(x) \coloneq E(k(x))$ the injective hull of the residue field of $x$.

\item Every injective in $\QCoh(X)$ is a direct sum of indecomposable injectives. In fact if $M \in \QCoh(X)$, then $E(M) = \bigoplus_{x \in X} E(x)^{\mu_0(x,M)}$ where
\[
\mu_0(x,M) = \mu_0^X(x,M) \coloneq \dim_{k(x)} \Hom_{\mathcal{O}_{X,x}}(k(x),M_x)
\]

\end{enumerate}
\end{theorem}

\begin{remark} \label{locNoethrem}
If $\mathcal{A}$ is any locally Noetherian Grothendieck category, then (2) and the first part of (4) hold for $\mathcal{A}$ by \cite[IV.2 Theorem 2, IV.2 Proposition 6]{Gab62} but the converse is not quite true. In fact, even if $X$ is locally Noetherian, $\QCoh(X)$ is not necessarily a locally Noetherian Grothendieck category (\cite[pp.\ 135]{Har66}). Although, $\Mod \mathcal{O}_X$ is locally Noetherian by \cite[Prop 7.8]{Har66} and if $X$ is additionally quasicompact, then $\QCoh(X)$ is locally Noetherian by \cite[VI.2]{Gab62}. 
\end{remark}

The following statements can be deduced from \cite[Chapter II Section 7]{Har66} and \cite{Ma58}. 

\begin{prop} \label{matlisprop}Let $X$ be a locally Noetherian scheme and $x,y \in X$.
\begin{enumerate}
\item If $y \not \in \overline{\{x\}}$ then $\Hom_X(k(x),E(y)) = \Hom_X(E(x),E(y)) = 0$

\item $\Hom_X(E(x),E(x)) \simeq \widehat{\mathcal{O}}_{X,x}$

\item  $\Hom_X(k(x),E(x)) = k(x)$

\item If $x \not\in \overline{\{ y \}}$, then $E(y)_x = 0$
\item If $x \in \overline{\{ y \}}$, $E(y)_x \simeq E_{\mathcal{O}_{X,x}}(y) \in \Mod \mathcal{O}_{X,x}$ where $E_{\mathcal{O}_{X,x}}(y)$ denotes the injective hull of the residue field $k(y)$ over $\mathcal{O}_{X,x}$.
\end{enumerate}
\end{prop}

As a consequence, the Auslander-Buchsbaum-Serre Theorem for commutative Noetherian rings extends to locally Noetherian schemes. 

\begin{theorem} \label{ABSscheme}
Suppose $X$ is a locally Noetherian scheme. Then $\QCoh(X)$ has global dimension $n$ if and only if $X$ is regular of dimension $n$.
\end{theorem}

This is surely well-known but we include a proof. If $I$ is a complex of injectives in a Grothendieck category $\mathcal{A}$, then we say $I$ is homotopically minimal if for each $n$ $I^n = E(\ker(d^n: I^n \to I^{n+1}))$. By \cite[Proposition B.2]{kra05} any complex of injectives is homotopy equivalent to a homotopically minimal complex of injectives. We say an injective resolution $I$ of an object $M \in \mathcal{A}$ is minimal if it is homotopically minimal.

\begin{lemma} \label{minimallemma}
Let $\mathcal{A}$ be a Grothendieck category, $I$ a homotopically minimal complex of injectives and $S \in \mathcal{A}$ a simple object. Then the differential on $\Hom_{\mathcal{A}}(S,I)$ vanishes.
\end{lemma}

\begin{proof}
It is enough to show that $d^n f = 0$ for any $f: S \to I^n$. As $S$ is simple, $f$ is mono. Since $\ker(d^n)$ is essential in $I^n$, $f(S) \cap \ker(d^n) \neq 0$. As $S$ is simple $f(S) \subseteq \ker(d^n)$ as required.
\end{proof}

\begin{lemma} \label{localismininj}
If $I$ is a homotopically minimal complex of injectives in $\QCoh(X)$, and $x \in X$, then $I_x$ is a homotopically minimal complex of injectives $\Mod \mathcal{O}_{X,x}$.
\end{lemma}

\begin{proof}
Let $Z^n$ denote the $n$-th kernel of $I$ so that $I^n = E(Z^n)$. It remains to show that $I^n_x = E(Z^n_x)$. By Theorem \ref{matlisthy} (4), $I^n \simeq \bigoplus_{y \in X} E(y)^{\mu^X_0(y,Z^n)}$. By Proposition \ref{matlisprop} (4) and (5), taking stalks gives 
\[
I^n_x \simeq \bigoplus_{x \in \overline{\{y\}}}  E_{\mathcal{O}_{X,x}}(y)^{\mu_0^X(y,Z^n)}
\]
Alternatively we can compute $E(Z^n_x)$ in $\Mod \mathcal{O}_{X,x}$ using the formula in Theorem \ref{matlisthy} (4). Identify $\Spec(\mathcal{O}_{X,x}) \simeq \{y \in X \mid x \in \overline{\{y\}}\}$ and note that $\Hom_{\mathcal{O}_{X,y}}(k(y),(Z^n_x)_y) \simeq \Hom_{\mathcal{O}_{X,y}}(k(y),Z^n_y)$. Therefore $\mu_0^{\mathcal{O}_{X,x}}(y,Z^n_x) = \mu_0^{X}(y,Z^n)$ and so $E(Z_x^n) \simeq I_x^n$.  
\end{proof}

\begin{lemma} \label{bassnumbers} Let $X$ be a locally Noetherian scheme, $M \in \QCoh(X)$ and suppose $M \to I$ is a minimal injective resolution of $M$ in $\QCoh(X)$. Then $I^i \simeq \bigoplus_{x \in X} E(x)^{\mu_i^X(x,M)}$ where 
\[
\mu_i^X(x,M) = \dim_{k(x)} \Ext^i_{\mathcal{O}_{X,x}}(k(x),M_x)
\]
\end{lemma}

\begin{proof} 
We have that $\Ext^i_{\mathcal{O}_{X,x}}(k(x),M_x) = H^i(\Hom_{\mathcal{O}_{X,x}}(k(x),I_x))$. By Lemmas \ref{minimallemma} and \ref{localismininj}, the differential on $\Hom_{\mathcal{O}_{X,x}}(k(x),I_x)$ vanishes. Therefore $\Ext^i_{\mathcal{O}_{X,x}}(k(x),M_x)$ is given by
\[
 \Hom_{\mathcal{O}_{X,x}}(k(x),\bigoplus_{x \in \overline{\{y\}}} E(y)^{\mu_i(y,M)}) = \Hom_{\mathcal{O}_{X,x}}(k(x),E(x)^{\mu_i(x,M)}) = k(x)^{\mu_i(x,M)}
\]
using Proposition \ref{matlisprop} (1) and (3).
\end{proof}

\begin{proof}[Proof of Theorem \ref{ABSscheme}]
If $X$ is regular of dimension at most $n$, then all the local rings have global dimension at most $n$. So for any $M$ and $x \in X$, $\mu_i(x,M) = 0$ for $i > n$. So by Lemma \ref{bassnumbers}, the minimal injective resolution of $M$ is of length at most $n$.  Conversely, if $\QCoh(X)$ has global dimension at most $n$, then each of the residue fields $k(x)$ admit injective resolutions of length at most $n$. Localization is exact and by Proposition \ref{matlisprop} 4) it preserves injectives, so $k(x)$ has finite injective dimension over $\mathcal{O}_{X,x}$. Hence each of the local rings are regular of dimension at most $n$. 
\end{proof}

It will be useful to note that the Bass numbers are computed stalkwise in the following sense.

\begin{prop} \label{stalkwise}
If $X$ is a locally Noetherian scheme, $M \in \QCoh(X)$ and $x,y \in X$ such that $x \in \overline{\{y\}}$. Then
\[
\mu_i^X(y,M) = \mu_i^{\mathcal{O}_{X,x}}(y,M_x)
\]
\end{prop}

\begin{proof}
Let $M \to I$ denote a minimal injective resolution of $M$, then $M_x \to I_x$ is a minimal injective resolution of $M_x$ by Lemma \ref{localismininj} and $I_x^i = \bigoplus_{x \in \overline{\{y\}}} E(y)^{\mu_i^X(y,M)}$. So by Lemma \ref{localismininj} applied to $\mathcal{O}_{X,x}$ we are done. 
\end{proof}

We introduce the following definition to abstract the properties of $\QCoh(X)$ that we will need later.

\begin{defn} \label{Matliscat}
We say a Grothendieck category $\mathcal{A}$ is Matlis if every injective in $\mathcal{A}$ is a direct sum of indecomposable injectives.
\end{defn}

\begin{remark} \label{matlisfacts}
\begin{enumerate}
\item Since the endomorphism rings of indecomposable injectives are local by \cite[II.5 Proposition 11]{Gab62}, the decompositions in a Matlis category are unique in the usual sense.

\item By Remark \ref{locNoethrem}, any locally Noetherian Grothendieck category is Matlis.

\item By Theorem \ref{matlisthy}, $\QCoh(X)$ is Matlis if $X$ is locally Noetherian.

\item By \cite[Theorem 2]{Ste68}, if $\mathcal{A}$ is Matlis, then direct sums of injectives in $\mathcal{A}$ are injective.

\item By \cite[IV.1 Theorem 1]{Gab62}, if $\mathcal{A}$ is a Grothendieck category with  Krull-Gabriel dimension, then direct sums of injective being injective implies $\mathcal{A}$ is Matlis. 

\item By \cite[Theorem 3]{Ste68}, if $\mathcal{A}$ is a Grothendieck category satisfying the AB6 axiom, then being Matlis is equivalent to being locally Noetherian. 

\end{enumerate}
\end{remark}

\section{Big Strong Generation by Injectives}
\label{bigstronggensec}

In this section, we show that one can extract the global dimension of a Grothendieck category from its unbounded derived category using strong generation by injectives. Our results are similar to those of \cite{KK05} but applied to injectives. Their results are related to those of \cite{Rou08} although in order to extract the precise dimension one needs to use the arguments of \cite{KK05}. Small strong generation by injectives is mentioned in an Artinian setting in \cite[Example 6.1]{BIKP24}. See also \cite{Ric19} for results concerning non-strong generation by injectives. We recall some notation for building objects in a triangulated category.

\begin{defn} \label{buildingnotation} Let $\mathcal{T}$ be a triangulated category, $S \subseteq \mathcal{T}$ a subcategory.
\begin{enumerate}
\item Let $\thick_0(S) \subseteq \mathcal{T}$ consist of finite sums of shifts of summands of objects of $S$. 

\item Let $\thick_n(S) \subseteq \mathcal{T}$ consist of all summands of objects $X$ that fit into a triangle
\[
X' \to X \to X''
\]
with $X' \in \thick_{n-1}(S)$ and $X'' \in \thick_0(S)$.

\item Let $\thick(S)$ be the union over $n\geq 0$ of all $\thick_n(S)$. 

\end{enumerate}
\end{defn}

\begin{remark} \ \begin{enumerate}
\item We note that $\thick(S)$ is the smallest thick subcategory containing $S$ i.e.\ a subcategory closed under summands, extensions and shifts. 
\item Beware: other authors have used different indexing conventions.
\end{enumerate}
\end{remark}

We let $\mathcal{D}(\mathcal{A})$ denote the derived category of a Grothendieck category $\mathcal{A}$ and $K(\mathcal{A})$ its homotopy category. We recall that Grothendieck categories have enough $K$-injectives and $\mathcal{D}(\mathcal{A})$ is a locally small category with small coproducts and products. Throughout we will write $\sigma_{\leq n}, \sigma_{\geq n}$ for the brutal truncations of a complex and $\tau_{\geq n}, \tau_{\leq n}$ for the smart truncations at an index $n \in \mathbb{Z}$. 

\begin{lemma} \label{injlift}
Let $\mathcal{A}$ be a Grothendieck category, $M \in \mathcal{D}(\mathcal{A})$ and $I \in \mathcal{A}$ an injective object then any map $H^i(M) \to I \in \mathcal{A}$ lifts to a map $\Sigma^i M \to I \in \mathcal{D}(\mathcal{A})$. 
\end{lemma}

\begin{proof}
We can assume that $i = 0$. If $N$ is any complex such that $N^0 = 0$, then $\Hom_{\mathcal{D}(A)}(N,I) = \Hom_{K(A)}(N,I) = 0$ since $I$ is injective. Hence the long exact sequences imply there are isomorphisms
\[
\Hom_{\mathcal{D}(\mathcal{A})}(M, I) \xrightarrow{\sim} \Hom_{\mathcal{D}(\mathcal{A})}(\tau_{\leq 0}M, I) \xleftarrow{\sim} \Hom_{\mathcal{D}(\mathcal{A})}(\tau_{\geq 0} \tau_{\leq 0} M,I) \simeq \Hom_{\mathcal{A}}(H^0(M),I) \qedhere
\]
\end{proof}

If $\mathcal{A}$ is an Abelian category, we let $\Inj(\mathcal{A})$ denote the injectives in $\mathcal{A}$. If $\mathcal{A}$ has enough injectives and $M \in \mathcal{A}$ we let $\idim_{\mathcal{A}} M$ denote the injective dimension of $M$. We let $\gldim(\mathcal{A})$ denote the global dimension of $\mathcal{A}$ i.e.\ the supremum of all the injective dimensions of objects of $\mathcal{A}$.

In our key example, $\mathcal{D}_{\QCoh}(X)$ for a locally Noetherian scheme $X$, we will be working with complexes whose cohomologies are contained in a weak Serre subcategory. Hence we introduce the following notation. Let $\mathcal{S}$ be a weak Serre subcategory of $\mathcal{A}$ such that $\Inj(\mathcal{S}) \subseteq \Inj(\mathcal{A})$. Let $\mathcal{D}_{\Inj(\mathcal{S})}(\mathcal{A}) \subseteq \mathcal{D}_{\mathcal{S}}(\mathcal{A}) \subseteq \mathcal{D}(\mathcal{A})$ denote complexes whose cohomologies are in $\Inj(\mathcal{S})$ and in $\mathcal{S}$ respectively. Note by Lemma \ref{injlift} that any complex with injective cohomologies is quasi-isomorphic to its cohomology.

\begin{lemma} \label{injdimtogen}
Let $\mathcal{S}$ be a weak Serre subcategory of a Grothendieck category $\mathcal{A}$ such that $\mathcal{S}$ has enough injectives and $\Inj(\mathcal{S}) \subseteq \Inj(\mathcal{A})$. Suppose $M \in \mathcal{D}_{\mathcal{S}}(\mathcal{A})$ and $\idim_{\mathcal{S}} H^i(M) \leq n$ for all $i$, then 
\[
M \in \thick_n\left(\mathcal{D}_{\Inj(\mathcal{S})}(\mathcal{A})\right) \subseteq {\mathcal{D}(\mathcal{A})}
\]

\end{lemma}

\begin{proof}
For $n = 0$, the cohomologies of $M$ are all injective in $\mathcal{S}$ and so also in $\mathcal{A}$. By Lemma \ref{injlift}, $M$ is quasi-isomorphic to its cohomology and so is contained in $\mathcal{D}_{\Inj(\mathcal{S})}(\mathcal{A})$. Suppose the result holds for smaller $n$. For every $i$, there exist $I^i \in \Inj(\mathcal{S})$ such that $H^i(M) \hookrightarrow I^i$. By Lemma \ref{injlift}, these lift and combine to a map $M \to \prod \Sigma^{-i} I^i \eqcolon I$ in $\mathcal{D}_
{\mathcal{S}}(\mathcal{A})$. Let $C$ be the cone of this map which is contained in $\mathcal{D}_
{\mathcal{S}}(\mathcal{A})$ since $\mathcal{S}$ is a weak Serre subcategory. This produces short exact sequences
\[
0 \to H^i(M) \to I^i \to H^i(C) \to 0
\]
in $\mathcal{S}$. So $\idim_{\mathcal{S}} H^i(C) \leq n-1$. By induction $C \in \thick_{n-1}(\mathcal{D}_{\Inj(S)}(\mathcal{A}))$ and we are done.
\end{proof}

The following argument is adapted from \cite[Lemma 2.4]{KK05}. 

\begin{lemma}  \label{gentoinjdim} Let $\mathcal{A}$ be a Grothendieck category. If $N \in \mathcal{A}$ and $N \in \thick_n(\mathcal{D}_{\Inj(\mathcal{A})}(\mathcal{A}))$ then $\idim(N) \leq n$.
\end{lemma}

\begin{proof}
We show that $\Ext_{\mathcal{A}}^{n+1}(M,N) = 0$ for any $M \in \mathcal{A}$. Let $\xi \in \Ext_{\mathcal{A}}^{n+1}(M,N)$ be represented by an extension
\[
0 \to N \to E_{n+1} \to E_n \to \dots \to E_1 \to M \to 0
\]
We can split this into $n+1$ short exact sequences
\[
\xi_i: 0 \to X_i \to E_i \to X_{i-1} \to 0 \text{ for }i = 1,\dots,n+1
\]
where $X_0 = M$ and $X_{n+1} = N$. Recall that for any objects $Y,Z \in \mathcal{A}$, there is an isomorphism $\Ext_{\mathcal{A}}^n(Y,Z) \simeq \Hom_{\mathcal{D}(\mathcal{A})}(Y,\Sigma^n Z)$. Under this isomorphism, the composition $(\Sigma^n \xi_{n+1}) \dots (\Sigma \xi_2)\xi_1$ is just $\xi: M \to \Sigma^{n+1} N \in \mathcal{D}(\mathcal{A})$. The short exact sequences $\xi_i$ produce triangles
\[
X_i \to E_i \to X_{i-1} \xrightarrow{\xi_i} \Sigma X_i
\]
For any injective $I$, exactness of $\Hom_{\mathcal{A}}(-,I)$ implies that $\Hom_{\mathcal{D}(\mathcal{A})}(\xi_i,I)$ vanishes. The natural transformations 
\[
\Sigma^{i-1} \xi_i: \Hom_{\mathcal{D}(\mathcal{A})}(\Sigma^i X_i,-) \to \Hom_{\mathcal{D}(\mathcal{A})}(\Sigma^{i-1} X_{i-1},-) 
\]
are between functors which preserve all products. Hence they vanish on $\mathcal{D}_{\Inj(S)}(\mathcal{A})$ which consist of products of shifts of injectives. The coGhost Lemma (the dual of Lemma 2.3 in \cite{KK05}) implies that the composite $\xi$ vanishes on $\thick_n(\mathcal{D}_{\Inj(S)}(\mathcal{A}))$. Since $M \in \thick_n(\mathcal{D}_{\Inj(S)}(\mathcal{A}))$ it follows  that $\xi = 0$. 
\end{proof}

\begin{lemma} \label{idimagrees}
Let $\mathcal{A}$ be a Grothendieck category and $\mathcal{S}$ a weak Serre subcategory of $\mathcal{A}$ which has enough injectives such that $\Inj(\mathcal{S}) = \Inj(\mathcal{A}) \cap \mathcal{S}$. Then for any $M \in \mathcal{S}$, we have $\idim_{\mathcal{S}}(M) = \idim_{\mathcal{A}}(M)$. 
\end{lemma}

\begin{proof}
Any injective resolution of $M$ in $\mathcal{S}$ is also an injective resolution of $M$ in $\mathcal{A}$ since $\Inj(\mathcal{S}) \subseteq \Inj(\mathcal{A})$. Therefore $\idim_{\mathcal{A}}(M) \leq \idim_{\mathcal{S}}(M)$. Suppose $\idim_{\mathcal{A}} M= n < \infty$. Take an injective resolution $M \to I$ in $\mathcal{S}$. As it is also an injective resolution in $\mathcal{A}$, we deduce that the kernel $K$ of the map $I^{n} \to I^{n+1}$ is injective in $\mathcal{A}$. As $\mathcal{S}$ is a weak Serre subcategory $K \in \mathcal{S}$ and so $K \in \Inj(\mathcal{A}) \cap \mathcal{S} = \Inj(\mathcal{S})$ and we are done. 
\end{proof}

\begin{theorem} \label{globdimvinj}
Let $\mathcal{A}$ be a Grothendieck category and $\mathcal{S}$ a weak Serre subcategory of $\mathcal{A}$ which has enough injectives such that $\Inj(\mathcal{S}) = \Inj(\mathcal{A}) \cap \mathcal{S}$. Then
\[
\mathcal{D}_{\mathcal{S}}(\mathcal{A}) = \thick_n(\mathcal{D}_{\Inj(\mathcal{S})}(\mathcal{A})) \text{ if and only if }\gldim(\mathcal{S}) \leq n.
\]
\end{theorem}

\begin{proof}
We always have $\thick_n(\mathcal{D}_{\Inj(\mathcal{S})}(\mathcal{A})) \subseteq \mathcal{D}_{\mathcal{S}}(\mathcal{A})$. Suppose $\mathcal{D}_{\mathcal{S}}(\mathcal{A}) = \thick_n(\mathcal{D}_{\Inj(\mathcal{S})}(\mathcal{A}))$. Let $M \in \mathcal{S}$, then 
\[
M \in \mathcal{D}_{\mathcal{S}}(\mathcal{A}) = \thick_n(\mathcal{D}_{\Inj(\mathcal{S})}(\mathcal{A})) \subseteq \thick_n(\mathcal{D}_{\Inj(\mathcal{A})}(\mathcal{A})).  
\]
So by Lemmas \ref{gentoinjdim} and \ref{idimagrees}, $\idim_{\mathcal{S}}(M) = \idim_{\mathcal{A}}(M) \leq n$. Conversely, if $\gldim(\mathcal{S}) \leq n$, Lemma \ref{injdimtogen} implies $\mathcal{D}_{\mathcal{S}}(\mathcal{A}) = \thick_n(\mathcal{D}_{\Inj(S)}(\mathcal{A}))$.
\end{proof}

\begin{cor} \label{stupidthing}
Let $\mathcal{A}$ be a Grothendieck category and $\mathcal{S}$ a weak Serre subcategory of $\mathcal{A}$ which is also a Grothendieck category such that $\Inj(\mathcal{S}) = \Inj(\mathcal{A}) \cap \mathcal{S}$. If $\mathcal{S}$ has finite global dimension, then the functor
\[
\mathcal{D}(\mathcal{S}) \to \mathcal{D}_{\mathcal{S}}(\mathcal{A})
\]
is an equivalence.
\end{cor}

\begin{proof}
Note the functor is fully faithful on $\mathcal{D}_{\Inj(\mathcal{S})}(\mathcal{S})$. Indeed, this follows since the functor is fully faithful on $\mathcal{S}$, there are no higher Ext's between injectives and any object in $\mathcal{D}_{\Inj(\mathcal{S})}(\mathcal{S})$ is a biproduct of shifts of injectives in $\mathcal{S}$. By Theorem \ref{globdimvinj}, $\mathcal{D}(\mathcal{S}) = \thick(\mathcal{D}_{\Inj(\mathcal{S})}(\mathcal{S}))$ and so the functor is fully faithful on all of $\mathcal{D}(\mathcal{S})$. By Lemma \ref{injdimtogen}, the functor is essentially surjective as its image generates. 
\end{proof}

\begin{remark}
The assumption that $\mathcal{S}$ is also a Grothendieck category in Corollary \ref{stupidthing} can be relaxed but we include it for simplicity. 
\end{remark}

Next we reduce to generation by a single injective. We consider two versions: one using products and one using coproducts. We need some more notation.

\begin{defn} \label{bigbuildingnotation} Let $\mathcal{T}$ be a triangulated category and $S \subseteq \mathcal{T}$.
\begin{enumerate}

\item  Suppose $\mathcal{T}$ has all coproducts. Let $\Add^\Sigma(S)$ denote all coproducts of shifts of summands of objects of $S$ and set $\loc_n(S) \coloneq \thick_n(\Add^\Sigma(S))$.

\item Suppose $\mathcal{T}$ has all products. Let $\Prod^\Sigma(S)$ denote all products of shifts of summands of objects of $\mathcal{S}$ and set $\coloc_n(S) \coloneq \thick_n(\Prod^\Sigma(S))$.

\end{enumerate}
\end{defn}

For a Grothendieck category $\mathcal{A}$, we let $\inj(\mathcal{A})$ denote the indecomposable injectives in $\mathcal{A}$. Recall from Definition \ref{Matliscat}, that $\mathcal{A}$ is Matlis if every injective is a direct sum of indecomposable injectives.

\begin{lemma} \label{loclemma}
Let $\mathcal{S}$ be a weak Serre subcategory of $\mathcal{A}$ such that $\Inj(\mathcal{S}) \subseteq \Inj(\mathcal{A})$. The following hold for any $n \geq 0$.
\begin{enumerate}

\item $\thick_n(\mathcal{D}_{\Inj(\mathcal{S})}(\mathcal{A})) = \coloc_n(\Inj(\mathcal{S})) \subseteq \mathcal{D}_{\mathcal{S}}(\mathcal{A})$ 

\item If $W$ is an injective cogenerator of $\mathcal{A}$, then $\coloc_n(\Inj(\mathcal{A})) = \coloc_n(W) \subseteq \mathcal{D}(\mathcal{A})$.

\item Suppose $\mathcal{A}$ is Matlis. Then there is an $I \in \Inj(\mathcal{A})$ such that 
\[
\coloc_n(\Inj(\mathcal{A})) = \loc_n(\Inj(\mathcal{A})) = \loc_n(\inj(\mathcal{A})) = \loc_n(I).
\]

\end{enumerate}
\end{lemma}

\begin{proof}
By Lemma \ref{injlift}, $\mathcal{D}_{\Inj(\mathcal{S})}(\mathcal{A})$ consists of complexes of injectives with zero differential. It follows that $\mathcal{D}_{\Inj(\mathcal{S})}(\mathcal{A}) = \Prod^\Sigma(\Inj(\mathcal{S}))$ and so (1) holds. 

If $W$ is an injective cogenerator, any injective is a summand of products of $W$. The functor $\mathcal{A} \to \mathcal{D}(\mathcal{A})$ preserves products of injectives as in the proof of Theorem 1.7 in \cite{Sto13}. Hence,
$\Prod^{\Sigma}(W) = \Prod^{\Sigma}(\Inj)$ as required. If $\mathcal{A}$ is Matlis, then direct sums of injectives are injective by Remark \ref{matlisfacts} (4) and every injective is a direct sum of indecomposable injectives by assumption. Let $I$ be the direct sum of the indecomposable injectives. The indecomposable injectives in a Grothendieck category form a set (e.g. \cite[pp. 9]{HPV25}) and so this is well defined. Then we have, 
\[
\Prod^{\Sigma}(\Inj) = \Add^\Sigma(\Inj) = \Add^\Sigma(\inj) = \Add^\Sigma(I)
\]
as required. 
\end{proof}

\begin{remark}
Recall by Remark \ref{locNoethrem}, the condition in (3) is satisfied by locally Noetherian Grothendieck categories and $\QCoh(X)$ where $X$ is a locally Noetherian scheme. 
\end{remark}

\begin{theorem} \label{bigproductglobdim}
Let $\mathcal{A}$ be a Grothendieck category. The following are equivalent
\begin{enumerate}
\item $\gldim \mathcal{A} \leq n$
\item $\mathcal{D}(\mathcal{A}) = \coloc_n(\Inj)$
\item There is an $I \in \Inj(\mathcal{A})$ such that $\mathcal{D}(\mathcal{A}) = \coloc_n(I)$
\end{enumerate}
and in this case $\mathcal{D}(\mathcal{A}) = \loc_n(\Inj)$. If $\mathcal{A}$ is Matlis, then the above are equivalent to the following.
\begin{enumerate}

\item[(4)] $\mathcal{D}(\mathcal{A}) = \loc_n(\inj)$

\item[(5)] There is an $I \in \Inj(\mathcal{A})$ such that $\mathcal{D}(\mathcal{A}) = \loc_n(I)$. 
\end{enumerate}
\end{theorem}

\begin{proof}
The equivalence of (1) and (2) is Theorem \ref{globdimvinj} and Lemma \ref{loclemma} (1). Also (2) implies (3) by Lemma \ref{loclemma} (2). And (3) implies (2) since $\coloc_n(I) \subseteq \coloc_n(\Inj)$. Note that $\Prod^\Sigma(\Inj) \subseteq \Add^{\Sigma}(\Inj)$ since the former consists of complexes of injectives with zero differential and the latter consists of complexes of summands of sums of injectives with zero differential. Hence there is an inclusion $\coloc_n(\Inj(\mathcal{A})) \subseteq \loc_n(\Inj(\mathcal{A}))$ from which the result follows. The proof of the second part is similar to the above using Lemma \ref{loclemma} (3).
\end{proof}

If $R$ is a (not necessarily commutative) ring, we let $\mathcal{D}(R)$ denote the unbounded derived category of all right $R$-modules.

\begin{theorem} The following are equivalent for any ring $R$. 
\begin{enumerate}
\item The right global dimension of $R$ is at most $n$.
\item $\mathcal{D}(R) = \loc_n(R)$.
\item $\mathcal{D}(R) = \coloc_n(\Inj(R))$
\item There is an injective $R$-module $I$ such that $\mathcal{D}(R) = \coloc_n(I)$.
\end{enumerate}
If $R$ is right Noetherian then the above are equivalent to 
\begin{enumerate}
\item $\mathcal{D}(R) = \loc_n(\inj(R))$ 
\item There is an injective $I$ such that  $\mathcal{D}(R) = \loc_n(I)$.
\end{enumerate}

\end{theorem}

\begin{proof}
The equivalence between (1) and (2) can be proved using big versions of the argument in \cite[Lemma 2.4]{KK05}. The proof is dual to that of Theorem \ref{globdimvinj} with injectives replaced with projectives.  The rest follow from Theorem \ref{bigproductglobdim} since if $R$ is right Noetherian, then the category of right $R$-modules is locally Noetherian and in particular Matlis. 
\end{proof}

For a scheme $X$ we recall that $\QCoh(X)$ is a weak Serre subcategory of $\Mod \mathcal{O}_X$ and let $\mathcal{D}_{\QCoh}(X) \coloneq \mathcal{D}_{\QCoh(X)}(\Mod \mathcal{O}_X)$. We let $\dim X$ denote the Krull dimension of $X$.

\begin{theorem} \label{Bigstronggenregular}
Suppose $X$ is a locally Noetherian scheme. Then the following are equivalent.

\begin{enumerate}
\item $X$ is regular and $\dim X \leq n$.
\item $\mathcal{D}(\QCoh(X)) = \loc_n(\inj(\QCoh(X)))$.
\item There is an $I \in \Inj(\QCoh(X)) $ such that $\mathcal{D}(\QCoh(X)) = \loc_n(I)$.
\item $\mathcal{D}_{\QCoh}(X) = \loc_n(\inj(\QCoh(X))$.
\item There is an $I \in \Inj(\QCoh(X))$ such that $\mathcal{D}_{\QCoh}(X) = \loc_n(I)$.
\end{enumerate}
and if these conditions hold the natural map is an equivalence $\mathcal{D}(\QCoh(X)) \xrightarrow{\sim} \mathcal{D}_{\QCoh}(X)$.

\end{theorem}

\begin{proof}
By Theorem \ref{matlisthy}, $\QCoh(X)$ is Matlis and so (1), (2) and (3) are equivalent by Theorems \ref{bigproductglobdim} and \ref{ABSscheme}. Since $\QCoh(X) \subseteq \Mod \mathcal{O}_X$ is closed under coproducts, similar arguments to the proof of Lemma \ref{loclemma} apply and 
\[
\thick_n(\mathcal{D}_{\Inj(\QCoh(X))}) = \coloc_n(\Inj(\QCoh(X)) = \loc_n(\Inj(\QCoh(X)) \subseteq \mathcal{D}_{\QCoh}(X)
\]
And similarly the above coincides with $\loc_n(\inj(\QCoh(X)) = \loc_n(I)$ for an injective $I$ in $\QCoh(X)$. Hence (1) is equivalent to (4) and (5) by Theorem \ref{globdimvinj}. If any of these conditions hold then there is an equivalence as claimed by Corollary \ref{stupidthing}.
\end{proof}

\begin{remark}
The last part of Theorem \ref{Bigstronggenregular} was conjectured by Hartshorne in \cite[pp.\ 136]{Har66}.
The comparison map $\mathcal{D}(\QCoh(X)) \to \mathcal{D}_{\QCoh}(X)$ is known to be an equivalence for Noetherian schemes (\cite[Tag 09T4]{stacks-project}), and for quasicompact semiseparated schemes (\cite[Corollary 5.5]{BN93}). 
\end{remark}

\begin{remark}
We note that if $X$ is any scheme, Theorem \ref{bigproductglobdim} still applies to relate the global dimension of $\QCoh(X)$ to strong generation of $\mathcal{D}(\QCoh(X))$ by injectives. However in this general setting, we cannot relate global dimension to regularity and have no control over the injectives.
\end{remark}

\section{Application: Global Dimension under Derived Equivalence}
\label{derivedequivsec}

As an application of the results of Section \ref{bigstronggensec} we consider when global dimension is preserved under derived equivalences between Grothendieck categories.

\begin{remark}
Throughout this section we assume that $\mathcal{A}$ and $\mathcal{B}$ are Grothendieck categories and there is a triangle equivalence $G: \mathcal{D}(\mathcal{A}) \xrightarrow{\sim} \mathcal{D}(\mathcal{B})$ with inverse $F$.
\end{remark}

For $m,n \in \mathbb{Z}$, we let $\mathcal{D}^{[m,n]}(\mathcal{A})$ denote complexes whose cohomology vanishes outside the interval $[m,n]$ and let $\mathcal{D}^b(\mathcal{A})$ denote complexes with bounded cohomology.  We recall that a Grothendieck category satisfies Roos's axiom $\Ab4^\ast$-$l$ for some $l \geq 0$ if the right derived functors of all infinite products vanish in degrees greater than $l$. So $\Ab4^\ast$-$0$ is the $\Ab4^\ast$ axiom stating that products are exact. See \cite{Roo06} for more detail.

\begin{lemma} \label{bigmoritalemma0}
Suppose that $G$ restricts to $\mathcal{D}^b(\mathcal{A}) \to \mathcal{D}^b(\mathcal{B})$ and that either
\begin{enumerate}
\item $\mathcal{A}$ is Matlis or
\item $\mathcal{B}$ satisfies $Ab^\ast$-$l$ for some $l \geq 0$
\end{enumerate}
then there exist $m,n \in \mathbb{Z}$ such that $G(\Inj(\mathcal{A})) \subseteq \mathcal{D}(\mathcal{B})^{[m,n]}$. 
\end{lemma}

\begin{proof}
First assume (1). Let $I \in \Inj(\mathcal{A})$ be the direct sum of all the indecomposable injectives (which is well-defined as the indecomposable injectives form a set). As we have assumed that $G$ sends bounded complexes to bounded complexes, there are $m,n$ such that $G(I) \in \mathcal{D}(\mathcal{B})^{[m,n]}$. As cohomology commutes with coproducts it follows that $G$ sends any injective to $\mathcal{D}(\mathcal{B})^{[m,n]}$.

Now suppose (2) holds and let $W$ be an injective cogenerator of $\mathcal{A}$. Then there are $m,n$ such that $G(W) \in \mathcal{D}(\mathcal{B})^{[m,n]}$. Any injective is a direct summand of a product of copies of $W$ taken in $\mathcal{A}$. The functor $\mathcal{A} \to \mathcal{D}(\mathcal{A})$ preserves products of injectives as in the proof of Theorem 1.7 in \cite{Sto13}. As $G$ preserves products, we see that $G$ sends any injective to a summand of products of $G(W)$. As $\mathcal{B}$ is $\Ab4^\ast$-$l$, we have that $H^\ast(\prod G(W))$ is non-zero only in degrees $[m,n + l]$  as required.
\end{proof}

If $S \subseteq \mathcal{T}$ is a collection of objects in a triangulated category, we let $\thick^{\geq 0}(S)$ denote the variant of $\thick(S)$ from Definition \ref{buildingnotation} where only non-negative shifts are allowed.  

\begin{remark} \label{milnor}
As $\mathcal{A}$ is a Grothendieck category, there is a short exact sequence of chain complexes over $\mathcal{A}$
\[
0 \to \bigoplus_{i \leq 0} \sigma_{\geq i} X \to \bigoplus_{i \leq 0} \sigma_{\geq i} X \to X \to 0 
\]
for any chain complex $X$ over $\mathcal{A}$. This produces the so called Milnor triangle in $\mathcal{D}(\mathcal{A})$ and we write that $X \simeq \hocolim(\sigma_{\geq i} X)$. 
\end{remark}

\begin{lemma} \label{bigmoritalemma1}
Suppose $\gldim(\mathcal{A}) = d < \infty$ and there is some $n \in \mathbb{Z}$ with $G(\Inj(\mathcal{A})) \subseteq \mathcal{D}(\mathcal{B})^{\leq n}$. Then $F$ sends $\mathcal{D}(\mathcal{B})^{\geq 0}$ to $\mathcal{D}(\mathcal{A})^{\geq -d-n}$. 
\end{lemma}

\begin{proof}
First let $X \in \mathcal{D}^b(\mathcal{A})^{\leq -d}$ and let $X \xrightarrow{\sim} J$ denote an injective resolution. We can assume that $J^i = 0$ when $i$ is less than the index of the smallest non-zero cohomology of $X$. Since $\gldim(\mathcal{A}) = d$ we can assume $J^i = 0$ for $i > 0$. Using the brutal truncations we see that $J \in \thick^{\geq 0}(\Inj(\mathcal{A}))$. As $\mathcal{D}(\mathcal{B})^{\leq n}$ is closed under positive shifts and extensions we deduce that $G(X) \simeq G(J) \in \mathcal{D}(\mathcal{B})^{\leq n}$. Now suppose $Y \in \mathcal{D}(\mathcal{B})^{\geq n}$, then 
\begin{equation} \label{eq1}
\Hom_{\mathcal{D}(\mathcal{A})}(X, \Sigma^{-1} F(Y)) \simeq \Hom_{\mathcal{D}(\mathcal{B})}(G(X),\Sigma^{-1} Y) = 0
\end{equation}
since $G(X) \in \mathcal{D}(\mathcal{B})^{\leq n}$ and $Y \in \mathcal{D}(\mathcal{B})^{\geq n}$ using the standard $t$-structure. If $X' \in \mathcal{D}(\mathcal{A})^{\leq -d}$, then $X' = \hocolim(\sigma_{\geq -i} X')$. Since $\sigma_{\geq -i} X' \in \mathcal{D}^b(\mathcal{A})^{\leq -d}$, (\ref{eq1}) implies we have $\Hom_{\mathcal{D}(\mathcal{A})}(\sigma_{\geq -i} X', \Sigma^{-1} F(Y)) =0$. The Milnor triangle for the homotopy colimit then implies that  
\[
\Hom_{\mathcal{D}(\mathcal{A})}(X',\Sigma^{-1}F(Y)) = 0
\]
Therefore $\Sigma^{-1}F(Y)$ is contained in the right perpendicular to $\mathcal{D}(\mathcal{A})^{\leq - d}$ and so $F(Y) \in \mathcal{D}(\mathcal{A})^{ \geq -d}$. By shifting the result follows.
\end{proof}

\begin{lemma}\label{bigmoritalemma2}
Suppose $\gldim(\mathcal{A}) = d < \infty$ and there exist $m,n$ such that $G(\Inj(\mathcal{A})) \subseteq \mathcal{D}(\mathcal{B})^{[m,n]}$. If $I \in \Inj(\mathcal{A})$ then $G(I) \in \thick_{d+n-m+1}(\Inj(\mathcal{B}))$.
\end{lemma}

\begin{proof}
Let $X = G(I) \in \mathcal{D}(\mathcal{B})^{[m,n]}$ and let $X \to J$ denote an injective resolution of $X$. Then we can assume $J^i = 0$ for $i < m$. Let $K$ be the kernel of the differential $J^{d+n+1} \to J^{d+n+2}$. Then 
\[
\Hom_{\mathcal{D}(\mathcal{B})}(K, \Sigma^{d+n+1} J) \simeq \Hom_{\mathcal{D}(\mathcal{A})}(F(K), \Sigma^{d+n+1} I) \simeq \Hom_{K(\mathcal{A})}(F(K), \Sigma^{d+n+1} I)
\]
By Lemma \ref{bigmoritalemma1}, $F(K) \in \mathcal{D}(\mathcal{A})^{\geq -d-n}$ and we can assume $F(K)^i = 0$ for $i<-d-n$. Therefore $\Hom_{K(\mathcal{A})}(F(K), \Sigma^{d+n+1} I) = 0$. It follows that the natural map $K \to \Sigma^{d+n+1} J$ is null homotopic. So there is a map $s: K \to J^{d+n}$ splitting the differential $J^{d+n} \to K$. Hence $K$ is a summand of $J^{d+n}$. Therefore $X \simeq \tau_{\leq d + n + 1} J$ and $\tau_{\leq d + n + 1} J$ is a complex of injectives which is non-zero only in degrees $[m, d + n +1]$. Hence the brutal truncation triangles imply that $X \in \thick_{d+n-m+1}(\Inj(\mathcal{B}))$. 
\end{proof}

\begin{lemma} \label{bigmoritalemma3}
Suppose $\gldim(\mathcal{A}) = d < \infty$ and there exist $m,n \in \mathbb{Z}$ with $G(\Inj(\mathcal{A})) \subseteq \mathcal{D}(\mathcal{B})^{[m,n]}$. Then $\gldim(\mathcal{B}) \leq (d+1)(d+n-m+1)$
\end{lemma}

\begin{proof}
By Lemma \ref{bigmoritalemma2}, $G(\Inj(\mathcal{A})) \subseteq \thick_{d + n -m+ 1}(\Inj(\mathcal{B}))$. By Theorem \ref{bigproductglobdim}, $\mathcal{D}(\mathcal{A}) = \coloc_d(\Inj(\mathcal{A}))$. Applying $G$ gives that
\begin{align*}
\mathcal{D}(\mathcal{B}) &= G(\coloc_d(\Inj(\mathcal{A}))) \\ 
&= \coloc_d(G(\Inj(\mathcal{A}))) \\
& = \coloc_d(\thick_{d+n-m+1}(\Inj(\mathcal{B}))) \\
& = \coloc_{(d+1)(d+n-m+1)}(\Inj(\mathcal{B}))
\end{align*}
Therefore $\mathcal{B}$ has global dimension at most $(d+1)(d+n-m+1)$ by Theorem \ref{bigproductglobdim}. 
\end{proof}

\begin{theorem} \label{bigmoritaglobdim}
Let $\mathcal{A}$ and $\mathcal{B}$ denote Grothendieck categories such that either of the conditions below hold
\begin{enumerate}
\item $\mathcal{A}$ is Matlis or

\item $\mathcal{B}$ satisfies $\Ab 4^\ast$-$l$ for some $l \geq 0$.

\end{enumerate}
Suppose there is an equivalence $G: \mathcal{D}(\mathcal{A}) \xrightarrow{\sim} \mathcal{D}(\mathcal{B})$ such that $G$ sends bounded complexes to bounded complexes. If $\gldim(\mathcal{A}) < \infty$, then $\gldim(\mathcal{B}) < \infty$. 
\end{theorem}

\begin{proof}
If any of the conditions are satisfied, Lemma \ref{bigmoritalemma0} applies to $G$. Therefore by Lemma \ref{bigmoritalemma3} the claim follows. 
\end{proof}

\begin{remark}
Note that in Theorem \ref{bigmoritaglobdim} we haven't insisted that $G$ restricts to an equivalence $G: \mathcal{D}^b(\mathcal{A}) \xrightarrow{\sim} \mathcal{D}^b(\mathcal{B})$ just that it restricts. Since $\mathcal{A}$ is assumed to have finite global dimension and $G$ is usually a derived functor which can be computed using resolutions, this condition is easy to check in practice.
\end{remark}

\begin{theorem} \label{bigmoritascheme}
Let $X,Y$ denote Noetherian separated schemes and suppose that $G: \mathcal{D}_{\QCoh}(X) \xrightarrow{\sim} \mathcal{D}_{\QCoh}(Y)$ is an equivalence. Then, there exist $m,n$ such that $G(E(x)) \in \mathcal{D}_{\QCoh}^{[m,n]}(Y)$ for every $x \in X$. If $X$ is regular of dimension $d$, then $Y$ is regular and of dimension at most $(d+1)(d+1 + n-m)$. 
\end{theorem}

\begin{proof}
As $X$ is quasicompact and separated, $\mathcal{D}_{\QCoh}(X) \simeq \mathcal{D}(\QCoh(X))$ is compactly generated (\cite[Proposition 2.5]{NeeGroth}) by the perfect complexes $\mathcal{D}^{\perf}(X)$. It follows from \cite[Proposition 7.0.2]{BNP13} that $M \in \mathcal{D}^b(\QCoh(X))$ if and only if $\Hom_{\mathcal{D}(\QCoh(X))}(P,\Sigma^i M)$ is bounded for all $P \in \mathcal{D}^{\perf}(X)$. The same results hold for $Y$ and so the equivalence restricts to $\mathcal{D}^b(\QCoh(X)) \simeq \mathcal{D}^b(\QCoh(Y))$ as it restricts to compacts $\mathcal{D}^{\perf}(X) \simeq \mathcal{D}^{\perf}(Y)$. By Theorem \ref{matlisthy}, every injective in $\QCoh(X)$ is a direct sum of the indecomposable injectives $E(x)$. Since $G(\bigoplus E(k(x)))$ is bounded, there exist $m,n$ as claimed. As in Lemma \ref{bigmoritalemma0}, this implies $G(\Inj(X)) \subseteq \mathcal{D}^{[m,n]}(\QCoh(Y))$. Hence the claim follows from Lemma \ref{bigmoritalemma3}.
\end{proof}

\begin{remark}
If $X,Y$ are as in Theorem \ref{bigmoritascheme}, then the equivalence restricts to compacts  $\mathcal{D}^{\perf}(X) \simeq \mathcal{D}^{\perf}(Y)$. By \cite[Theorem 0.5]{Nee21b}, it follows that $X$ is regular and finite-dimensional if and only if $Y$ is. However, there is no such bound on the dimensions as given in Theorem \ref{bigmoritascheme}. Note also the Noetherian hypothesis is missing from \cite[Theorem 0.5]{Nee21b} but it is necessary as noted in \cite{Ste25}.
\end{remark}

\begin{remark}
If $X,Y$ are quasicompact and separated but not necessarily Noetherian, then by \cite[Theorem 1.4]{Pos26}, $\QCoh(X)$ is $\Ab4^\ast$-$n$ where $n$ is one less than the number of affines needed to cover $X$. So using similar arguments to Theorem \ref{bigmoritascheme} but with Lemma \ref{bigmoritalemma0} 2), we see that if $\mathcal{D}_{\QCoh}(X) \simeq \mathcal{D}_{\QCoh}(Y)$, then $\gldim(\QCoh(X)) < \infty$ if and only if $\gldim(\QCoh(Y)) < \infty$. However without Auslander-Buchsbaum-Serre, the geometric meaning of this is somewhat unclear. 
\end{remark}

\begin{remark}
Recall that a tilting object is a compact generator $T \in \mathcal{D}(\mathcal{A})$ such that $\Hom_{\mathcal{D}(\mathcal{A})}(T,\Sigma^i T) = 0$ for $i \neq 0$. Let $\Lambda \coloneq \End_{\mathcal{D}(\mathcal{A})}(T)$, then there is an equivalence
\[
G\coloneq \RHom_{\mathcal{A}}(T,-): \mathcal{D}(\mathcal{A}) \xrightarrow{\sim} \mathcal{D}( \Lambda)
\]
Here $\RHom_{\mathcal{A}}(T,-)$ denotes the right derived functor which exists since $\mathcal{A}$ has enough $K$-injectives. Again we let $F$ denote the inverse of $G$ and fix this notation for the remainder of the section.
\end{remark}

Our arguments almost recover Theorem 2 of \cite{BK20} except that our bound on the global dimension of the endomorphism ring of the tilting object is worse. However we can also prove the converse using a result of \cite{Pos26}.

\begin{lemma} \label{boundrestrict}
If $\gldim \Lambda < \infty$, then $G$ restricts to an equivalence 
\[
\mathcal{D}^b(\mathcal{A}) \xrightarrow{\sim} \mathcal{D}^b(\Lambda)
\]
\end{lemma}

\begin{proof}
We note that $T \in \mathcal{D}^b(\mathcal{A})$ by e.g.\ \cite[Lemma 7]{BK20}. As $\Lambda$ has finite global dimension, every complex in $\mathcal{D}^b(\Lambda)$ is quasi-isomorphic to a bounded complex of projectives. Since $F(\Lambda) \simeq T \in \mathcal{D}^b(\mathcal{A})$, it follows that $F$ restricts to $\mathcal{D}^b(\Lambda) \to \mathcal{D}^b(\mathcal{A})$ (c.f. the proof of \cite[Lemma 8]{BK20}). So it remains to show that $G$ restricts. Let $M \in \mathcal{D}^{\geq n}(\mathcal{A})$ for some $n$. As some shift of $T$ is contained in $\mathcal{D}(\mathcal{A})^{< n}$, there is an $n'$ such that 
\[
0 = \Hom_{\mathcal{D}(\mathcal{A})}(T, \Sigma^i M) = H^i(G(M))
\]
for $i < n'$. Hence $G(M) \in \mathcal{D}(\Lambda)^{\geq n'}$. 

Suppose $N \in \mathcal{D}^{\geq n}(\Lambda)$ for some $n$. Then $N \simeq \hocolim(\tau_{\leq m} N)$, the smart truncations of $N$ (by e.g.\ \cite[Lemma 4]{BK20}). Hence $F(N) \simeq \hocolim(F(\tau_{\leq m} N))$. Since each $\tau_{\leq m} N \in \mathcal{D}(\Lambda)^{[n,m]}$, there is an $n'$ such that $F(\tau_{\leq m} N) \in \mathcal{D}^{\geq n'}(\mathcal{A})$ for every $m$. Indeed this follows from the proof above that $F$ restricts to bounded complexes. The Milnor triangle for the homotopy colimit then implies that $H^i(F(M)) = 0$ for $i < n'-1$. Therefore $F(\mathcal{D}(\Lambda)^{\geq n}) \subseteq \mathcal{D}(\mathcal{A})^{\geq n' -1}$. Since $\mathcal{D}(\Lambda)^{> n}$ is the right perpendicular to $\mathcal{D}(\Lambda)^{\leq n}$, it follows that $G(\mathcal{D}(\mathcal{A})^{< n'-1}) \subseteq \mathcal{D}(\Lambda)^{< n}$. Therefore $G$ preserves bounded above and bounded below complexes and so it restricts to $\mathcal{D}^b(\mathcal{A}) \xrightarrow{\sim} \mathcal{D}^b(\Lambda)$. 
\end{proof}

\begin{theorem} \label{tiltingglobdim}
Let $\mathcal{A}$ denote a Grothendieck category, $T \in \mathcal{D}(\mathcal{A})$ a tilting object. Then $\gldim \mathcal{A} < \infty$ if and only if $\gldim \End_\mathcal{A}(T) < \infty$. 
\end{theorem}

\begin{proof}
If $\mathcal{A}$ has finite global dimension, then any $M \in \mathcal{D}^b(\mathcal{A})$ has a finite injective resolution $M \xrightarrow{\sim} I$. Then $G(M) \simeq \Hom_{\mathcal{A}}(T,I)$ is bounded since $T$ and $I$ are. Since products in $\Mod \Lambda$ are exact, condition 2) of Lemma \ref{bigmoritalemma0} is met. So by Theorem \ref{bigmoritaglobdim}, the global dimension of $\Lambda$ is finite. 

Conversely, suppose $\gldim(\Lambda) < \infty$. By Lemma \ref{boundrestrict}, $G$ restricts to an equivalence $\mathcal{D}^b(\mathcal{A}) \xrightarrow{\sim} \mathcal{D}^b(\Lambda)$. By Theorem 4.7 of \cite{Pos26}, and since products in $\Mod \Lambda$ are exact, there is an $n$ such that $\mathcal{A}$ is $\Ab 4^\ast$-$n$. Then the inverse of $G$ meets the conditions of Theorem \ref{bigmoritaglobdim}. 
\end{proof}

\section{Small Strong Generation by Injectives and Immaculate complexes}

\label{immacsec}

In Section \ref{bigstronggensec}, we showed that big strong generation by injectives of the derived category can be used to detect global dimension. We now consider a small version. The results here are morally dual to the fact that strong generation of the perfects is a measure of regularity (\cite[Theorem 3.14]{VB02}, \cite[Theorem 0.5]{Nee21b}) with the caveat that we now work with a set of generators instead of a single generator. For a Grothendieck category $\mathcal{A}$ we recall that $\inj = \inj(\mathcal{A}) \subseteq \mathcal{A}$ denotes the set of indecomposable injectives.

\begin{defn} \label{immacdefn}Let $\mathcal{A}$ be a Grothendieck category.
\begin{enumerate}
\item A complex in $\mathcal{D}(\mathcal{A})$ is immaculate if it is quasi-isomorphic to a bounded complex of finite sums of indecomposable injectives. Let $\mathcal{D}^{\immac}(\mathcal{A})$ denote the  subcategory of immaculate complexes.

\item A complex is $\inj$-pseudocoherent if it is quasi-isomorphic to a bounded below complex of finite sums of indecomposable injectives. Let $\mathcal{D}_{\ipscoh}(\mathcal{A})$ denote the $\inj$-pseudocoherent complexes and, $\mathcal{D}_{\ipscoh}^b(\mathcal{A})$ the $\inj$-pseudocoherent complexes with bounded cohomology.

\item An object of $\mathcal{A}$ is $\inj$-pseudocoherent if it is $\inj$-pseudocoherent as an object of $\mathcal{D}(\mathcal{A})$. 

\item An object  $M \in \mathcal{A}$ is $\inj$-finitely presented if there is an exact sequence
\[
0 \to M \to I \to J 
\]
with $I,J$ finite sums of indecomposable injectives. We let $\ifp(\mathcal{A}) \subseteq \mathcal{A}$ denote the $\inj$-finitely presented objects. We let $\mathcal{D}^b_{\ifp}(\mathcal{A}) \subseteq \mathcal{D}^+_{\ifp}(\mathcal{A}) \subseteq \mathcal{D}(\mathcal{A})$ denote the complexes whose cohomologies are all $\inj$-finitely presented and bounded or bounded below respectively.
\end{enumerate}
\end{defn}

We prove the following in Appendix \ref{thickimmac}.

\begin{prop} \label{thingsarethick}
If $\mathcal{A}$ is a Grothendieck category, then $\mathcal{D}_{\ipscoh}(\mathcal{A})$ and $\mathcal{D}^{\immac}(\mathcal{A})$ are both thick and $\mathcal{D}^{\immac}(\mathcal{A}) = \thick(\inj)$. 
\end{prop}

\begin{remark}
We note that $\mathcal{D}^{\immac}(\mathcal{A})$ is essentially small since the isomorphism classes of indecomposable injectives in a Grothendieck category form a set. 
\end{remark}

We now relate injective dimension to generation of the immaculates.

\begin{lemma} \label{vanishinj}
If $f\colon X \to Y \in \mathcal{D}(\mathcal{A})$ vanishes on cohomology then for any injective object $I$,  $\Hom_{\mathcal{D}(\mathcal{A})}(f,I) = 0$.
\end{lemma}

\begin{proof}
The proof of Lemma \ref{injlift} shows that there is a commuting diagram
\[
\begin{tikzcd}[column sep = huge]
\Hom_{\mathcal{D}(\mathcal{A})}(Y,I) \arrow[r,"{\Hom(f,I)}"] \arrow[d,"H^0"] &  \Hom_{\mathcal{D}(\mathcal{A})}(X,I) \arrow[d,"H^0"] \\
\Hom_{\mathcal{A}}(H^0(Y),I) \arrow[r,"{\Hom(H^0(f),I)}"] & \Hom_{\mathcal{A}}(H^0(X),I)
\end{tikzcd}
\]
where the vertical maps are isomorphisms.
\end{proof}

The proof of Proposition \ref{smallinjgen} is dual to \cite[Proposition A.1]{Kra24b}.

\begin{prop} \label{smallinjgen} Let $\mathcal{A}$ be a Grothendieck category.
\begin{enumerate}
\item If $\mathcal{D}^{\immac}(\mathcal{A}) = \mathcal{D}^b_{\ipscoh}(\mathcal{A})$, then every $\inj$-pseudocoherent object has finite injective dimension.

\item If there is an $n \geq 0$ such that $\mathcal{D}^{\immac}(\mathcal{A}) = \thick_n(\inj)$, then $\mathcal{D}^{\immac}(\mathcal{A}) = \mathcal{D}^b_{\ipscoh}(\mathcal{A})$ and every $\inj$-pseudocoherent object of $\mathcal{A}$ has injective dimension at most $n$. 
\end{enumerate}
\end{prop}

\begin{proof}
For (1), if $M \in \mathcal{A}$ is $\inj$-pseudocoherent, then $M \in \mathcal{D}^b_{\ipscoh}(\mathcal{A}) = \mathcal{D}^{\immac}(\mathcal{A})$ and so $M \in \thick_n(\inj) \subseteq \coloc_n(\Inj)$ for some $n$. Hence $M$ has finite injective dimension by Lemma \ref{gentoinjdim}.

For (2), suppose $\mathcal{D}^{\immac}(\mathcal{A}) = \thick_n(\inj)$. Let $M \in \mathcal{D}^b_{\ipscoh}(\mathcal{A})$ which we may assume is a bounded below complex of finite sums of indecomposable injectives with bounded cohomology. So  $H^i(M) = 0$ for $i \gg 0$. Then by taking the brutal truncations we may assume that $M$ is bounded below at zero and has cohomology only in degree zero. Consider the maps
\[
\sigma_{[n+1,2n+3]}M \xrightarrow{\phi_0} \sigma_{[n,2n+2]} M \to \dots \to  \sigma_{[1,n+3]} M \xrightarrow{\phi_n} \sigma_{[0,n+2]} M
\]
between the brutal truncations where the subscripts denote the intervals which are non-zero. Each of the $\sigma_{[i,i+n+2]} M$ has cohomology only in degrees $i$ and $i+n+2$. Therefore each of the maps $\phi_i$ vanish on cohomology. By Lemma \ref{vanishinj}, it follows that for any injective module $J$, the induced map $\Hom_{\mathcal{D}(\mathcal{A})}(\phi_i,J) = 0$. Hence $\Hom_{\mathcal{D}(\mathcal{A})}(\phi_i,-)$ vanishes on $\thick_0(\inj)$ so by the coGhost Lemma, the composite $\phi= \phi_n \cdots \phi_0$ vanishes on $\thick_{n}(\inj) = \mathcal{D}^{\immac
}(\mathcal{A})$. Since $\phi$ is a morphism between objects in $\mathcal{D}^{\immac
}(\mathcal{A})$, this implies that $\phi = 0$. Since the complexes are $K$-injective, we deduce that $\phi$ is null homotopic. So there is a map $s: \im(d^n) \to I^{n+1}$ such that $d^n s^{n+1} = \id_{\im(d^n)}$. Therefore $\im(d^n)$ is a summand of $M^n$ and thus it is a finite sum of indecomposable injectives. Therefore the smart truncation $\tau_{\leq n+1} M$, which is quasi-isomorphic to $M$, is a bounded complex of finite sums of indecomposable injectives. Therefore $M$ is immaculate and so $\mathcal{D}^b_{\ipscoh}(\mathcal{A}) = \mathcal{D}^{\immac}(\mathcal{A})$. If $M \in \mathcal{A}$ is any $\inj$-pseudocoherent object then as in (1) $M \in \mathcal{D}^b_{\ipscoh}(\mathcal{A}) = \mathcal{D}^{\immac}(\mathcal{A}) = \thick_n(\inj)$ and $M$ has injective dimension at most $n$. 
\end{proof}

In order to say more, we need the following coherence assumption on $\mathcal{A}$. 

\begin{defn} \label{injcohdefn}
A Grothendieck category $\mathcal{A}$ is inj-coherent if for any map $f: I \to J$ between finite sums of indecomposable injectives, the injective hull $E(\coker(f))$ is a finite sum of indecomposable injectives.
\end{defn}

We note the following characterisation of $\inj$-coherence. The proof is similar to the analogous and well known result for coherent rings. 
Recall that the derived category of any Grothendieck category $\mathcal{D}(\mathcal{A})$ admits the standard $t$-structure (see \cite{BBD82}) whose aisle and coaisle consist of complexes with vanishing cohomology in positive and negative degrees respectively.

\begin{prop} \label{injcohequiv}
The following are equivalent.
\begin{enumerate}
\item $\mathcal{A}$ is $\inj$-coherent.
\item Any $M \in \ifp(\mathcal{A})$ admits an injective resolution $I$ with each $I^i$ a finite sum of indecomposable injectives for $i \geq 0$.

\item The standard $t$-structure restricts to $\mathcal{D}_{\ipscoh}(\mathcal{A})$ (or to $\mathcal{D}^b_{\ipscoh}(\mathcal{A})$).

\item We have $\mathcal{D}^b_{\ipscoh}(\mathcal{A}) = \mathcal{D}^b_{\ifp}(\mathcal{A})$ (or $\mathcal{D}_{\ipscoh}(\mathcal{A}) = \mathcal{D}_{\ifp}^+(\mathcal{A})$).
 
\end{enumerate}
\end{prop}

\begin{remark} 
In a hereditary (global dimension 1) Grothendieck category $\mathcal{A}$, the quotient of any injective is injective and it follows that $\mathcal{A}$ is $\inj$-coherent. Also since every complex is quasi-isomorphic to its cohomology, 
\[
\mathcal{D}^b_{\ipscoh}(\mathcal{A}) = \mathcal{D}^b_{\ifp}(\mathcal{A}) = \mathcal{D}^{\immac}(\mathcal{A}) = \thick_1(\inj).
\]
\end{remark}

\begin{remark}
One could also define an ``$\inj$'' version of being Noetherian in an obvious way. However, this property rarely holds. For example, the indecomposable injective $\mathbb{Q}$ over $\mathbb{Z}$ has a quotient $\mathbb{Q}/\mathbb{Z}$ whose injective hull (which is itself) is not a finite direct sum of indecomposable injectives. 
\end{remark}

\begin{remark}
The main result of Section \ref{immaculateschemes} is that $\QCoh(X)$ is $\inj$-coherent where $X$ is a locally Noetherian scheme. It would be interesting to consider this property for other naturally occurring Grothendieck categories. If $A$ is a right Artinian ring, then the indecomposable injectives are the injective hulls of the simple modules. It follows that $\Mod A$ is $\inj$-coherent. We expect that there are right Noetherian rings whose module categories are not $\inj$-coherent but we do not know any examples. 
\end{remark}

\begin{lemma} \label{smallglobdimtostronggen}
Suppose $\mathcal{A}$ is $\inj$-coherent. If $X \in \mathcal{D}^b_{\ipscoh}(\mathcal{A}) = \mathcal{D}^b_{\ifp}(\mathcal{A})$ and $\idim H^i(X) \leq n$ for all $i$, then $X \in \thick_n(\inj)$.
\end{lemma}

\begin{proof}
The proof is similar to Lemma \ref{injdimtogen}. If $n = 0$, then Lemma \ref{injlift} implies that $X \simeq \bigoplus \Sigma^i H^i(X)$. Each $H^i(X)$ is an injective $\inj$-finitely presented module and thus is a finite sum of indecomposable injectives. Hence $X \in \thick_0(\inj)$. Suppose the result holds for smaller $n$. Since $\mathcal{A}$ is $\inj$-coherent, each $H^i(X)$ is $\inj$-finitely presented by Proposition \ref{injcohequiv}. So there are monomorphisms $H^i(X) \hookrightarrow I^i$ such that $I^i$ is a finite sum of indecomposable injectives. Let $C$ denote the cone of the map $X \to \prod \Sigma^{-i} I^i$ defined using Lemma \ref{injlift}. Then $C \in \mathcal{D}^b_{\ipscoh}(\mathcal{A})$ and $\idim H^i(C) \leq n-1$. So by induction $C \in \thick_{n-1}(\inj)$. Therefore $X \in \thick_n(\inj)$ as required. 
\end{proof}

\begin{theorem} \label{smallglobdim}
Let $\mathcal{A}$ be an $\inj$-coherent Grothendieck category. 
\begin{enumerate}
\item Then $\mathcal{D}^{\immac}(\mathcal{A}) = \thick_n(\inj)$ if and only if every $\inj$-finitely presented object has injective dimension at most $n$.

\item $\mathcal{D}_{\ifp}^b(\mathcal{A}) = \mathcal{D}^{\immac}(\mathcal{A})$ if and only if every $\inj$-finitely presented object has finite injective dimension.
\end{enumerate}
\end{theorem}

\begin{proof}
Both statements follow immediately from Lemma \ref{smallglobdimtostronggen} and Proposition \ref{smallinjgen}.
\end{proof}

 We consider when the standard $t$-structure restricts to $\mathcal{D}^{\immac}(\mathcal{A})$.

\begin{theorem} \label{t-structregular} Let $\mathcal{A}$ denote a Grothendieck category. The standard $t$-structure restricts to $\mathcal{D}^{\immac}(\mathcal{A})$ if and only if $\mathcal{A}$ is $\inj$-coherent and every $\inj$-finitely presented module has finite injective dimension. 
\end{theorem}

\begin{proof}
Suppose the standard $t$-structure restricts. Let $I \in \mathcal{D}^b_{\ipscoh}(\mathcal{A})$ which we may assume to be a bounded below complex of finite sums of indecomposable injectives. Let $n >> 0$ be such that $H^i(I) = 0$ for $i > n$. Then 
\[
I \simeq \tau_{\leq n} I \simeq \tau_{\leq n} \sigma_{\leq n+1} I
\]
Since $\sigma_{\leq n+1} I$ is immaculate, so is $\tau_{\leq n} \sigma_{\leq n+1} I \simeq I$ by assumption. Hence $\mathcal{D}^{\immac}(\mathcal{A}) = \mathcal{D}^b_{\ipscoh}(\mathcal{A})$. Therefore the standard $t$-structure restricts to $\mathcal{D}^b_{\ipscoh}(\mathcal{A})$ and $\mathcal{A}$ is $\inj$-coherent by Proposition \ref{injcohequiv}. Then by Theorem \ref{smallglobdim} (2), every $\inj$-finitely presented object has finite injective dimension. 

The converse follows from Theorem \ref{smallglobdim} (2) since the standard $t$-structure restricts to $\mathcal{D}^b_{\ifp}(\mathcal{A})$ by Proposition \ref{injcohequiv}.
\end{proof}

\section{Immaculate complexes as Modules over a Preadditive Category}
\label{preadditivesec}

In this section we note that the immaculate complexes are equivalent to the perfect derived category of a small preadditive category.

A preadditive category $\mathcal{C}$ is a category enriched in $\Mod \mathbb{Z}$. The category of right $\mathcal{C}$-modules $\Mod \mathcal{C}$ consist of additive functors $\mathcal{C}^{op} \to \Mod \mathbb{Z}$. Recall $\Mod \mathcal{C}$ is a Grothendieck category with a set of compact projective generators given by the representables. We write $\mathcal{D}(\mathcal{C}) \coloneq \mathcal{D}(\Mod \mathcal{C})$ and let $\mathcal{D}^{\perf}(\mathcal{C}) \subseteq \mathcal{D}(\mathcal{C})$ denote the bounded complexes of finitely generated projectives which coincide with the compact objects in $\mathcal{D}(\mathcal{C})$. We let $\mathcal{D}_{\pc}(\mathcal{C})$ denote the pseudocoherent complexes i.e.\ bounded above complexes of finitely generated projectives. We will let $\fdmod \mathcal{C}$ denote the finitely presented objects in $\Mod \mathcal{C}$ and $\mathcal{D}_{\fdmod}(\mathcal{C})$ the subcategory of complexes with finitely presented cohomology. We will use the superscripts $\{b,-\}$ to denote complexes with bounded and bounded above cohomology respectively.
\\

We say $\mathcal{C}$ has finite right global dimension if $\Mod \mathcal{C}$ has finite global dimension. We say $\mathcal{C}$ is right coherent if the kernel of a map between finitely generated projective right  modules is finitely generated (in the sense that it admits an epimorphism from a finitely generated projective). The same argument for a ring shows that, if $\mathcal{C}$ is right coherent, then $\mathcal{D}_{\pc}(\mathcal{C}) = \mathcal{D}^-_{\fdmod}(\mathcal{C})$ and the standard $t$-structure restricts to this category. As usual we will use the word left instead of right if the property holds for $\mathcal{C}^{op}$. 
\\

If $\mathcal{A}$ is a Grothendieck category, then $\inj = \inj(\mathcal{A})$ denotes the indecomposable injectives viewed as a preadditive subcategory of $\mathcal{A}$.

\begin{remark}
Note that the restricted Yoneda embedding 
\[
\mathcal{A}^{op} \to \Mod(\inj^{op}); M \mapsto \Hom_{\mathcal{A}}(M,-)
\]
is both left and right exact.
\end{remark}

\begin{remark}
The covariant functor $\mathcal{A} \to \Mod(\inj); M \mapsto \Hom_{\mathcal{A}}(-,M)$ is left exact and derives to a functor  
\[
Y: \mathcal{D}(\mathcal{A}) \to \mathcal{D}(\inj)
\]
The usual tilting arguments show $Y$ restricts to an equivalence
\[
\mathcal{D}^{\immac}(\mathcal{A}) \xrightarrow{\sim} \mathcal{D}^{\perf}(\inj)
\]
which identifies $\mathcal{D}(\inj)$ as the Ind-completion of $\mathcal{D}^{\immac}(\mathcal{A})$. This covariant functor is related to being Gorenstein and we plan to explore it further in later work. However the results below show the contravariant functor is better behaved in general. 
\end{remark}

\begin{theorem} \label{tiltingequiv1}
If $\mathcal{A}$ is a Grothendieck category, the restricted Yoneda embedding $\mathcal{A}^{op} \to \Mod (\inj^{op})$ induces equivalences
\[
\begin{tikzcd}
\mathcal{D}_{\ipscoh}(\mathcal{A})^{op} \arrow[r,"\sim"] & \mathcal{D}_{\pc}(\inj^{op}) \\
\mathcal{D}^{\immac}(\mathcal{A})^{op} \arrow[u,hook] \arrow[r,"\sim"] & \mathcal{D}^{\perf}(\inj^{op}) \arrow[u,hook]
\end{tikzcd}
\]

\end{theorem}

\begin{proof}
The restricted Yoneda embedding induces a functor which sends an $\inj$-pseudocoherent complex $I$ to the complex
\[
\dots \to \Hom_{\mathcal{A}}(I^{n-1},-) \to \Hom_{\mathcal{A}}(I^{n},-) \to 0
\]
Since each $I^i$ is a finite sum of indecomposable injectives, each $\Hom_{\mathcal{A}}(I^i,-)$ is a finite sum of corepresentables in $\Mod (\inj^{op})$ and thus a finitely generated projective. Therefore the functor restricts to $F: \mathcal{D}_{\ipscoh}(\mathcal{A})^{op} \to \mathcal{D}_{\pc}(\inj^{op})$. Any morphism $f: F(J) \to F(I)$ in $\mathcal{D}_{\pc}(\inj^{op})$ is a morphism between $K$-projectives and so may be represented by a chain map. So we have maps $f^i: \Hom_\mathcal{A}(J^i,-) \to \Hom_\mathcal{A}(I^i,-)$ in $\Mod \inj^{op}$. As the restricted Yoneda embedding is fully faithful on $\inj$, it is also fully faithful on finite sums of indecomposable injectives. Therefore each $f^i$ may be lifted to maps $\tilde{f}^i: I^i \to J^i$ such that $f^i = \Hom_{\mathcal{A}}(\tilde{f^i},-)$. As $f$ is a chain map we have that $\Hom_{\mathcal{A}}(\tilde{f^{i+1}}d_I^i,-) = \Hom_{\mathcal{A}}(d_J^i \tilde{f^{i}},-)$. Faithfulness of the restricted Yoneda embedding on finite sums of indecomposable injectives then implies that $\tilde{f}$ is a chain map. Hence $F$ is full. 

If $f: I \to J$ is sent to zero in $\mathcal{D}_{\ipscoh}(\inj^{op})$ then $F(f)$ must be null-homotopic as $F(I)$ and $F(J)$ are $K$-projectives. The components of any null homotopy can be lifted to maps $\tilde{s}^i: I^i \to J^{i-1}$ and as before fully faithfulness of the restricted Yoneda embedding implies that these $\tilde{s}^i$ form a null homotopy of $f$. Therefore $F$ is faithful. 

If $P$ is a bounded above complex of finitely generated projective left $\inj$-modules, then $P^i \simeq \Hom_{\mathcal{A}}(I^i,-)$ for a finite sum of indecomposable injectives $I^i$. Fully faithfulness of the restricted Yoneda embedding means the differentials $P^i \to P^{i+1}$ can be lifted to differentials $I^{i+1} \to I^{i}$. Hence $F$ is essentially surjective. It is clear from the proof that $F(I)$ is perfect if and only if $I$ is immaculate and thus the equivalence restricts as stated.
\end{proof}

\begin{theorem} \label{tiltingequiv2} Let $\mathcal{A}$ be an $\inj$-coherent Grothendieck category. Then $\inj$ is a left coherent preadditive category and there are $t$-exact equivalences
\[
\begin{tikzcd}
\mathcal{D}^{+}_{\ifp}(\mathcal{A})^{op} \arrow[r,"\sim"] & \mathcal{D}^{-}_{\fdmod}(\inj^{op}) \\
\mathcal{D}^b_{\ifp}(\mathcal{A})^{op} \arrow[u,hook] \arrow[r,"\sim"] &  \mathcal{D}^b_{\fdmod}(\inj^{op}) \arrow[u,hook]
\end{tikzcd}
\]
with respect to the standard $t$-structures on both sides. 
\end{theorem}

\begin{proof}
Let $F: \mathcal{D}_{\ipscoh}(\mathcal{A})^{op} \xrightarrow{\sim} \mathcal{D}_{\pc}(\inj^{op})$ denote the equivalence of Theorem \ref{tiltingequiv1}.  Since $F$ is induced by the restricted Yoneda embedding, which is exact,  we have that for any $I \in \mathcal{D}_{\ipscoh}(\mathcal{A})$,  $H^i(F(I)) \simeq \Hom_{\mathcal{A}}(H^{-i}(I),-)$. If $H^i(F(I)) = 0$, then as $\mathcal{A}$ is $\inj$-coherent, $\Hom_{\mathcal{A}}(H^{-i}(I),E(H^{-i}(I))) = 0$ and so $H^{-i}(I) = 0$. Hence the equivalence restricts to $\mathcal{D}^{\leq 0}_{\ipscoh}(\mathcal{A}) \simeq \mathcal{D}^{\geq 0}_{\pc}(\inj^{op})^{op}$ and $\mathcal{D}^{\geq 0}_{\ipscoh}(\mathcal{A}) \simeq \mathcal{D}^{\leq 0}_{\pc}(\inj^{op})^{op}$. As the standard $t$-structure restricts to $\mathcal{D}_{\ipscoh}(\mathcal{A})$ by Proposition \ref{injcohequiv}, it follows that the standard $t$-structure restricts to $\mathcal{D}_{\pc}(\inj^{op})$. The latter then implies that $\inj$ is left coherent. It is then clear that the equivalence restricts as shown and is $t$-exact. 
\end{proof}
We note a partial converse.

\begin{prop} \label{partialconverse}
Suppose $\mathcal{A}$ is a Matlis category. If $\inj$ is left coherent, then $\mathcal{A}$ is $\inj$-coherent.
\end{prop}

\begin{proof}
If $\inj$ is left coherent, the standard $t$-structure restricts to $\mathcal{D}_{\pc}(\inj^{op})$. The corresponding $t$-structure on $\mathcal{D}_{\ipscoh}(\mathcal{A})$ according to Theorem \ref{tiltingequiv1} has aisle consisting of complexes $I$ such that $\Hom_{\mathcal{A}}(H^{-i}(I),E) = 0$ for all $E \in \inj$ and $i \leq 0$. Since $E(H^{-i}(I))$ is a direct sum of indecomposable injectives, this occurs if and only if $I$ is in the aisle of the standard $t$-structure. Hence the standard $t$-structure restricts to $\mathcal{D}_{\ipscoh}(\mathcal{A})$ and $\mathcal{A}$ is $\inj$-coherent.
\end{proof}

\begin{cor} \label{Abelianequiv}
Suppose $\mathcal{A}$ is an $\inj$-coherent Grothendieck category. Then the Yoneda embedding induces an equivalence
\[
\ifp(\mathcal{A})^{op} \xrightarrow{\sim} \fdmod(\inj^{op})
\]
\end{cor}

\begin{proof}
This follows as the $t$-exact equivalence of Theorem \ref{tiltingequiv2} restricts to hearts. 
\end{proof}

Using this equivalence we can relate injective dimension of $\inj$-finitely presented objects to projective dimensions of the corresponding left $\inj$-modules. We recall that weak dimension is defined similarly to global dimension but projective resolutions are replaced with flat resolutions. 

\begin{prop} \label{idimvsinjweak}
Let $\mathcal{A}$ be an $\inj$-coherent Grothendieck category. Then
\begin{enumerate}
\item Every $\inj$-finitely presented object has finite injective dimension if and only if every finitely presented left $\inj$-module has finite projective dimension.

\item Every $\inj$-finitely presented object has injective dimension at most $n$ if and only if $\inj$ has weak dimension at most $n$.
\end{enumerate}
\end{prop}

\begin{proof}
An object $M \in \ifp(\mathcal{A})$ has finite injective dimension as an object of $\mathcal{A}$ if and only if it admits a finite injective resolution by finite sums of indecomposable injectives. This is equivalent, by Corollary \ref{Abelianequiv}, to its image in $\fdmod(\inj^{op})$ admitting a projective resolution of the same length. Hence (1) holds and the first part of (2) is equivalent to every finitely presented $\inj^{op}$ module having projective dimension at most $n$. As is the case for rings, this is the same has having flat dimension at most $n$. 
\end{proof}

\begin{remark}
Proposition \ref{notnoetheg} and Remark \ref{notnoetheg2} shows that $\inj$ can easily fail to be left or right Noetherian and so weak dimension is a more natural invariant to consider than global dimension. In fact, Proposition \ref{dvreg} shows that one cannot replace weak dimension with left global dimension in the statement of Proposition \ref{idimvsinjweak}.
\end{remark}

\section{Regularity and Dimension via Immaculate Complexes}
\label{immaculateschemes}

In this section we investigate the immaculate complexes $\mathcal{D}^{\immac}(X) \coloneq \mathcal{D}^{\immac}(\QCoh(X))$ over a locally Noetherian scheme $X$. We quickly note that it doesn't matter if we work in $\mathcal{D}(\QCoh(X))$ or $\mathcal{D}_{\QCoh}(X)$. The result below follows from \cite[Corollary 7.19]{Har66}.

\begin{prop}
Let $X$ be a locally Noetherian scheme. There are equivalences
\begin{align*}
\mathcal{D}^{\immac}(\QCoh(X)) & \xrightarrow{\sim} \mathcal{D}^{\immac}(\Mod \mathcal{O}_X) \cap \mathcal{D}_{\QCoh}(X) \\
\mathcal{D}_{\ipscoh}(\QCoh(X)) & \xrightarrow{\sim} \mathcal{D}_{\ipscoh}(\Mod \mathcal{O}_X) \cap \mathcal{D}_{\QCoh}(X) 
\end{align*}
\end{prop}

Matlis theory gives us a good understanding of the immaculate complexes. By Theorem \ref{matlisthy},
\[
\mathcal{D}^{\immac}(X) = \thick(E(x) \mid x \in X )
\]
where $E(x) \coloneq E(k(x))$ denotes the injective hull of the residue field at $x$. In this section we will prove the following in order to apply the results of Section \ref{immacsec} and \ref{preadditivesec}. See Definition \ref{injcohdefn}
for the notion of $\inj$-coherence.

\begin{theorem} \label{injcohscheme}
Let $X$ be a locally Noetherian scheme, then $\QCoh(X)$ is $\inj$-coherent.
\end{theorem}

\begin{notation} If $I$ is an injective over a locally Noetherian scheme $X$, we will say that $x \in X$ appears in $I$ to mean that $E(x)$ is a summand of $I$. We also let $\max(I)$ denote the points $x$ appearing in $I$ such that if $x$ specialises to some $y \neq x$, then $y$ does not appear in $I$. By the height of a point $x \in X$ we mean $ \hgt(x) \coloneq \dim \mathcal{O}_{X,x} = \codim_X \overline{\{x\}}$. 
\end{notation}

\begin{lemma} \label{weakinjcohlemma}
Let $R$ be a commutative Noetherian ring and $f: I \to J$ a map between finite sums of indecomposable injectives. If a prime $p$ appears in $E(\coker(f))$, then $p$ is contained in a prime appearing in $J$.  
\end{lemma}

\begin{proof}
We can assume that every prime appearing in $I$ is contained in a prime appearing in $J$. Indeed if $p$ is not contained in any prime appearing in $J$, then $\Hom_R(E(R/p),J) = 0$ by Proposition \ref{matlisprop}. So $f$ must be zero restricted to any such summand of $I$ and $K \coloneq \coker(f)$ is unchanged by removing it.

Let $S = R \setminus \bigcup_{p \in \max(J)} p$. Any $s \in S$ is not contained in any prime appearing in $I$ or $J$ and so acts invertibly on $I$ and $J$. Hence $f$ is a map of $S^{-1}R$-modules and so $K$ and $E(K)$ are too. Since $\max(J)$ is finite, $S^{-1}R$ is a semilocal ring with maximal ideals $\max(J)$. Hence the only primes which can appear in $E(K)$ are contained in primes in $\max(J)$.
\end{proof}

\begin{lemma} \label{stronginjcohlemma}
Let $X$ be a locally Noetherian scheme and $f: I \to J$ a map between finite sums of indecomposable injectives. If $x \in X$ appears in $E(\coker(f))$, then $x$ appears in $J$. So if $M \in \ifp(X)$ with minimal injective resolution $I$, then there is a finite set $\mathcal{P} \subseteq X$ such that if $x$ appears in $I^n$ for any $n \geq 0$, then $x \in \mathcal{P}$. 
\end{lemma}

\begin{proof}
As in the proof of Lemma \ref{weakinjcohlemma}, we may assume that every point appearing in $I$ specialises to a point appearing in $J$. Let $J = \bigoplus E(y)^{\mu(y)}$ and $K \coloneq E(\coker(f))$. Suppose $x$ is a point appearing in $E(K)$. There is an exact sequence of $\mathcal{O}_{X,x}$-modules
\[
I_x \to J_x \to K_x \to 0
\]
and $J_x \simeq \bigoplus_{p_y \subseteq p_x} E(R/p_y)^{\mu(y)}$ where $p_y$ is the prime of $R \coloneq \mathcal{O}_{X,x}$ corresponding to $y \in X$ for any $y$ with $x \in \overline{\{y\}}$. Since $x$ appears in $E(K)$, the maximal ideal $p_x$ of $\mathcal{O}_{X,x}$ appears in $E(K)_x = E(K_x)$. Hence by Lemma \ref{weakinjcohlemma}, there is a prime $q$ appearing in $J_x$ such that $p_x \subseteq q$. However every prime appearing in $J_x$ is contained in $p_x$ and so we must have that $p_x = q$. Hence $p_x$ appears in $J_x$. Therefore $\mu(x) \neq 0$. 
\end{proof}

\begin{remark} 
Lemma \ref{stronginjcohlemma} shows that if $M \in \ifp(R)$ and $p \in \Spec(R)$ appears in its minimal injective resolution $I$ at some index $I^n$ for $n \geq 1$, then $p$ appears in $I^1$. This is clearly very different to the usual situation for injective resolutions of finitely generated modules. We note in Section \ref{perfimmacsec} that very few finitely generated modules are $\inj$-finitely presented.  
\end{remark}

\begin{lemma} \label{injcohlemma}
Let $X$ be a locally Noetherian scheme and $x \in X$ and $f: E(x)^{\mu} \to E(x)^{\lambda}$ for finite $\mu,\lambda \geq 0$. Then $E(\coker(f))$ is a finite sum of copies of $E(x)$. 
\end{lemma}

\begin{proof}
Let $K \coloneq \coker(f)$. Then note that $f$ is in the image of the fully faithful exact functor $i_\ast: \Mod \mathcal{O}_{X,x} \hookrightarrow \QCoh(X)$ which is the right adjoint to taking stalks. Since $E(x)$ is an Artinian $\mathcal{O}_{X,x}$ module (e.g.\ \cite[Theorem 3.4.3]{EJ11}), so is $K_x$. Hence $E(K_x) = E(K)_x$ is a finite sum of copies of $E(x)$. By Lemma \ref{stronginjcohlemma}, no other points can appear in $E(K)$ and so we are done.
\end{proof}

\begin{proof}[Proof of Theorem \ref{injcohscheme}]
Suppose $f: I \to J$ is map between finite sums of indecomposable injectives and let $K \coloneq \coker(f)$. As in the proof of Lemma \ref{weakinjcohlemma}, we can assume that every point appearing in $I$ specialises to a point appearing in $J$. Let $I \simeq \bigoplus E(x)^{\mu(x)}$ and  $J \simeq \bigoplus E(x)^{\lambda(x)}$. 

We show that $E(K)$ is a finite sum of indecomposable injectives by inducting on the heights of the points appearing in $J$. Suppose all points appearing in $J$ (and so also in $I$) are of height zero. If $\hgt(x) = \hgt(y)$ but $x \neq y$ then $\Hom_X(E(x),E(y)) = 0$ by Proposition \ref{matlisprop}. So $f: I \to J$ splits as a finite direct sum of maps $f_x : E(x)^{\mu(x)} \to E(x)^{\lambda(x)}$. By Lemma \ref{injcohlemma}, we see that each $E(\coker(f_x))$ is a finite sum of indecomposable injectives and that $E(K) = \bigoplus E(\coker(f_x))$ is too.

Let $n$ be the maximum height of any point appearing in $J$ and suppose the claim holds for strictly smaller $n$. By Lemma \ref{stronginjcohlemma}, every point appearing in $E(K)$ appears in $J$ and so it is enough to show that $E(K)_x$ is a finite sum of indecomposable injectives for each $x \in \max(J)$. Indeed, if $y$ appears in $E(K)$, then $y$ specialises to some $x \in \max(J)$. Then the multiplicity of $E(y)$ appearing in $E(K)$ is the same as the multiplicity of $E(y)$ appearing in $E(K)_x$ by Proposition \ref{matlisprop}. Therefore we may fix some $x \in \max(J)$ and assume $R = \mathcal{O}_{X,x}, I = I_x, J = J_x, K = K_x$ and let $m$ denote the maximal ideal of $R$. If $\hgt(x) < n$, then we are done by induction and so we assume that $\hgt(x) = \dim(R) = n$. 

Now Let $I^m = E(R/m)^{\mu(m)} \subseteq I$ and $I^{<m} = \bigoplus_{p \subset m} E(R/p)^{\mu(p)} \subseteq I$ and similarly define $J^m$ and $J^{<m}$. Let $
M \coloneq f^{-1}(J^m) \cap I^{< m} \subseteq I^{<m}$. Then $f(M \bigoplus I^m) \subseteq J^m$ and so $f$ descends to a map $\overline{f}: I^{<m}/M \to J/J^m = J^{<m}$. Note that $\overline{f}$ is mono since any element of $I^{<m}$ sent to $J^m$ is contained in $M$. We let  $K^{<m} \coloneq \coker(\overline{f})$. A diagram chase shows that the following diagram commutes and has exact rows and columns. 
\[
\begin{tikzcd}[row sep = small]
  & 0 \arrow[d] & 0 \arrow[d] & 0 \arrow[d] & \\
& M \bigoplus I^m \arrow[d]
\arrow[r] & J^{m} \arrow[d] \arrow[r] & K^m \arrow[r]  \arrow[d] & 0  \\
& I \arrow[d] \arrow[r,"f"]  & J \arrow[r] \arrow[d] & K \arrow[r] \arrow[d] & 0 \\
0 \arrow[r]  & I^{<m}/ M \arrow[r,"\overline{f}"] \arrow[d] & \arrow[d] J^{<m} \arrow[r]& K^{<m} \arrow[r] \arrow[d] & 0 \\
  &  0   & 0   & 0 
\end{tikzcd}
\]
where we define $K^m$ to be the kernel of the obvious epi $K \to K^{<m}$. Exactness of the top row at $J^m$ uses that $\overline{f}$ is mono. 

As $J^m$ is an Artinian $R$-module, $K^m$ is too and thus $E(K_m)$ is a finite sum of indecomposable injectives (e.g.\ by \cite[Theorem 3.4.3]{EJ11}). As $f(I^m) \subseteq J^m$, $f$ descends to a map $f^{<m}:  I^{<m} = I/I^m \to  J/J^m = J^{<m}$. Since $f^{<m}(M) = 0$, we see that $f^{<m}$ factors through $\overline{f}$ via the quotient $I^{<m} \to I^{<m}/M$. Hence $K^{<m} \simeq \coker(f^{<m})$. Now $f^{<m}$ is a map between finite sums of indecomposable injectives and the height of any prime appearing in $J^{<m}$ is strictly less than $n$. By induction $E(K^{<m}) \simeq E(\coker(f^{<m}))$ is a finite sum of indecomposable injectives. By the Horseshoe Lemma applied to the right most column in the above diagram, $K$ embeds into $E(K^m) \bigoplus E(K^{<m})$ which is a finite sum of indecomposable injectives. Hence $E(K)$ is also a finite sum of indecomposable injectives as required. 
\end{proof}

\begin{lemma}\label{residueifp}
Let $X$ be a locally Noetherian scheme and $x \in X$. Then $k(x) \in \mathcal{D}^b_{\ifp}(X)$. 
\end{lemma}

\begin{proof}
Consider the minimal injective resolution $I$ of $k(x)$. If $y \not\in \overline{\{x\}}$, then $k(x)_y = 0$ and so by Lemma \ref{bassnumbers}, $y$ does not appear in $I$. Otherwise we have that $\Ext^i_{\mathcal{O}_{X,y}}(k(y),k(x)_y) = \Ext^i_{\mathcal{O}_{X,x}}(k(y)_x,k(x))$ which is non-zero only if $y = x$.  So only $x$ appears in the minimal injective resolution of $k(x)$. Since $k(x)$ is finite over $\mathcal{O}_{X,x}$, Lemma \ref{bassnumbers} implies the terms of $I$ are finite sums of copies of $E(x)$. 
\end{proof}

\begin{theorem} \label{regularscheme}
Suppose $X$ is a locally Noetherian scheme. 
\begin{enumerate}
\item The standard $t$-structure restricts to $\mathcal{D}_{\ipscoh}(X) = \mathcal{D}^+_{\ifp}(X)$ and to $\mathcal{D}^b_{\ipscoh}(X) = \mathcal{D}^b_{\ifp}(X)$.

\item $X$ is regular and of dimension at most $n$ if and only if $\mathcal{D}^{\immac}(X) = \thick_n(\inj)$. 

\item $X$ is regular if and only if $\mathcal{D}^{\immac}(X) = \mathcal{D}^b_{\ifp}(X)$ if and only if the standard $t$-structure restricts to $\mathcal{D}^{\immac}(X)$.
\end{enumerate}
\end{theorem}

\begin{proof}
By Theorem \ref{injcohscheme}, $\QCoh(X)$ is $\inj$-coherent. Proposition \ref{injcohequiv} implies (1). By Theorem \ref{smallglobdim} (1), $\mathcal{D}^{\immac}(X) = \thick_n(\inj)$ if and only if every $\inj$-finitely presented object has injective dimension at most $n$. By Lemma \ref{residueifp}, this implies that the residue fields each have injective dimension at most $n$. Then each local ring, and so also $X$, is regular of dimension at most $n$. The converse follows from Theorem \ref{ABSscheme} and so (2) holds.

By Theorems \ref{smallglobdim} (2) and \ref{t-structregular}, it remains to show that $X$ is regular if and only if every $\inj$-finitely presented object has finite injective dimension. Again by Lemma \ref{residueifp}, the latter condition implies $X$ is regular. Suppose $X$ is regular and $M \in \ifp(X)$ with minimal injective resolution $I$. By Lemma \ref{stronginjcohlemma}, only finitely many points $\mathcal{P} \subseteq X$ appear in $I$. So $\mu_i(x,M) = 0$ for $x \not\in\mathcal{P}$ and $i \geq 0$. Each local ring is regular of finite-dimension and let $N$ denote the supremum of the dimensions of $\mathcal{O}_{X,x}$ for $x \in \mathcal{P}$. By Lemma \ref{bassnumbers}, $\mu_i(x,M) = 0$ for $i > N$ and $x \in \mathcal{P}$. Therefore $I^i = 0$ for $i > N$ and $M$ has finite injective dimension. 
\end{proof}

We can also phrase our statements in terms of the small preadditive category $\inj$ of Section \ref{preadditivesec}.

\begin{remark}
If $X$ is a locally Noetherian scheme, then the small preadditive category $\inj(X)$ has objects indexed by $X$ and morphisms given by $\Hom_X(E(x),E(y))$. By Proposition \ref{matlisprop}, $\Hom_X(E(x),E(y)) = 0$ if $y \not \in \overline{ \{x\}}$ and so $\inj(X)$ is directed according to the topology of $X$. 
\end{remark}

\begin{theorem} \label{weakdimscheme}
Let $X$ be a locally Noetherian scheme. 
\begin{enumerate}
\item $\inj(X)$ is left coherent.

\item $X$ is regular of dimension $n$ if and only if $\inj(X)$ has weak dimension $n$.

\item $X$ is regular if and only if every finitely presented left $\inj(X)$ module has finite projective dimension.
\end{enumerate}
\end{theorem}

\begin{proof}
By Theorems \ref{injcohscheme} and \ref{tiltingequiv2}, $\inj(X)$ is left coherent. Parts (2) and (3) follow from Proposition \ref{idimvsinjweak} and the proof of Theorem \ref{regularscheme}.
\end{proof}

\begin{remark} \label{openqsscheme}We consider a few open problems.
\begin{enumerate}

\item There is a similar structure theory for injectives over sheaves of noncommutative algebras over Noetherian schemes due to \cite{BD22}. We expect similar arguments to show these categories are also $\inj$-coherent. This would also include group algebras of finite groups over commutative Noetherian rings. 

\item We don't know if $\inj(X)$ is also right coherent. This is particularly relevant for studying the immaculates over Gorenstein schemes. Here it is natural to ask for a duality $\mathcal{D}^b_{\fdmod}(\inj(X))^{op} \simeq \mathcal{D}^b_{\fdmod}(\inj(X)^{op})$ and to study the covariant Yoneda embedding, $\mathcal{D}^b_{\ifp}(X) \to \mathcal{D}^b_{\fdmod}(\inj(X))$.

\item Theorem \ref{regularscheme} shows that the vanishing of the Verdier quotient $\mathcal{D}^b_{\ifp}(X)/\mathcal{D}^{\immac}(X)$ measures singularity of $X$ and thus provides an immaculate version of the singularity category $\mathcal{D}^{b}_{\fdmod}(X)/\mathcal{D}^{\perf}(X)$. One can ask if analogous properties to those studied in \cite{Buch21} hold. For example, is it equivalent to the stable category of a Frobenius exact category if $X$ is Gorenstein ? See also Corollary \ref{idempotentcomplete}.

\item It was shown in \cite{Nee24t-struct} that $X$ is regular if and only if $\mathcal{D}^{\perf}(X)$ admits a bounded $t$-structure. Does the analogous statement hold for $\mathcal{D}^{\immac}(X)$? Note Theorem \ref{regularscheme} (3) refers only to the standard $t$-structure and the approximability methods of \cite{Nee24t-struct} do not immediately apply as $\mathcal{D}^{\immac}(X)$ does not usually have a single generator. 
\end{enumerate}

\end{remark}

\section{Perfect and Immaculate complexes}
\label{perfimmacsec}

In this short section we compute the intersection of the immaculate and the perfect complexes. Let $\fsupp(X) \subseteq \coh(X)$ denote the coherent sheaves with finite support over a locally Noetherian scheme $X$. The support of any sheaf in $\fsupp(X)$ can be at most 1 dimensional and any point in the support must be either closed or specialise to only finitely many closed points. We let $\supp$ denote the support of a sheaf.

\begin{lemma} \label{imperflem1}
Let $X$ be a locally Noetherian scheme and $M \in \QCoh(X)$. 
\begin{enumerate}
\item If $x \in \supp(M)$, then there is a $y \in \supp(M)$ appearing in $E(M)$ which specializes to $x$.

\item If $x$ appears in the injective resolution of $M$, then $x \in \supp(X)$. 
\end{enumerate}

\end{lemma}

\begin{proof}
Consider $M_x \in \Mod R$ where $R = \mathcal{O}_{X,x}$. By \cite[Tag 05BV]{stacks-project}, the prime corresponding to $x$ contains an associated prime $q$ of $M_x$. Hence $R/q \subseteq M_x \subseteq E(M_x)$. So $q$ appears in $E(M_x)$ and by Proposition \ref{stalkwise}, the point $y \in X$ corresponding to $q$ appears in $E(M)$. If $\mu_n(x,M) \neq 0$ for some $n$, then by Lemma \ref{bassnumbers}, we see that $M_x \neq 0$. 
\end{proof}

\begin{lemma}\label{imperflem2}
Let $X$ be a locally Noetherian scheme and $M \in \coh(X)$. If $x$ appears in the injective resolution of $M$ and $x$ specialises to $y \in X$, then $y$ appears in the injective resolution of $M$.
\end{lemma}

\begin{proof}
Induct on $\hgt(y) - \hgt(x)$. If they coincide there is nothing to prove. Otherwise we can consider a point $x'$ such that $x$ specialises to $x'$ and $x'$ specialises to $y$ such that $\hgt(y) = \hgt(x') + 1$. By induction $x'$ appears in the injective resolution of $M$. So there is some $n$ such that $\mu_n^X(x',M) = \mu_n^{\mathcal{O}_{X,y}}(x',M_y) \neq 0$ by Proposition \ref{stalkwise}. But then by  \cite[Proposition 9.2.13]{EJ11}, $\mu_{n+1}^X(y,M) = \mu_{n+1}^{\mathcal{O}_{X,y}}(y,M_y) \neq 0$. Hence $y$ appears in the injective resolution of $M$. 
\end{proof}

\begin{lemma} \label{imperflem3}
Let $X$ be a locally Noetherian scheme and $M \in  \coh(X)$. A point appears in the injective resolution of $M$ if and only if it is in the support of $M$.
\end{lemma}

\begin{proof}
If $x \in \supp(M)$, then by Lemma \ref{imperflem1} (1), there is a point $y$ specialising to $x$ appearing in $E(M)$. So by Lemma \ref{imperflem2}, $x$ appears in the injective resolution of $M$. The converse is Lemma  \ref{imperflem1} (2).
\end{proof}

\begin{prop} \label{intersectionchar}
Let $X$ be a locally Noetherian scheme. Then $\ifp(X) \cap \coh(X) = \fsupp(X)$ and $\mathcal{D}^b_{\ifp}(X) \cap \mathcal{D}^b_{\coh}(X) = \mathcal{D}^b_{\fsupp}(X)$.
\end{prop}

\begin{proof}
If $M$ is a coherent sheaf, then by Lemma \ref{bassnumbers} each $\mu_i(x,M)$ is finite. Hence $M \in \ifp(X)$ if and only if only finitely many points appear in its injective resolution. By Lemma \ref{imperflem3} this occurs if and only if it has finite support. The second statement then follows immediately. 
\end{proof}

\begin{remark}
Using Proposition \ref{intersectionchar}, one can describe the other intersections such as $\mathcal{D}^{\immac}(X) \cap \mathcal{D}^{\perf}(X)$ using the characterisations of $\mathcal{D}^{\immac}(X)$ within $\mathcal{D}^b_{\ifp}(X)$ from Lemma \ref{immacchar} and $\mathcal{D}^{\perf}(X)$ within $\mathcal{D}^b_{\coh}(X)$ as in \cite[08CQ]{stacks-project}.
\end{remark}

\section{Semi-orthogonal Decompositions and Additive Invariants}
\label{SODsec}

In this section we note that if $X$ is a locally Noetherian scheme, then Matlis theory gives us a very good structural understanding of $\mathcal{D}^{\immac}(X)$ using semi-orthogonal decompositions. See \cite{Kuz14} for a survey and the relevant definitions. Recall a subset $ \mathcal{P} \subseteq X$ is specialisation closed if whenever $y \in \overline{\{x\}}$ and $x \in \mathcal{P}$, then $y \in \mathcal{P}$. 

\begin{lemma} \label{SODlemma}
If $\mathcal{P}$ is a specialisation closed subset of a locally Noetherian scheme $X$, then there is a semi-orthogonal decomposition
\[
\mathcal{D}^{\immac}(X) = \langle \thick(E(y) \mid y \not\in \mathcal{P}), \thick(E(x) \mid x \in \mathcal{P}) \rangle
\]
\end{lemma}

\begin{proof}
We must show there are no morphisms from right to left. If $x \in \mathcal{P}$ and $y \not\in \mathcal{P}$ then $y \not\in \overline{\{x\}}$. So by Proposition \ref{matlisprop} and since they are injective, $\Hom_{\mathcal{D}^{\immac}(X)}(E(x),\Sigma^i E(y)) = 0$ for all $i$. It follows that whenever $M \in \thick(E(x) \mid x \in \mathcal{P})$ and $N \in \thick(E(y) \mid y \notin \mathcal{P})$, then $\Hom_{\mathcal{D}^{\immac}(X)}(M,N) = 0$. By Proposition \ref{immacsarethick}, the indecomposable injectives generate $\mathcal{D}^{\im}(X)$ as a triangulated category as required. 
\end{proof}

\begin{remark} \label{explicitapprox}
If $M \in \mathcal{D}^{\immac}(X)$ is a bounded complex of finite sums of indecomposable injectives, let $M(\mathcal{P}) \subseteq M$ denote the subcomplex of $M$ consisting of the injectives of the form $E(x)$ with $x \in \mathcal{P}$. Then the approximation triangle associated to the semi-orthogonal decomposition of Lemma \ref{SODlemma} comes from the short exact sequence
\[
0 \to M(\mathcal{P}) \to M \to M/M(\mathcal{P}) \to 0
\]
\end{remark}

\begin{remark}
There are semi-orthogonal decompositions of $\mathcal{D}^{\immac}(X)$ which are not of the form in Lemma \ref{SODlemma}. See Proposition \ref{SODeg}.
\end{remark}

\begin{remark} \label{closedlocalisation}
If $x$ is a closed point of $X$, then taking $\mathcal{P} = \{x\}$ in Lemma \ref{SODlemma} implies there is a localisation 
\[
\mathcal{D}^{\immac}(X) \twoheadrightarrow \thick(E(x)) \simeq \mathcal{D}^{\perf}(\widehat{\mathcal{O}}_{X,x})
\]
where the equivalence follows from Proposition \ref{matlisprop} (2) and the usual tilting arguments. This gives an alternative proof that strong generation of $\mathcal{D}^{\immac}(X)$ implies $X$ is regular and of finite dimension. 
\end{remark}

Recall that if $X$ is a quasicompact quasiseparated scheme and $U\subseteq X$ is a quasicompact open, then there is a localisation-up-to-summands functor $\mathcal{D}^{\perf}(X) \twoheadrightarrow \mathcal{D}^{\perf}(U)$ (due to \cite{TT90}). For immaculates this becomes a Bousfield localisation i.e. there is a fully faithful right adjoint. This can be viewed as a special case of Lemma \ref{SODlemma} with $\mathcal{P} = X \setminus U$.

\begin{prop} Let $X$ be a locally Noetherian scheme and $j:U \hookrightarrow X$ an open immersion. Then there is an adjoint pair
\[
\begin{tikzcd}
\mathcal{D}^{\immac}(X) \arrow[r,two heads,shift left = 1ex,"j^\ast"] & \arrow[l,hook',shift left = 1ex,"\mathbb{R}j_\ast"] \mathcal{D}^{\immac}(U)
\end{tikzcd}
\]
where $\mathbb{R}j_\ast$ is fully faithful with image $\thick(E(x)\mid x \in U)$ and $j^\ast$ is a localisation.
\end{prop}

\begin{proof}
Recall that $\mathbb{R}j_\ast$ is fully faithful and right adjoint to $j^\ast$ as functors between $\mathcal{D}(\Mod \mathcal{O}_X)$ and $\mathcal{D}(\Mod \mathcal{O}_U)$. Let $i_U$ and $i_X$ denote the inclusions of $\Spec(\mathcal{O}_{X,x})$ into $U$ and $X$ respectively. Then 
\[
\mathbb{R}j_\ast E_U(x) = j_\ast E_U(x) = j_\ast (i_U)_{\ast} E_{\mathcal{O}_{U,x}}(x) = (i_X)_{\ast} E_{\mathcal{O}_{X,x}}(x) = E_X(x)
\]
Therefore $\mathbb{R}j_\ast$ restricts to a fully faithful functor on the immaculates with image as described. If $x \in U$, then $j^\ast E_X(x) = j^\ast \mathbb{R}j_\ast E_U(x) \simeq E_U(x)$. If $x \not\in U$, then $j^\ast E_X(x) = 0$ since $E(x)$ has support in the complement of $U$. Therefore $j^\ast$ restricts to the immaculates as well. As it admits a fully faithful right adjoint, $j^\ast$ is a localisation functor. 
\end{proof}

\begin{remark}
We note that unlike for the perfects, the functor $j^\ast: \mathcal{D}^{\immac}(X) \to \mathcal{D}^{\immac}(U)$ is essentially surjective without having to take direct summands. 
\end{remark}

Recall that if $x \in X$ is a point of a locally Noetherian scheme, $\hgt(x) =  \dim(\mathcal{O}_{X,x})$.

\begin{theorem} \label{SODfiltration}
Let $X$ be a locally Noetherian scheme. Let $\mathcal{P}_{i} \subseteq X$ consist of all points $x$ such $\hgt(x) = i$. For any $i$, let
\[
\mathcal{F}_i = \thick(E(x) \mid x \in \mathcal{P}_i) \quad \mathcal{F}_{>i} = \thick(E(x) \mid x \in \mathcal{P}_{j} \text{ for } {j >i})
\]
For any $n>0$, there is a semi-orthogonal decomposition
\[
\mathcal{D}^{\immac}(X) = \langle \mathcal{F}_0,\dots ,\mathcal{F}_n,\mathcal{F}_{>n} \rangle
\] 
and $\mathcal{F}_i \simeq \bigoplus_{x \in \mathcal{P}_i} \mathcal{D}^{\perf}(\widehat{\mathcal{O}}_{X,x})$. In particular if $X$ is finite-dimensional, then there is a semi-orthogonal decomposition
\[
\mathcal{D}^{\immac}(X) = \langle \mathcal{F}_0,\dots ,\mathcal{F}_{\dim(X)} \rangle
\] 
\end{theorem}

\begin{proof}
The existence of the semi-orthogonal decomposition is similar to Lemma \ref{SODlemma}. For each $i$, there is an orthogonal decomposition
$\mathcal{F}_i \simeq \bigoplus_{x \in \mathcal{P}_i} \thick(E(x))$ since if $x \neq y$ but $\hgt(x) = \hgt(y)$, then Proposition \ref{matlisprop} (1) implies $\Hom_X(E(x),E(y)) = 0$. The equivalences follow from standard tilting arguments since $\Hom_X(E(x),E(x)) \simeq \widehat{\mathcal{O}}_{X,x}$ by Proposition \ref{matlisprop} (2). 
\end{proof}

\begin{remark}
To consider invariants, we recall that $\mathcal{D}^{\immac}(X)$ admits an enhancement as a stable $\infty$-category using the natural enhancement on $\mathcal{D}_{\QCoh}(X)$. As an abuse of notation we will write $\mathcal{D}^{\immac}(X)$ for the stable $\infty$-category as well as the triangulated category. See \cite{Lur09} for more details. 
\end{remark}

\begin{remark}
Suppose $\mathcal{D}$ is a stable $\infty$-category such as spectra or the derived category of a ring. A finitary additive invariant taking values in $\mathcal{D}$ associates to any small stable $\infty$-category an object of $\mathcal{D}$ in a way which sends semi-orthogonal decompositions to direct sum decompositions and preserves filtered colimits. Examples include (non-) connective $K$-theory and (topological) Hochschild homology. See \cite{BGT13} for precise definitions but note they use the term additive invariants for what we call finitary additive invariants.
\end{remark}

\begin{theorem} \label{additivethm}
Let $U$ be any finitary additive invariant and $X$ a locally Noetherian scheme. Then 
\[
U(\mathcal{D}^{\immac}(X)) \simeq \bigoplus_{x \in X} U(\mathcal{D}^{\perf}(\widehat{\mathcal{O}}_{X,x}))
\]
\end{theorem}

\begin{proof}
For $n \geq 0$, let $\mathcal{F}_{n}$ be as in Theorem \ref{SODfiltration} and $\mathcal{F}_{\leq n} = \thick(E(x) \mid \hgt(x) \leq n)$. Then $\mathcal{D}^{\immac}(X) = \bigcup_{n \geq 0} \mathcal{F}_{\leq n}$ and this is a filtered colimit of stable $\infty$-categories. Each $\mathcal{F}_{\leq n}$ admits a semi-orthogonal decomposition $\mathcal{F}_{\leq n} = \langle \mathcal{F}_{\leq n-1}, \mathcal{F}_n\rangle$ by the same arguments as in Lemma \ref{SODlemma}. Since $U$ is finitary and additive we have that
\[
U(\mathcal{D}^{\immac}(X)) \simeq \colim_n U(\mathcal{F}_{\leq n}) \simeq \colim_n \bigoplus_{0 \leq i \leq n} U(\mathcal{F}_i) \simeq \bigoplus_{n \geq 0} U(\mathcal{F}_n)
\]
As in Theorem \ref{SODfiltration}, each $\mathcal{F}_n$ is the coproduct $\bigoplus_{\hgt(x) = n} \mathcal{D}^{\perf}(\widehat{\mathcal{O}}_{X,x})$ of small stable $\infty$-categories. This infinite coproduct is a filtered colimit of finite coproducts. So by the same argument as above $U(\mathcal{F}_n) \simeq \bigoplus_{\hgt(x) = n} U(\mathcal{D}^{\perf}(\widehat{\mathcal{O}}_{X,x}))$. Since every point of $X$ has a finite height, we see that $\bigoplus_{n \geq 0} \bigoplus_{\hgt(x) = n} U(\mathcal{D}^{\perf}(\widehat{\mathcal{O}}_{X,x})) = \bigoplus_{x \in X} U(\mathcal{D}^{\perf}(\widehat{\mathcal{O}}_{X,x}))$ as required.
\end{proof}

\begin{remark}  \label{Kremark}
In particular Theorem \ref{additivethm} applies to $K$-theory. See \cite[Section 9]{BGT13} for a construction of the nonconnective $K$-theory spectrum $\mathbb{K}(\mathcal{D})$ of a stable $\infty$-category $\mathcal{D}$. If $R$ is a ring and $n \in \mathbb{Z}$, then $\pi_n\mathbb{K}(\mathcal{D}^{\perf}(R)) \simeq K_n(R)$ where the latter denotes the $K$-groups of $R$ which admit a number of classical descriptions as in \cite{Wei13}. We write $K_n(\mathcal{D}^{\immac}(X)) \coloneq \pi_n \mathbb{K}(\mathcal{D}^{\immac}(X))$. Also recall that $K_0(\mathcal{D}^{\immac}(X))$ coincides with the triangulated Grothendieck group of $\mathcal{D}^{\immac}(X)$.
\end{remark}

\begin{cor} \label{K-theory}
Let $X$ be a locally Noetherian scheme, there are isomorphisms for all $n \in \mathbb{Z}$
\[
K_n(\mathcal{D}^{\immac}(X)) \simeq \bigoplus_{x \in X} K_n(\widehat{{\mathcal{O}}_{X,x}})
\]
and in particular 
\begin{enumerate}
\item $K_0(\mathcal{D}^{\immac}(X))$ is the free Abelian group of rank given by the cardinality of $X$.
\item If $X$ is regular, then $K_n(\mathcal{D}^{\immac}(X)) = 0$ for $n < 0$.

\item If $X$ is finite-dimensional $K_n(\mathcal{D}^{\immac}(X)) = 0$ for $ n < - \dim(X)$.
\end{enumerate}
\end{cor}

\begin{proof}
By \cite[Theorem 9.8]{BGT13}, nonconnective $K$-theory is a finitary additive invariant so Theorem \ref{additivethm} applies. By Remark \ref{Kremark}, $\pi_n\mathbb{K}(\mathcal{D}^{\perf}(R)) \simeq K_n(R)$ for any ring. Since homotopy groups commute with coproducts the first claim holds. The Grothendieck group of any local ring is $\mathbb{Z}$ and so (1) holds. If $X$ is regular, then so is $\widehat{\mathcal{O}}_{X,x}$ for each $x \in X$. The negative $K$-groups of any regular ring vanish (e.g.\ see \cite[III.4]{Wei13}) and so (2) holds. By \cite[Theorem B]{KST19}, the negative $K$-groups of any Noetherian ring $R$ vanish below $\dim(R)$ and so (3) holds. 
\end{proof}

Using similar ideas to \cite{Orl11}, we can use these negative $K$-theory results to discuss idempotent completion of the immaculate singularity category $\mathcal{D}_{\isg}(X) \coloneq \mathcal{D}^b_{\ifp}(X)/\mathcal{D}^{\immac}(X)$. Recall that $\mathcal{D}_{\sg}(R) \coloneq \mathcal{D}^b_{\fdmod}(R)/ \mathcal{D}^{\perf}(R)$ denotes the singularity category of a Noetherian ring. 

\begin{cor} \label{idempotentcomplete} Let $X$ be a locally Noetherian scheme. 
\begin{enumerate}

\item $\mathcal{D}_{\isg}(X)$ is idempotent complete if and only if $\mathcal{D}_{sg}(\widehat{\mathcal{O}}_{X,x})$ is idempotent complete for every $x \in X$. 

\item If $X$ is Gorenstein, then $\mathcal{D}_{\isg}(X)$ is idempotent complete.
\end{enumerate}
\end{cor}

\begin{proof}
By \cite[Theorem 9.1]{Sch06} (or \cite[Theorem 1.1]{AGH19}) and since $\mathcal{D}^b_{\ifp}(X)$ has a bounded $t$-structure by Theorem \ref{regularscheme}, $K_{-1}(\mathcal{D}^b_{\ifp}(X)) = 0$. Hence the same arguments as \cite[Section 4]{Orl11} apply. Therefore $\mathcal{D}_{\isg}(X)$ is idempotent complete if and only if 
\[
\bigoplus_{x \in X} K_{-1}(\widehat{\mathcal{O}}_{X,x}) \simeq K_{-1}(\mathcal{D}^{\immac}(X)) = 0
\]
where we have used Corollary \ref{regularscheme}. By \cite[Section 4]{Orl11}, $K_{-1}(\widehat{\mathcal{O}}_{X,x}) = 0$ if and only if $\mathcal{D}_{\sg}(\widehat{\mathcal{O}}_{X,x})$ is idempotent complete. If $X$ is Gorenstein, then each $\widehat{\mathcal{O}}_{X,x}$ is a Gorenstein complete local ring. By \cite[Lemma 3.1]{BK12}, $\mathcal{D}_{\sg}(\widehat{\mathcal{O}}_{X,x})$ is idempotent complete. Hence (2) follows from (1). 
\end{proof}

\begin{ex}
If $C \subseteq \mathbb{P}^2_k$ is the nodal cubic curve over a field, then $\mathcal{D}_{\sg}(C)$ is not idempotent complete. However $C$ is Gorenstein and so $\mathcal{D}_{\isg}(C)$ is idempotent complete by Corollary \ref{idempotentcomplete}.
\end{ex}

\begin{remark}
Hochschild homology is also a finitary additive invariant but this is probably not the best invariant of $\mathcal{D}^{\immac}(X)$ to consider even if $X$ is smooth and projective over a field. Indeed the uncompleted Hochschild homology of $\widehat{\mathcal{O}}_{X,x}$ is seldom well behaved. For example, one is forced to consider the non-Noetherian ring $k[[x]] \otimes_k k[[x]]$. 
\end{remark}

\section{Immaculate Equivalences}
\label{immacequivsec}

The characterisations of regularity and dimension in Theorem \ref{regularscheme} require us to keep track of an infinite set of generators. So a priori, it is not clear these characterisations depend only on the category of immaculate complexes. We now address this problem by showing that regularity and the precise dimension are invariant under equivalences between categories of immaculates. 

\begin{remark}
We recall that $\mathcal{D}^{\immac}(X)$ inherits a natural enhancement from the DG-category or stable $\infty$-category $\mathcal{D}_{\QCoh}(X)$ and we will consider our equivalences to be enhanced. The only point where we use this is to show that a functor $\mathcal{D}^{\perf}(\inj(X)^{op}) \to \mathcal{D}^{\perf}(\inj(Y)^{op})$ extends to one on the unbounded derived categories. See \cite{kel06} for a background on DG-categories. Again as an abuse of notation we use the same notation for the category of immaculates and its enhancement. 
\end{remark}

\begin{defn}
We say that schemes $X,Y$ are immaculate equivalent if there is a quasi-equivalence of DG-categories $\mathcal{D}^{\immac}(X) \simeq \mathcal{D}^{\immac}(Y)$. 
\end{defn}

\begin{remark}
If $X$ and $Y$ are schemes over a commutative ring $k$, one should ask for an equivalence of $k$-linear DG-categories. We will leave the base implicit since it does not affect any of our arguments. 
\end{remark}

The goal of this section is prove the following result.

\begin{theorem} \label{dimsame}
Suppose $X$ and $Y$ are immaculate equivalent locally Noetherian schemes. 
\begin{enumerate}
\item Then $X$ is regular if and only if $Y$ is regular.
\item In this case, $X$ has finite dimension if and only if $Y$ has finite dimension and then $\dim(X) = \dim(Y)$.
\end{enumerate} 
\end{theorem}

\begin{remark}
Theorems \ref{dimsame} and \ref{regularscheme} combined show that the minimum generation time of $\mathcal{D}^{\immac}(X)$ by the indecomposable injectives is invariant under equivalence. We find this surprising as (1) there are usually infinitely many indecomposable injectives and (2) even if one considers one-object strong generators of a triangulated category, one expects different generators to have different generation times. So (2) in particular shows that the indecomposable injectives are some special kind of generating set. If this could be made precise it may shed light on Orlov's conjecture \cite[Conjecture 10]{Orl09}. 
\end{remark}

A famous result of \cite{Ric89}, states that if $A$ and $B$ are right coherent rings, then $\mathcal{D}^{\perf}(A) \simeq \mathcal{D}^{\perf}(B)$ if and only if $\mathcal{D}^b_{\fdmod}(A) \simeq  \mathcal{D}^b_{\fdmod}(B)$. We will prove a version of this for the immaculates. A key difference in this version is that there is usually no single compact generator. However we will use the following finiteness property.

\begin{prop} \label{finitenessprop} \label{finitenesspropcor}
Let $X$ be a locally Noetherian scheme. 

\begin{enumerate}
\item If $M \in \mathcal{D}^b_{\ifp}(X)$. Then there is a finite set $\mathcal{P} 
\subseteq X$ and an injective resolution $M \xrightarrow{\sim} I$ such that if $E(x)$ appears in $I^n$ for some $n \in \mathbb{Z}$, then $x \in \mathcal{P}$. 

\item Every $M \in \mathcal{D}^b_{\fdmod}(\inj(X)^{op})$ admits a projective resolution with only finitely many isomorphism classes of indecomposable projectives appearing.
\end{enumerate}
\end{prop}

\begin{proof}
Take an injective resolution $M \xrightarrow{\sim} I$ of $M$. As the cohomology of $M$ is bounded, there is some $n \geq 0$ such that $\sigma_{\leq n}I$ is immaculate and the cohomology of $\sigma_{>n} I$ is concentrated in a single degree. Hence it remains to show that every $M \in \ifp(X)$ admits an injective resolution with only finitely many points appearing across all terms of the resolution. But this is Lemma \ref{stronginjcohlemma}. The second part then follows from Theorem \ref{tiltingequiv2}.
\end{proof}

\begin{notation} If $X$ is a locally Noetherian scheme and $x \in X$, we let $P(x) \in \fdmod(\inj(X)^{op})$ denote the indecomposable projective corepresented by $E(x)$. If $\mathcal{T}$ is a triangulated category with coproducts and $S \subseteq \mathcal{T}$, then recall that $\thick^{\geq 0}(S)$ denotes the variant of Definition \ref{buildingnotation} where only positive shifts are included. Equivalently, $\thick^{\geq 0}(S)$ is the smallest subcategory of $\mathcal{T}$ closed under extensions, summands and positive shifts. We let $\loc^{\geq 0}(S)$ be the smallest subcategory of $\mathcal{T}$ closed under extensions, positive shifts and all coproducts ($\loc^{\geq 0}(S)$ is automatically closed under summands). 

\end{notation}
\begin{lemma} \label{lemma1equiv}
Let $X,Y$ be locally Noetherian schemes. Suppose 
\[
F: \mathcal{D}(\inj(X)^{op}) \xrightarrow{\sim} \mathcal{D}(\inj(Y)^{op})
\]
is an equivalence. If $M \in \mathcal{D}^b_{\fdmod}(\inj(X)^{op})$, then $F(M) \in \mathcal{D}^-_{\fdmod}(\inj(Y)^{op})$ and there is a finite $\mathcal{Q} \subseteq Y$ such that $F(M) \in \loc(P(y) \mid y \in \mathcal{Q})$.
\end{lemma}

\begin{proof}
Let $M \in \mathcal{D}^b_{\fdmod}(\inj(X)^{op})$, then by Proposition \ref{finitenesspropcor}, we can assume $M$ is a bounded above complex of finitely generated projectives with only finitely many isomorphism classes of indecomposable projectives appearing. Let $\mathcal{P}$ denote the corresponding finite set of points of $X$. Note $F$ sends compacts to compacts and so for any $x \in X$, 
\[
F(P(x)) \in \mathcal{D}^{\perf}(\inj(Y)^{op}) = \thick(P(y) \mid y \in Y )
\]
Hence each $F(P(x))$ can be built using only finitely many corepresentable $\inj(Y)$-modules. Therefore there is a finite set $\mathcal{Q} \subseteq Y$ such that 
\[
\{ F(P(x)) \mid x \in \mathcal{P} \} \subseteq \thick(P(y) \mid y \in \mathcal{Q}).
\]
Furthermore since each $F(P(x))$ is built finitely from the $P(y)$ for $y \in \mathcal{Q}$, there exist $N_y \geq 0$ for each $y \in \mathcal{Q}$, such that 
\[
\{ F(P(x)) \mid x \in \mathcal{P} \} \subseteq \thick^{ \geq 0}(\Sigma^{-N_y} P(y) \mid y \in \mathcal{Q}).
\]
Let $N = \max(N_y \mid y \in \mathcal{Q})$. Suppose $m \in \mathbb{Z}$, then consider the brutal truncations
\[
\sigma_{> -m -1 -N} M \to M \to \sigma_{\leq -m -1 -N} M
\]
so that $\sigma_{> -m -1 -N} M \in \mathcal{D}^{\perf}(\inj(X)^{op})$. Using the triangle from Remark \ref{milnor}, we see that $\sigma_{\leq -m -1 -N} M \in \loc^{\geq 0}(\Sigma^{m+1+N} P(x) \mid x \in \mathcal{P})$. Therefore $F(\sigma_{> -m -1 -N} M) \in \mathcal{D}^{\perf}(\inj(Y)^{op})$ and 
\begin{align*}
F(\sigma_{\leq -m -1 -N} M) & \in \loc^{\geq 0}(\Sigma^{m+1+N} F(P(x)) \mid x \in \mathcal{P})  \\ 
& \subseteq \loc^{\geq 0}(\Sigma^{m+1} P(y) \mid y \in \mathcal{Q}) \\
& \subseteq \mathcal{D}^{\leq -m-1}(\inj(Y)^{op}).
\end{align*}
The vanishing on cohomology in degrees at least $-m-1$ implies the map $F(\sigma_{> -m-1-N} M) \to F(M)$ is an isomorphism on $H^i$ for $i > m$ and an epi on $H^m$. Hence $F(M)$ is $m$-pseudocoherent for every $m$ (in the sense of \cite[Tag 064N]{stacks-project}). A similar proof to \cite[Tag 0640]{stacks-project} or the dual to Lemma \ref{injfinchar} implies that $F(M)$ is a bounded above complex of finitely generated projectives. As $\inj(Y)$ is left coherent by Theorem \ref{weakdimscheme}, $F(M) \in \mathcal{D}^-_{\fdmod}(\inj(Y)^{op})$. Finally since $M \in \loc(P(x) \mid x \in \mathcal{P})$, 
\[
F(M) \in \loc(F(P(x)) \mid x \in \mathcal{P}) \subseteq \loc(P(y) \mid y \in \mathcal{Q}).\qedhere
\]
\end{proof}

\begin{lemma} \label{flatlemma}
Let $\mathcal{C}$ be a small preadditive category and $c \in \mathcal{C}$ and let $R_c := \mathcal{C}(c,c)$. Suppose $\mathcal{C}(c,b)$ is a flat right $R_c$-module for every $b \in \mathcal{C}$. Then there is an equivalence
\[
\mathcal{D}(R_c^{op}) \xrightarrow{\sim} \loc(\mathcal{C}(c,-)) \subseteq \mathcal{D}(\mathcal{C}^{op})
\]
which sends $\mathcal{D}^b(R_c)$ to $\mathcal{D}^b(\mathcal{C})$. 
\end{lemma}

\begin{proof}
Let $P$ be the  $\mathcal{C}$-$R_c$ bimodule which sends $b \mapsto \mathcal{C}(c,b)$. Then consider the functor
\[
P \otimes_{R_c} -: \Mod R_c^{op} \to \Mod \mathcal{C}^{op}
\]
where $P \otimes_{R_c} M \in \Mod \mathcal{C}^{op}$ sends $b \mapsto \mathcal{C}(c,b) \otimes_{R_c} M$. Since $\mathcal{C}(c,b)$ is a flat $R_c$ module for every $b \in \mathcal{C}$, we see that $P \otimes_{R_c} -$ is exact. The functor it induces 
\[
\mathcal{D}(R_c^{op}) \hookrightarrow \mathcal{D}(\mathcal{C}^{op})
\]
is fully faithful as it is coproduct preserving and fully faithful at $R_c$. As it sends $R_c$ to $\mathcal{C}(c,-)$ its image consists of the localising subcategory generated by $\mathcal{C}(c,-)$. Since $- \otimes_{R_c} P$ is exact it sends bounded complexes to bounded complexes. 
\end{proof}

\begin{lemma} \label{presbound}
Let $Y$ be a locally Noetherian scheme, $y \in Y$, then there is an equivalence
\[
\mathcal{D}(\widehat{\mathcal{O}}_{Y,y}) \xrightarrow{\sim} \loc(P_y) \subseteq \mathcal{D}(\inj(Y)^{op})
\]
which sends $\mathcal{D}^b(\widehat{\mathcal{O}}_{Y,y})$ to $\mathcal{D}^b(\inj(Y)^{op})$. 
\end{lemma}

\begin{proof}
We will verify the conditions of Lemma \ref{flatlemma}. Suppose $x \in Y$. By Proposition \ref{matlisprop} (1), we can assume $x \in \overline{\{y\}}$. By  \cite[II.7]{Har66} $E(x) = i_\ast E(x)$ where $i$ is the inclusion of $\Spec(\mathcal{O}_{X,x})$ into $X$. As $i_\ast$ is the right adjoint to taking stalks, 
\[
\Hom_Y(E(y),E(x)) \simeq \Hom_{\mathcal{O}_{Y,x}}(E(y),E(x))
\]
Now by the proof of \cite[Theorem 3.3.15]{EJ11} and \cite[Theorem 3.4.1 (7)]{EJ11} we have equivalences of $\mathcal{O}_{Y,y}$ modules
\[
\Hom_{\mathcal{O}_{Y,x}}(E(y),E(x)) \simeq \Hom_{\mathcal{O}_{Y,y}}\left(E(y), \bigoplus E(y)\right) \simeq \left( \bigoplus \mathcal{O}_{Y,y} \right)^{\wedge}
\]
By \cite[Tag 0AGW]{stacks-project}, this is a flat $\widehat{\mathcal{O}}_{Y,y}$ module.
\end{proof}

\begin{lemma} \label{boundedchar}
Let $\mathcal{Q} \subseteq Y$ be a finite set and ${M \in \loc(P(y) \mid y \in \mathcal{Q})} \subseteq \mathcal{D}(\inj(Y)^{op})$. Suppose for every $N \in \mathcal{D}^{\perf}(\inj(Y)^{op})$, there exists an $n \geq 0$ such that $\Hom_{\mathcal{D}(\inj(Y)^{op})}(N,\Sigma^i M) = 0$ for $\lvert i \rvert \geq n$.  Then $M \in \mathcal{D}^b(\inj(Y)^{op})$. 
\end{lemma}

\begin{proof}
We induct on $\lvert \mathcal{Q} \rvert$. Suppose $\mathcal{Q} = \{y\}$. Then $M$ is in the image of the functor 
\[
F: \mathcal{D}(\widehat{\mathcal{O}}_{Y,y}) \xrightarrow{\sim} \loc(P(y)) \subseteq \mathcal{D}(\inj(Y)^{op})
\]
from Lemma \ref{presbound}. So let $M \simeq F(L)$. By assumption 
\[
\Hom_{\mathcal{D}(\inj(Y)^{op})}(P(y),\Sigma^i M) \simeq \Hom_{\mathcal{D}(\widehat{\mathcal{O}}_{Y,y})}\left(\widehat{\mathcal{O}}_{Y,y},\Sigma^i L\right) \simeq H^i(L) = 0
\]
for $\lvert i \rvert >> 0$. By Lemma \ref{presbound}, $F$ preserves bounded complexes and so $M \in \mathcal{D}^b(\inj(Y)^{op})$. 

Suppose the claim holds for smaller sets than $\mathcal{Q}$. Since $Y \in \loc(P(y) \mid y \in \mathcal{Q})$, we can assume that $M$ is a complex of copies of $P(y)$ with $y \in \mathcal{Q}$. Let $y \in \mathcal{Q}$ be such that $y \not \in \overline{\{x\}}$ for all $y \neq x  \in \mathcal{Q}$ i.e.\ let $y$ be any of the minimal elements of $\mathcal{Q}$ with respect to the partial ordering by height. Let $M(y) \subseteq M$ consist only of copies of $P(y)$. If $y \neq x \in \mathcal{Q}$, then by Proposition \ref{matlisprop} (1) and the Yoneda Lemma,
\[
0 = \Hom_Y(E(x),E(y)) \simeq \Hom_{\inj(Y)^{op}}(P(y),P(x)).
\]
So $M(y) \subseteq M$ is a subcomplex. Let $M' \coloneq M/M(y)$, then $\Hom_{\inj(Y)^{op}}(P(y),\Sigma^i M') = 0$ for all $i$ since $M' \in \loc(P(x) \mid x \in \mathcal{Q} \setminus y)$. Therefore the long exact sequence implies there is an equivalence
\[
\Hom_{\mathcal{D}(\inj(Y)^{op})}(P(y),\Sigma^i M(y))  \simeq \Hom_{\mathcal{D}(\inj(Y)^{op})}(P(y),\Sigma^i M) 
\]
which is bounded by assumption. Therefore by the above argument $M(y) \in \mathcal{D}^b(\inj(Y)^{op})$. So there exists an $n \geq 0$ such that for any $z \in Y$, 
\[
\Hom_{\mathcal{D}(\inj(Y)^{op})}(P(z),\Sigma^i M(y)) = H^i(M(y))(z) = 0
\]
for $\lvert i \rvert \geq n$. So for any $N \in \mathcal{D}^{\perf}(\inj(Y)^{op})$ we have that $\Hom_{\mathcal{D}(\inj(Y)^{op})}(N,\Sigma^i M(y))$ is bounded. By assumption the same holds for $M$ and so by the long exact sequence we see that $\Hom_{\mathcal{D}(\inj(Y)^{op})}(N,\Sigma^i M') = 0$ for $\lvert i \rvert >> 0$ for any $N \in \mathcal{D}^{\perf}(\inj(Y)^{op})$. Therefore by induction $M' \in \mathcal{D}^b(\inj(Y)^{op})$. So since $M(y), M' \in \mathcal{D}^b(\inj(Y)^{op})$ we are done. 
\end{proof}

\begin{theorem} \label{determineschemes}
Let $X,Y$ be locally Noetherian schemes. If they are immaculate equivalent, then there is a diagram of equivalences
\[
\begin{tikzcd}
\mathcal{D}^b_{\ifp}(X) \arrow[r,"\sim"] &\mathcal{D}^b_{\ifp}(Y) \\
\mathcal{D}^{\immac}(X) \arrow[u,hook] \arrow[r,"\sim"] & \mathcal{D}^{\immac}(Y) \arrow[u,hook]
\end{tikzcd}
\]
\end{theorem}

\begin{proof} Suppose $X,Y$ are immaculate equivalent. By Theorem \ref{tiltingequiv1}, there is an enhanced equivalence 
\[
\mathcal{D}^{\perf}(\inj(X)^{op}) \simeq \mathcal{D}^{\immac}(X)^{op} \simeq \mathcal{D}^{\immac}(Y)^{op} \simeq \mathcal{D}^{\perf}(\inj(Y)^{op})
\]
which extends to an equivalence $F: \mathcal{D}(\inj(X)^{op}) \xrightarrow{\sim} \mathcal{D}(\inj(Y)^{op})$ on the Ind-completions. We claim this restricts to an equivalence $\mathcal{D}^b_{\fdmod}(\inj(X)^{op}) \xrightarrow{\sim} \mathcal{D}^b_{\fdmod}(\inj(Y)^{op})$. If $M \in \mathcal{D}^b_{\fdmod}(\inj(X)^{op})$, then by Lemma \ref{lemma1equiv}, $F(M) \in \mathcal{D}^-_{\fdmod}(\inj(X)^{op})$ and there is a finite $\mathcal{Q} \subseteq Y$ such that $F(M) \in \loc(P(y) \mid y \in \mathcal{Q})$. Since $M \in \mathcal{D}^b(\inj(X)^{op})$, we have that 
\[
\Hom_{\mathcal{D}(\inj(X)^{op})}(P(x), \Sigma^i M) \simeq H^i(M)(x) = 0
\]
for $\lvert i \rvert >> 0$ and every $x \in X$. Therefore $\Hom_{\mathcal{D}(\inj(X))^{op}}(N,\Sigma^i M)$ is bounded for every $N \in \mathcal{D}^{\perf}(\inj(X)^{op})$. As $F$ restricts to an equivalence on perfects, we have that $\Hom_{\mathcal{D}(\inj(Y)^{op})}(N,\Sigma^i F(M))$ is bounded for every $N \in \mathcal{D}^{\perf}(\inj(Y)^{op})$. Therefore by Lemma \ref{boundedchar}, $F(M) \in \mathcal{D}^b_{\fdmod}(\inj(Y)^{op})$. By symmetry, we see that the inverse of $F$ sends $\mathcal{D}^b_{\fdmod}(\inj(Y)^{op})$ to $\mathcal{D}^b_{\fdmod}(\inj(X)^{op})$. Therefore $F$ restricts to an enhanced equivalence $\mathcal{D}^b_{\fdmod}(\inj(Y)^{op}) \simeq  \mathcal{D}^b_{\fdmod}(\inj(X)^{op})$. By Theorem \ref{tiltingequiv2}, there are equivalences
\[
\mathcal{D}^b_{\ifp}(X) \simeq \mathcal{D}^b_{\fdmod}(\inj(Y)^{op})^{op} \simeq  \mathcal{D}^b_{\fdmod}(\inj(X)^{op})^{op} \simeq \mathcal{D}^b_{\ifp}(Y)
\]
making the above diagram commute.
\end{proof}

The converse will follow from a few lemmas.

\begin{lemma} \label{hmin}
Let $X$ be a locally Noetherian scheme,  $\mathcal{P} \subseteq X$ a specialisation closed subset and $I$ a bounded below homotopically minimal complex of injectives in $\QCoh(X)$. Let $I(\mathcal{P}) \subseteq I$ denote the subcomplex consisting of copies of $E(x)$ with $x \in \mathcal{P}$. Then $I(\mathcal{P})$ is also homotopically minimal. 
\end{lemma}

\begin{proof}
The differential restricts to $d': I(\mathcal{P})^n \to I(\mathcal{P})^{n+1}$ by Proposition \ref{matlisprop} (1) and since $\mathcal{P}$ is specialisation closed. Since $\ker(d)$ is essential in $I^n$, it is clear that $\ker(d') = \ker(d) \cap I(\mathcal{P})^n$ is essential in $I(\mathcal{P})^n$. Hence $I(\mathcal{P})^n = E(\ker(d'))$ as required. 
\end{proof}

\begin{lemma} \label{basslemma}
Suppose $X$ is a locally Noetherian scheme, $x \in X$ and $I$ a bounded below homotopically minimal complex of injectives such that for each $n$, $I^n = \bigoplus_{x \in \overline{\{y\}}} E(y)^{\mu_n(y)}$ for some cardinalities $\mu_n(y)$. Then $\Hom_{\mathcal{D}_{\QCoh(X)}}(k(x), \Sigma^n I) = k(x)^{\mu_n(y)}$ 
\end{lemma}

\begin{proof}
Consider the subcomplex $I(x) \subseteq I$ consisting of injectives with only $x$ appearing. Then by Proposition \ref{matlisprop} (1), $\Hom_X(k(x),I/I(x)) = 0$ and so  for each $n$ $\Hom_{\mathcal{D}_{\QCoh}(X)}(k(x),\Sigma^n I/I(x))$ also vanishes. It follows that 
\[
\Hom_{\mathcal{D}_{\QCoh}(X)}(k(x), \Sigma^n I) = \Hom_{\mathcal{D}_{\QCoh(X)}}(k(x),\Sigma^n I(x)) 
\]
Consider the fully faithful functor $i_\ast: \Mod \mathcal{O}_{X,x} \hookrightarrow \QCoh(X)$ which is the exact right adjoint to taking stalks. Then $I(x) = i_\ast I(x)_x = \mathbb{R} i_\ast I(x)_x$. Therefore,
\[
\Hom_{\mathcal{D}_{\QCoh}(X)}(k(x),\Sigma^n I(x)) \simeq \Hom_{\mathcal{D}(\mathcal{O}_{X,x})}(k(x),\Sigma^n I(x)) \simeq \Hom_{\mathcal{O}_{X,x}}(k(x), E(x)^{\mu_n(x)}).
\]
The second equivalence follows from Lemmas \ref{hmin} and \ref{minimallemma}. Now the claim follows from Proposition \ref{matlisprop} (3).
\end{proof}

\begin{lemma} \label{conversedetermine}
Let $X$ be a locally Noetherian scheme, then $M \in \mathcal{D}^b_{\ifp}(X)$ is immaculate if and only if for every $ N \in \mathcal{D}^b_{\ifp}(X)$ there is some $n \geq 0$ such that  
\[
\Hom_{\mathcal{D}_{\QCoh}(X)}(N,\Sigma^i M) = 0 \text{ for } \lvert i \rvert \geq n. 
\]
\end{lemma}

\begin{proof}
Suppose that $M$ satisfies the condition. By Proposition \ref{finitenessprop} (1), we can assume that $M$ is a bounded below complex of finite sums of indecomposable injectives in which only finitely many points $\mathcal{P} \subseteq X$ appear. By \cite[Proposition B.2]{kra05}, we can assume $M$ is also homotopically minimal. By assumption, there exists an $n \geq 0$ such that $\Hom_{\mathcal{D}_{\QCoh}(X)}(k(x),\Sigma^i M) = 0$  for  $\lvert i \rvert \geq n$ for all $x \in \mathcal{P}$. We induct on the size of $\mathcal{P}$. Let $M^i = \bigoplus E(x)^{\mu_i(x)}$. If $x$ is maximal in $\mathcal{P}$, then every point appearing in $I_x$ specialises to $x$ and Lemma \ref{basslemma} applies. Hence $\mu_i(x) = 0$ for $i > n$. So $I(x) \subseteq I$, the subcomplex consisting of copies of $E(x)$, is immaculate. Now consider $I/I(x)$. Since $I$ and $I(x)$ are in $\mathcal{D}^b_{\ifp}(X)$ so is $I/I(x)$. As $I(x)$ is immaculate the vanishing statement holds for $I(x)$ and $I$. Hence the vanishing holds for $I/I(x)$. Now by induction $I/I(x)$ is immaculate. Hence $I$ is immaculate. 
\end{proof}

\begin{remark}
Note the condition in Lemma \ref{conversedetermine} is a priori weaker than the condition in Lemma \ref{immacchar} since the vanishing bound is not required to uniform even when restricted to objects of $\ifp(X)$.  
\end{remark}

To prove (2) of Theorem \ref{dimsame}, we use the following result. 

\begin{lemma} \label{dimensionlemma}
Let $X$ be a regular scheme of dimension $n$ and $\mathcal{Q} \subseteq X$ any subset. Then $\thick(E(x) \mid x \in \mathcal{Q}) = \thick_n(E(x) \mid x \in \mathcal{Q})$. 
\end{lemma}

\begin{proof} 
Let $M \in \thick(E(x) \mid x \in \mathcal{Q})$, then we may assume that $M$ is a bounded complex of finite sums of copies of $E(x)$ with $x \in \mathcal{Q}$. Let $K = \coker(d^{i-1}: M^{i-1} \to M^{i})$. Then by Lemma \ref{stronginjcohlemma}, $E(K)$ is a finite sum of copies of $E(x)$ with $x \in \mathcal{Q}$. Then we have
\[
H^i(M) = \ker(d^i)/\im(d^{i-1}) \subseteq M^i/\im(d^{i-1}) =K \subseteq E(K)
\]
Therefore $E(H^i(M))$ is a summand of $E(K)$ and so also a finite sum of copies of $E(x)$ with $x \in \mathcal{Q}$. We now proceed as in Lemma \ref{smallglobdimtostronggen} by inducting on the injective dimension of $H^\ast(M) \leq n$. If $H^\ast(M)$ is injective then $M \simeq H^\ast(M) \simeq E(H^\ast(M))$ is a finite sum of copies of $E(x)$ with $x \in \mathcal{Q}$. Otherwise, we lift the maps $H^i(M) \to E(H^i(M))$ to a map $M \to \bigoplus \Sigma^{i} E(H^i(M))$ where the latter is a finite sum of copies of $E(x)$ with $x \in Q$. Let $C$ denote the cone of this map so $C \in \thick(E(x) \mid x \in \mathcal{Q})$. Then the long exact sequence on cohomology implies that $\idim H^\ast(C) \leq \idim H^\ast(M) - 1$. So by induction $C \in \thick_{ \idim M -1}(E(x) \mid x \in \mathcal{Q})$. Hence $M \in \thick_{\idim H^\ast(M)}(E(x) \mid x \in \mathcal{Q})$. Since $\idim H^\ast(M) \leq n$ we are done. 
\end{proof}

\begin{proof}[Proof of Theorem \ref{dimsame}]
By Theorem \ref{determineschemes} and Lemma \ref{conversedetermine}, we see that $X$ and $Y$ are immaculate equivalent if and only if there is an enhanced equivalence $\mathcal{D}^b_{\ifp}(X) \simeq \mathcal{D}^b_{\ifp}(Y)$ which is compatible with the inclusions of the immaculates. So by Theorem \ref{regularscheme} (3), $X$ is regular if and only if $Y$ is regular. 

Suppose $X$ is regular and finite-dimensional. We have an equivalence $F: \mathcal{D}^{\immac}(X) \xrightarrow{\sim} \mathcal{D}^{\immac}(Y)$ with inverse $G$. Let $y \in Y$ be a closed point and $n \coloneq \dim(X)$. Then $G(E(y)) \in \thick(E(x) \mid x \in \mathcal{Q})$ for some finite $\mathcal{Q} \subseteq X$. Consider the composition below where the last map is the localisation functor of Remark \ref{closedlocalisation}.
\[
H: \thick(E(x) \mid x \in \mathcal{Q}) \hookrightarrow \mathcal{D}^{\immac}(X) \xrightarrow{\sim} \mathcal{D}^{\immac}(Y) \twoheadrightarrow \thick(E(y)) \simeq \mathcal{D}^{\perf}(\widehat{\mathcal{O}}_{Y,y})
\]
Note $H$ is essentially surjective since $E(y) \simeq FG(E(y)) \in \thick(F(E(x)) \mid x \in \mathcal{Q})$. By Lemma \ref{dimensionlemma}, $\thick(E(x) \mid x \in \mathcal{Q}) = \thick_n(E(x) \mid x \in \mathcal{Q})$. Therefore as $H$ is essentially surjective we have that $\thick_n(H(\bigoplus_{x \in \mathcal{Q}} E(x))) = \mathcal{D}^{\perf}(\widehat{\mathcal{O}}_{Y,y})$. As $\mathcal{D}^{\perf}(\widehat{\mathcal{O}}_{Y,y})$ admits a strong generator, $\widehat{\mathcal{O}}_{Y,y}$ is regular and of finite dimension by \cite[Proposition 7.5]{Rou08}.  By \cite[Theorem A]{Let25}, we see that $n \geq \dim(\widehat{\mathcal{O}}_{Y,y}) = \dim(\mathcal{O}_{Y,y})$. Therefore every closed point of $Y$ is regular and has dimension at most $\dim(X)$ and so $\dim(Y) \leq \dim(X)$. But then $Y$ is regular and finite-dimensional and so by symmetry we see that $\dim(X) \leq \dim(Y)$. 
\end{proof}

\begin{remark} We make a couple of remarks on the proof of Theorem \ref{dimsame}.
\begin{enumerate}
\item We needed to use the full power of \cite[Theorem A]{Let25} and not the earlier versions of 
\cite{Rou08} which have finite type assumptions. 

\item The proof that $X$ being regular and finite dimensional implies that $Y$ is regular and $\dim(Y) = \dim(X)$ does not use that the equivalence $\mathcal{D}^{\immac}(X) \simeq \mathcal{D}^{\immac}(Y)$ is enhanced. 
\end{enumerate}
\end{remark}

\section{Reconstruction for Affine Schemes}
\label{reconstructionsec}

In this section we show that a commutative Noetherian ring can be reconstructed from its category of immaculate complexes. Recall that the centre $Z(\mathcal{C})$ of a preadditive category $\mathcal{C}$ is the ring of natural endotransformations of the identity functor. The following result is surely well known. For example a similar statement appears in \cite[Theorem 6.9]{AS22}. We include a proof for convenience.

\begin{lemma} \label{centreiscentre}
Let $\mathcal{A}$ be any Matlis category. Then
\[
Z(\mathcal{A}) \simeq Z(\Inj(\mathcal{A})) \simeq Z(\inj(\mathcal{A})).
\]
\end{lemma}

\begin{proof}
Restriction produces ring maps $Z(\mathcal{A}) \to Z(\Inj) \to Z(\inj)$. If $I \in \Inj$, then $I \simeq \bigoplus E_i$ for $E_i \in \inj$ and naturality implies that for any $\alpha \in Z(\Inj)$, $\alpha_{I} = \bigoplus \alpha_{E_i}$. So the map $Z(\Inj) \to Z(\inj)$ is mono. If $\beta \in Z(\inj)$, and $I \in \Inj$, then choose any decomposition $I \simeq \bigoplus E_i$ into indecomposables and set $\alpha_I \coloneq \bigoplus \beta_{E_i}$. Uniqueness of the decomposition implies $\alpha_I$ is well defined. Given a map of injectives $f: I = \bigoplus E_i \to J = \bigoplus E_j$, then we must check $\alpha_J f = f \alpha_J$. As $\mathcal{A}$ is a Grothendieck category, it is enough to check the two maps $E_i \to I \to J \to E_j$ agree when composed with restriction and projection for all $i$ and $j$. But this follows by naturality of $\beta$. Therefore $Z(\Inj) \simeq Z(\inj)$. If $M \in \mathcal{A}$, then it admits an injective presentation
\[
0 \to M \to I^0 \to I^1.
\]
If $\alpha \in Z(\mathcal{A})$ such that $\alpha_{I^0} = 0$, then $\alpha_{M} = 0$. Given $\beta \in Z(\Inj)$, we see that $\beta_{I^0}$ restricts to a map $\alpha_M: M \to M$. Again one checks this doesn't depend on the presentation and produces a natural transformation $\beta \in Z(\mathcal{A})$. Therefore $Z(\mathcal{A}) \simeq Z(\Inj)$. 
\end{proof}

\begin{remark} Suppose $\mathcal{D}$ is a DG-category over a commutative ring $k$. The Hochschild cohomology of $\mathcal{D}$ (relative to $k$), denoted $HH^\ast_k(\mathcal{D})$ is an invariant generalising the usual notion for $k$-algebras. See \cite{LVdB05} for more background. By definition $HH^\ast_k(\mathcal{D})$ is the cohomology of the Hochschild cochain complex of a $k$-cofibrant resolution of $\mathcal{D}$. 
\end{remark}

If $k$ is a commutative ring, recall that a $k$-linear category is a category enriched in $\Mod k$. For example a preadditive category is a $\mathbb{Z}$-linear category. Any $k$-linear category can be viewed as a DG-category over $k$ with morphisms concentrated in degree zero. If $\mathcal{A}$ is a $k$-linear Grothendieck category, then $\mathcal{D}(\mathcal{A})$ and its subcategories are $k$-linear DG-categories.

\begin{theorem} \label{centreinvariant}
Let $k$ be a commutative ring and $\mathcal{A}$ a $k$-linear Matlis category. Then 
\[
HH^0_k(\mathcal{D}^{\immac}(\mathcal{A})) \simeq Z(\mathcal{A}). 
\]
\end{theorem}

\begin{proof}
Since $\mathcal{D}^{\immac}(\mathcal{A}) \simeq \mathcal{D}^{\perf}(\inj(\mathcal{A}))$ and Hochschild cohomology is Morita invariant, we have that $HH^0_k(\mathcal{D}^{\immac}(\mathcal{A})) \simeq HH^0_k(\inj(\mathcal{A}))$. The zero-th Hochschild cohomology of any $k$-linear category coincides with its centre. This is well known if $k$ is a field, and otherwise follows from a similar argument to \cite[pp.\ 182]{Shuk61} where it is shown for $k$-algebras i.e.\ one object $k$-linear categories. Then the claim follows from Lemma \ref{centreiscentre}.
\end{proof}

\begin{cor} \label{recoveraffine}
Suppose $A$ and $B$ are right Noetherian rings and there is a quasi-equivalence of DG-categories $\mathcal{D}^{\immac}(A) \simeq \mathcal{D}^{\immac}(B)$. Then their centres are isomorphic $Z(A) \cong Z(B)$. In particular if $R,S$ are commutative Noetherian rings which are immaculate equivalent then $R \cong S$.
\end{cor}

\begin{proof}
As $A$ and $B$ are right Noetherian their module categories are Matlis categories. Hochschild cohomology is invariant under quasi-equivalence of DG-categories and so by Theorem \ref{centreinvariant} $Z(\Mod A) \simeq Z(\Mod B)$. It is well known that the centre of a ring coincides with the centre of its module category.
\end{proof}

\begin{remark}
In \cite{Ric89}, it was shown that derived equivalent rings have isomorphic centres. Theorem \ref{recoveraffine} is an immaculate version of this result.  
\end{remark}

\begin{remark}
We can also make statements in the non-affine case. For simplicity we work over a field $k$. The Hochschild cohomology of a separated $k$-scheme $X$ is defined as $HH^n(X) \coloneq \Ext^n_{X \times_k X}(\Delta_* \mathcal{O}_X,\Delta_*\mathcal{O}_X)$ where $\Delta: X \to X \times_k X$ is the diagonal. 
\end{remark}

\begin{theorem} \label{H^0invariant} Let $k$ be a field and suppose $X, Y$ are immaculate equivalent Noetherian separated $k$-schemes. Then $H^0(X,\mathcal{O}_X) \simeq H^0(Y,\mathcal{O}_Y)$. 
\end{theorem}

\begin{proof}
Ignoring some size issues we have by Lemma  \ref{centreiscentre}  that
\[
HH^0_k(\Inj(X)) \simeq Z(\Inj(X))) \simeq Z(\inj(X)) \simeq HH^0_k(\mathcal{D}^{\immac}(X))
\]
By \cite[Corollary 7.8.2]{LVdB05}, we have that $HH^\ast_k(\Inj(X)) \simeq HH^\ast_k(X)$. Now we note that
\[
HH^0(X) \simeq \Hom_{\mathcal{D}_{\QCoh}(X)}(\mathbb{L} \Delta^\ast \Delta_\ast \mathcal{O}_X, \mathcal{O}_X) \simeq \Hom_{\mathcal{D}_{\QCoh}(X)}(\mathcal{O}_X, \mathcal{O}_X) \simeq H^0(X,\mathcal{O}_X)
\]
where the second equivalence uses the standard $t$-structure and the fact that $\mathcal{H}^0(\mathbb{L} \Delta^\ast \Delta_\ast \mathcal{O}_X) \simeq \Delta^\ast \Delta_\ast \mathcal{O}_X \simeq \mathcal{O}_X$. 
\end{proof}

These results and the abundance of Hochschild cohomology isomorphisms (for example \cite[Theorem 4.4.1]{LVdB05}, \cite[Proposition 2.4]{kel19}) suggests the following conjecture. 

\begin{conj} \label{HHconj}
Let $X$ be a separated scheme over a field $k$. Then $HH_k^\ast(X) \simeq HH^\ast_k(\inj(X))$
\end{conj}

\begin{remark}
\begin{enumerate}
\item By \cite[Corollary 7.8.2]{LVdB05},  $HH_k^\ast(X) \simeq HH^\ast_k(\Inj(X))$ where $\Inj(X)$ is viewed as a $k$-linear category. So to produce an isomorphism on all Hochschild cohomology it is enough to show the inclusion $\inj(X) \hookrightarrow \Inj(X)$ induces an isomorphism on Hochschild cohomologies of the $k$-linear categories.

\item If Conjecture \ref{HHconj} holds, then the semi-orthogonal decompositions of Section \ref{SODsec} and the long exact sequences of \cite[Section 4.5]{kel03} would give a relationship between the Hochschild cohomology of $X$ and of those of the completions of its local rings. 

\end{enumerate} 
\end{remark}

It is natural to ask about the non-affine reconstruction problem (motivated by the results of \cite{BO01} for the perfects). 

\begin{quest} \label{recoverscheme}
What are some conditions on schemes $X,Y$ such that $\mathcal{D}^{\immac}(X) \simeq \mathcal{D}^{\immac}(Y)$ implies $X \simeq Y$?
\end{quest}

\begin{remark} We finish with some speculative remarks.
\begin{enumerate}
\item Using Corollary \ref{K-theory} (1), we can always recover the underlying set of $X$ from $\mathcal{D}^{\immac}(X)$ and so  immaculate equivalent schemes have the same cardinality.

\item Theorem \ref{H^0invariant} implies a proper scheme cannot be immaculate equivalent to a non-Artinian affine scheme. We don't know if being proper over a field can be detected from $\mathcal{D}^{\immac}(X)$. 

\item One can show using Theorem \ref{regularscheme} and since there are only finitely many points, that any locally Noetherian scheme immaculate equivalent to a discrete valuation ring $R$ must in fact be isomorphic to $\Spec(R)$. 

\item If $X$ and $Y$ are Noetherian schemes such that $\inj(X) \simeq \inj(Y)$, then they have the same underlying topological spaces. This was proved in \cite[VI.3]{Gab62}.

\item The category $\inj(R)$ seems somehow related to the chromatic fracture squares of \cite{GB20}. Here $\mathcal{D}^{\perf}(X)$ is shown to be a homotopy limit of the derived categories of the completions of its local rings. Although, it is not clear how to make this relationship to the immaculates precise. 

\item It was shown in \cite[Theorem 6]{Roos69} that one can recover $X$ from $\End_X(\bigoplus E(x))$ if the latter ring is given a certain topology. This is also clearly a similar object to $\inj(X)$. 

\item We don't expect $\mathcal{D}^{\immac}(X)$ to determine $X$ completely but unfortunately we don't have any examples of immaculate equivalences to disprove this. Although one could interpret Theorem \ref{dimsame} as some (extremely weak) evidence in this direction.
\end{enumerate}
\end{remark}

\section{Examples: Curves} \label{curvessec}

In this section we completely describe the  preadditive category $\inj(X)$ when $X$ is a one-dimensional locally Noetherian scheme which is either Gorenstein or integral. By Theorem \ref{matlisthy} and Proposition \ref{matlisprop}, this amounts to describing $\Hom_X(E(x),E(y))$ as a $\widehat{\mathcal{O}}_{X,y}-\widehat{\mathcal{O}}_{X,x}$ bimodule for all points $x,y \in X$ with $y \in \overline{\{x\}}$.

\begin{lemma} \label{gorensteinhom}
Let $R$ be a Gorenstein ring of dimension one, $p$ a minimal prime and $m$ a maximal ideal with $p \subset m$. Then 
\[
\Hom_R(E(R/p),E(R/m)) \simeq \widehat{R}_m \otimes_{R_m}  R_p. 
\]
\end{lemma}

\begin{proof}
As $R$ is Gorenstein, the minimal injective resolution of $R_m$ is 
\[
0 \to R_m \to \bigoplus_{q \subset m} E(R/q) \to E(R/m) \to 0
\]
by e.g.\ \cite[Theorem 9.2.27]{EJ11}.
Applying $\Hom_R(-,E(R/m))$ and Proposition \ref{matlisprop} gives an exact sequence
\[
0 \to \widehat{R}_m \to \bigoplus_{q \subset m} \Hom_R(E(R/q),E(R/m)) \to E(R/m) \to 0.
\]
Since $\Hom_R(E(R/q),E(R/m))$ is an $R_q$-module, we have that $\Hom_R(E(R/q),E(R/m))_p = 0$ for all minimal primes $q \neq p$ and $\Hom_R(E(R/p),E(R/m))_p \simeq \Hom_R(E(R/p),E(R/m))$. Then since $E(R/m)$ is supported only at $m$, localising the above sequence at $p$ gives an isomorphism
\[
\Hom_R(E(R/p),E(R/m)) \simeq \widehat{R}_m \otimes_R R_p \simeq \widehat{R}_m \otimes_{R_m} R_p.
\]
\end{proof}

\begin{theorem} \label{gorensteineg}
If $X$ is a Gorenstein scheme of dimension $1$, then $\inj(X)$ has the form
\[
\begin{tikzcd}[column sep = huge]
x \arrow[r,"\widehat{\mathcal{O}}_{X,y} \otimes_{\mathcal{O}_{X,y}} \mathcal{O}_{X,x}"] \arrow[loop left,"\mathcal{O}_{X,x}"] & y \arrow[loop right, "\widehat{\mathcal{O}}_{X,y}"]
\end{tikzcd}
\]
where $x$ varies over all height zero points and $y$ all closed points. By convention we set $\widehat{\mathcal{O}}_{X,y} \otimes_{\mathcal{O}_{X,y}} \mathcal{O}_{X,x} = 0$ if $y \not \in \overline{\{x\}}$.
\end{theorem}

\begin{proof}
By Theorem \ref{matlisthy} and Proposition \ref{matlisprop}, the objects of $\inj(X)$ correspond to those described and the only morphisms are endomorphisms and morphisms $\Hom_X(E(x),E(y))$ where $x$ specialises to $y$. Since $E(x)$ and $E(y)$ lie in the image of the fully faithful functor induced by $i: \Spec(\mathcal{O}_{X,y}) \to X$, we can reduce to working over $\mathcal{O}_{X,y}$. Hence the result follows from Lemma \ref{gorensteinhom}.
\end{proof}

In the case of integral schemes we can use a result of Schenzel to describe the morphisms between indecomposable injectives.

\begin{theorem} \label{integraleg}
Let $X$ be an integral scheme of dimension 1 with generic point $\eta$ and let $K(X) \coloneq \mathcal{O}_{X,\eta}$. For $y \in X$ a closed point, let $q^{y}_i$ denote the associated primes of $\widehat{\mathcal{O}}_{X,y}$. Then $\inj(X)$ has the form 
\[
\begin{tikzcd}[column sep = huge]
\eta \arrow[r,"{\bigoplus_i E_{\widehat{\mathcal{O}}_{X,y}}\left(\widehat{\mathcal{O}}_{X,y}/q^y_i\right)}"] \arrow[loop left,"K(X)"] & y \arrow[loop right, "\widehat{\mathcal{O}}_{X,y}"]
\end{tikzcd}
\]
where $y$ varies over closed points.
\end{theorem}

\begin{proof}
As in Theorem \ref{gorensteineg}, we can reduce to the case $X = \Spec(R)$ where $R$ is a local one-dimensional domain. Let $q_i$ denote the associated primes of the completion of $R$ at the maximal ideal, $\widehat{R}$. Then by \cite[Theorem 2.5]{Sch15}
\[
\Hom_R(E(R),E(R/m)) \simeq \bigoplus E_{\widehat{R}}(\widehat{R}/q_i). \qedhere
\]
\end{proof}

\begin{remark}
We can simplify the descriptions in Theorems \ref{gorensteineg} and \ref{integraleg} under a mild assumption. Recall a locally Noetherian scheme $X$ is analytically unramified if $\widehat{\mathcal{O}}_{X,x}$ is reduced for each $x \in X$. This implies that $X$ is reduced and includes all naturally occurring examples of reduced schemes. For example, any reduced quasi-excellent scheme is analytically unramified. We let $Q(R)$ denote the total ring of fractions of a ring $R$. If $R$ is an analytically unramified local ring and $p$ is a minimal prime of $R$, then we let
\[
Q_p(\widehat{R}) = \prod_{q \in \min(\widehat{R})  \ q \cap R = p} k(q) 
\]
where $\min(\widehat{R})$ consists of the minimal $\widehat{R}$ primes. The total ring of fractions decomposes as $Q(\widehat{R}) \simeq \prod_{p \in \min(R)} Q_p(\widehat{R})$.
\end{remark}

\begin{lemma} \label{unramgorensteinlem}
Let $R$ be an analytically unramified Noetherian local ring, $p$ a minimal prime of $R$. Then 
\[
k(p) \otimes_R \widehat{R} \simeq Q_p(\widehat{R}).
\] 
\end{lemma}

\begin{proof}
Since $\widehat{R}$ is reduced there is an embedding $\widehat{R} \hookrightarrow \prod_{q \in \min(\widehat{R})} k(q)$
As $k(p) = R_p$ is flat this produces an embedding
\[
\widehat{R} \otimes_R k(p) \hookrightarrow \prod_{q \in \min(\widehat{R})} k(q) \otimes_R k(p) \simeq \prod_{q \in \min(\widehat{R}) \ q \cap R = p} k(q)
\]
To show the above map is also epi it is enough to show that the maps $\widehat{R} \otimes_R k(p) \to k(q)$ are all non-zero. But if one of these compositions were zero, then the composition $\widehat{R} \to \prod k(q) \to k(q)$ would be zero. This cannot occur since $\widehat{R}$ is essential in $\prod_{q \in \min(\widehat{R})} k(q) = Q(\widehat{R})$.
\end{proof}

\begin{prop} \label{unramgorielem}
Let $X$ be an analytically unramified one-dimensional Gorenstein scheme. Then, $\inj(X)$ has the form
\[
\begin{tikzcd}[column sep = huge]
x \arrow[r,"{Q_x\left(\widehat{\mathcal{O}}_{X,y}\right)}"] \arrow[loop left,"\mathcal{O}_{X,x}"] & y \arrow[loop right, "\widehat{\mathcal{O}}_{X,y}"]
\end{tikzcd}
\]
where $y$ varies over closed points and $x$ over all height zero points and by convention $Q_x(\widehat{\mathcal{O}}_{X,y}) = 0$ if $y \not \in \overline{\{x\}}$.
\end{prop}

\begin{proof}
This follows from Theorem \ref{gorensteineg} and \ref{unramgorensteinlem}. 
\end{proof}

\begin{prop} \label{unramintegralprop}
Let $X$ be an analytically unramified one-dimensional integral locally Noetherian scheme. Then $\inj(X)$ has the form
\[
\begin{tikzcd}[column sep = huge]
\eta \arrow[r,"{Q\left(\widehat{\mathcal{O}}_{X,y}\right)}"] \arrow[loop left,"K(X)"] & y \arrow[loop right, "\widehat{\mathcal{O}}_{X,y}"]
\end{tikzcd}
\]
where $\eta$ is the generic point and $y$ varies over all closed points. 
\end{prop}

\begin{proof}
We use Theorem \ref{integraleg}. As $X$ is analytically unramified, $\widehat{\mathcal{O}}_{X,y}$ is reduced and so its associated primes are all minimal. If $R$ is a reduced Noetherian ring then $E(R) \simeq Q(R) \simeq \prod_{p \in \min(R)} k(p)$ is the injective hull of $R$. So for any minimal prime $E(R/q) \simeq E(k(q)) \simeq k(q)$. Therefore $\bigoplus E(\widehat{\mathcal{O}}_{X,y}/q_i^y) \simeq Q(\widehat{\mathcal{O}}_{X,y})$. 
\end{proof}

\begin{ex} \label{injegslist}
We give a few examples to show one can very explicitly describe $\inj(X)$. Let $\Frac(R)$ denote the field of fractions of a domain and $k$ an algebraically closed field. 
\begin{enumerate}
\item If $R$ is a Dedekind domain, $\inj(R)$ has the form
\[
\begin{tikzcd}[column sep = huge]
0 \arrow[r,"\Frac(\widehat{R}_m)  "] \arrow[loop left,"\Frac(R)"] & m \arrow[loop right, "\widehat{R}_m"]
\end{tikzcd}
\]
where $m$ varies over all height one primes. 

\item In particular, $\inj(\mathbb{A}^1_k)$ has the form
\[
\begin{tikzcd}[column sep = huge]
\eta \arrow[r,"{k((t- \lambda))}"] \arrow[loop left,"k(t)"] & \lambda \arrow[loop right, "{k[[t-\lambda]]}"] 
\end{tikzcd}.
\]

\item Let $C$ be a smooth one-dimensional scheme over $k$. Then $\inj(C)$ has the form 
\[
\begin{tikzcd}[column sep = huge]
\eta \arrow[r,"{k((t_x))}"] \arrow[loop left,"K(C)"] & x \arrow[loop right, "{k[[t_x]]}"]
\end{tikzcd}
\]
where $x$ varies over closed points and $t_x$ are uniformizers of $\widehat{\mathcal{O}}_{X,x}$. 

\item In particular, $\inj(\mathbb{P}^1_k)$ is given by
\[
\begin{tikzcd}[column sep = huge]
\eta \arrow[r,"{k((t- \lambda))}"] \arrow[loop left,"k(t)"] \arrow[dr,"{k((t^{-1}))}",swap] & \lambda \arrow[loop right, "{k[[t-\lambda]]}"]  \\
& \infty \arrow[loop right, "{k[[t^{-1}]]}"] 
\end{tikzcd}
\]
where $\lambda \in k$.
\end{enumerate}
\end{ex}

We now consider the case where $\inj(X)$ has finitely many objects and so is Morita equivalent to a ring. For simplicity we will assume $X$ is the spectrum of a local one-dimensional ring but this assumption could be replaced by any scheme with finitely many points. 

\begin{prop} \label{oneobjegs}
Suppose $R$ is a local one-dimensional Gorenstein ring. Let $R_{\min} = \prod_{ p \in \min{R}} R_p$. Define the ring
\[
A \coloneq \begin{pmatrix}
R_{\min} & 0 \\
\widehat{R} \otimes_R R_{\min} & \widehat{R} 
\end{pmatrix}
\]
Then $\Mod \inj(R) \simeq \Mod A$. If $R$ is analytically unramified, then 
\[
A \simeq \begin{pmatrix}
Q(R) & 0 \\
Q(\widehat{R}) & \widehat{R} 
\end{pmatrix}
\]
Let $S$ be a one-dimensional Noetherian analytically unramified domain. 
Then $\Mod \inj(S) \simeq \Mod B$ where
\[
B \coloneq \begin{pmatrix}
\Frac(S) & 0 \\
Q(\widehat{S}) & \widehat{S} 
\end{pmatrix}
\]

\end{prop}

\begin{proof}
Since $\Spec(R)$ is finite, $\inj(R)$ has finitely many objects and $\Mod \inj(R)$ has a compact projective generator given by the direct sum of the representables. Hence $\Mod \inj(R) \simeq \Mod A$ where $A = \End_R(\bigoplus E(R/p))$. The description of $A$ then follows from Theorem \ref{gorensteineg} and Proposition \ref{unramgorielem}. The result for $S$ is similar and follows from Proposition \ref{unramintegralprop}.
\end{proof}

We can verify some other properties in these small examples.

\begin{prop} \label{notnoetheg}
Let $R$ be a one-dimensional Noetherian local ring. Suppose $R$ is either Gorenstein or an analytically unramified domain. Then
\begin{enumerate}
\item $\inj(R)$ is a right and left coherent preadditive category.

\item If $R$ is complete local, then $\inj(R)$ is right Noetherian.

\item $\inj(R)$ is never left Noetherian.
\end{enumerate}
\end{prop}

\begin{proof}
In either case $\widehat{R} \otimes_R R_{\min}$ and $Q(\widehat{R})$ are flat $\widehat{R}$-modules. Hence the rings $A,B$ in Proposition \ref{oneobjegs} are right coherent by \cite[Theorem 4.2]{HMV12}. Then (1) follows from this and Theorem \ref{weakdimscheme}.

If $R$ is a complete local one-dimensional Noetherian ring, then  $\widehat{R} \otimes_R R_{\min} \simeq R_{\min}$ is a finitely generated $R_{min}$ module and so $A$ is right Noetherian. Similarly $Q(\widehat{R})$ is finitely generated $Q(\widehat{R})$ module module and so $B$ is right Noetherian.

$A$ is left Noetherian if and only if $\widehat{R} \otimes_R R_{\min}$ is a finitely generated $\widehat{R}$-module. As $\widehat{R}$ is faithfully flat, this implies $R_{\min}$ is finitely generated over $R$ which can only occur if $R$ has dimension zero. Now $B$ is left Noetherian if and only if $Q(\widehat{R})$ is finitely generated over $\widehat{R}$. This implies $\widehat{R}$ is also zero dimensional. 
\end{proof}

\begin{remark} \label{notnoetheg2}
If $R$ is not complete local, then we doubt that $\inj(R)$ can be right Noetherian. For example, if $R$ is a discrete valuation ring, $\inj(R)$ is right Noetherian if and only if $\dim_{\Frac(R)} \Frac(\widehat{R}) < \infty$. This fails for cardinality reasons in most cases. 
\end{remark}

If $R$ is a discrete valuation ring, Theorem \ref{weakdimscheme} shows that the weak dimension of $\inj(R)$ is equal to $1$. Below we consider the right and left global dimensions.

\begin{prop} \label{dvreg}
Suppose $R$ is a discrete valuation ring and let $K = \Frac(R)$ and $\widehat{K} = \Frac(\widehat{R})$ and let $A$ denote the ring. 
\[
\begin{pmatrix}
K & 0 \\
\widehat{K} & \widehat{R} 
\end{pmatrix}
\]
Then $\Mod \inj(R) \simeq \Mod A$ and $A$ right hereditary, left semihereditary and never left hereditary. The right global dimension of $\inj(R)$ is 1 and the left global dimension is strictly greater than 1.
\end{prop}

\begin{proof}
The description of $A$ follows from Proposition \ref{oneobjegs}. Using the classification of finitely generated modules over discrete valuation rings, one can show that every ideal of $A$ is isomorphic to a direct sum of copies of $P_1$ and $P_2$ where 
\[
P_1 \coloneq \begin{pmatrix}
K & 0 \\
0 & 0
\end{pmatrix} \quad P_2 \coloneq \begin{pmatrix}
0 & 0 \\
\widehat{K} & \widehat{R}
\end{pmatrix}
\]
and so all right ideals of $A$ are projective and $A$ is right hereditary. As $A$ has weak dimension one by Theorem \ref{weakdimscheme}, $\inj(R)$ is left semihereditary (alternatively one can compute all the finitely generated left ideals and show they are projective). Consider the left ideals
\[
X \coloneq  \begin{pmatrix}
0 & 0 \\
\widehat{K} & 0
\end{pmatrix} \quad P \coloneq \begin{pmatrix}
0 & 0 \\
0 & \widehat{R}
\end{pmatrix}.
\]
If $X$ were projective, the surjective map $\bigoplus_{\widehat{K}} P \to X$ would split. But this cannot be the case since $\Hom_{\widehat{R}}(\widehat{K},\widehat{R}) = 0$. Therefore $A$ is not left hereditary. 
\end{proof}

\begin{remark}
Proposition \ref{dvreg} shows that one cannot replace weak dimension with left global dimension in the statement of Theorem \ref{weakdimscheme} since $A$ has left global dimension strictly greater than one. 
\end{remark}

Finally we note that there are more semi-orthogonal decompositions than those which appear in Lemma \ref{SODlemma}.

\begin{prop} \label{SODeg}
Suppose $R$ is a complete local Gorenstein one-dimensional ring. There is a semi-orthogonal decomposition
\[
\mathcal{D}^{\immac}(R) = \langle R, E(R) \rangle
\]
\end{prop}

\begin{proof}
The minimal injective resolution of $R$ implies that $R$ and $E(R)$ collectively generate. To see that there are no morphisms from right to left we compute $\Ext^i_R(E(R),R)$ using the minimal injective resolution of $R$ to give
\[
0 \to \Hom_R(E(R),E(R)) \to \Hom_R(E(R),E(R/m)) \to 0
\]
Lemma \ref{gorensteinhom} allows us to identify this map with $\prod_{p \in \min(R)} R_p \to \prod_{p \in \min(R)} R_p \otimes_R \widehat{R}$ which is an isomorphism since $R = \widehat{R}$. Hence $\Ext^i(E(R),R) = 0$ for all $i$ as required. 
\end{proof}

\section{A Variant: Maximal immaculates}
\label{mimmacsec}

In this section we consider a variant of the immaculate complexes which we call the maximal immaculates. We relate strong generation of this category to a classically studied homological invariant for a wide class of noncommutative rings.

Let $\mathcal{A}$ denote a Grothendieck category. If $S$ is a simple object then the injective hull of $S$ is indecomposable by \cite[Corollary 1]{Gab58}. Let $\minj \subseteq \inj$ denote the injective hulls of the simples. 

\begin{defn}
A complex $I \in \mathcal{D}^{\immac}(\mathcal{A})$ is maximal immaculate if it is a bounded complex of finite sums of injective hulls of the simple modules. 
Let $\mathcal{D}^{\mimmac}(\mathcal{A}) \subseteq \mathcal{D}^{\immac}(\mathcal{A})$ denote the subcategory consisting of the maximal immaculates. 
\end{defn}

\begin{remark}
A similar proof to Proposition \ref{immacsarethick}, shows that 
\[
\mathcal{D}^{\mimmac}(X) = \thick(E(S) \mid S \in \mathcal{A} \text{ simple}).
\]
\end{remark}

\begin{ex} \label{Dmimscheme}
If $X$ is a locally Noetherian scheme, then as in Theorem \ref{SODfiltration}
\[
\mathcal{D}^{\mimmac}(X) = \thick(E(x)\mid x \text{ closed} ) \simeq \bigoplus_{x \text{ closed}} \mathcal{D}^{\perf}(\widehat{\mathcal{O}}_{X,x}).
\]
\end{ex}

Finitely embedded modules were introduced in \cite{Vam68} as a dual notion to finitely generated modules. The following invariant was studied in \cite{Fac82}.

\begin{defn} \label{deltainvariant}
An object $M \in \mathcal{A}$ is finitely embedded if $E(M)$ is a finite sum of injective hulls of simples. Let 
\[
\delta(\mathcal{A}) \coloneq \sup\{ \idim(M) \mid M \text{ finitely embedded} \}.
\]
\end{defn}

\begin{ex} 
If $R$ is a commutative ring, then $\delta(R) = 0$ if and only if $R$ is Von Neumann regular. A characterisation of commutative rings with $\delta(R) = 1$ is given in \cite{Fac82}. 
\end{ex}

The associated notion of being Noetherian was studied in \cite{Jan69}.

\begin{defn}
We say $\mathcal{A}$ is co-Noetherian if the quotient of any finitely embedded module is finitely embedded. 
\end{defn}

\begin{ex} \
\begin{enumerate}
\item By \cite[Theorem 1]{Jan69} and \cite[Theorem 2]{Vam68}, a commutative ring is co-Noetherian if and only if $R_m$ is Noetherian for every maximal ideal $m$. 

\item \cite[Corollary 2.3]{HIR00} A module-finite algebra over a commutative Noetherian ring is co-Noetherian.

\item In \cite[Example 2.1]{HIR00} The Weyl algebra over a characteristic zero field is shown not to be co-Noetherian. 
\end{enumerate}
\end{ex}

Similar arguments to Section \ref{immacsec} apply to any subset of the indecomposable injectives and prove the following.

\begin{theorem} 
Suppose $\mathcal{A}$ is co-Noetherian, then $\mathcal{D}^{\mimmac}(\mathcal{A}) = \thick_n(\minj)$ if and only if $\delta(\mathcal{A}) \leq n$.  
\end{theorem}

The following result also follows from the equivalence in Example \ref{Dmimscheme}.

\begin{theorem}
Suppose $X$ is a locally Noetherian scheme. Then $X$ is regular of dimension at most $n$ if and only if $\mathcal{D}^{\mimmac}(X) = \thick_n(\minj)$.
\end{theorem}

\appendix

\section{The immaculate complexes form a thick subcategory} \label{thickimmac}

We show the category of immaculates over a Grothendieck category is thick. The proofs in the section are dual versions to standard arguments. See \cite[Tag 064N]{stacks-project}.

\begin{defn} For $M \in \mathcal{D}(\mathcal{A})$ and $n \in \mathbb{Z}$, we say $M$ is $n$-$\inj$-pseudocoherent if there exists an immaculate complex $I$ and a map $\alpha: M \to I$ such that $H^i(\alpha)$ is an isomorphism for $i < n$  and $H^n(\alpha)$ is a monomorphism. 
\end{defn}

\begin{remark}
Note that if $M$ is $n$-$\inj$-pseudocoherent, then $\Sigma M$ is $n-1$-$\inj$-pseudocoherent. 
\end{remark}

\begin{lemma} \label{triainjfin} Suppose there is a triangle $X \to Y \to Z \in \mathcal{D}(\mathcal{A})$.
\begin{enumerate}
\item If $X$ and $Y$ are $m$-$\inj$-pseudocoherent, then $Z$ is $m-1$-$\inj$-pseudocoherent.

\item If $X$ is $m$-$\inj$-pseudocoherent and $Z$ is $m-1$-$\inj$-pseudocoherent then $Z$ is $m-1$-$\inj$-pseudocoherent.

\item If $Y$ and $Z$ are $m$-$\inj$-pseudocoherent then $X$ is $m$--$\inj$-pseudocoherent.
\end{enumerate}
\end{lemma}

\begin{proof}
We prove (1), the others follow from rotating triangles. Suppose $\alpha: X \to I$ and $\beta: Y \to J$ are such that $H^i(\alpha)$ and $H^i(\beta)$ are isomorphisms in degrees $i < m$ and monomorphisms for $i = m$, and $I$ and $J$ are immaculate. By composing with the map $J \to \sigma_{\leq m} J$ we can assume that $J^i= 0$ for $i > m$. Let $C$ be the cocone of $\alpha$, then the long exact sequence implies that $H^i(C) = 0$ for $i \leq m$. Therefore $C \simeq \tau_{\geq m+1} C$ which is zero in degrees less than or equal to $m$. Since $J$ is $K$-injective, we have
\[
\Hom_{\mathcal{D}(\mathcal{A})}(C,J) = \Hom_{K(\mathcal{A})}(\tau_{\geq m+1}C,J) = 0
\]
So the composite $C \to X \to Y \to J$ vanishes and there is a map $I \to J$ which fits into the following diagram.
\[
\begin{tikzcd}
X \arrow[r] \arrow[d] & Y \arrow[d] \arrow[r] & Z \arrow[d,"\gamma"] \\
I \arrow[r] & J \arrow[r] & K
\end{tikzcd}
\]
The long exact sequence and the five Lemma imply that $H^i(\gamma)$ is an isomorphism for $i \leq m-2$ and the four lemma implies that $H^{m-1}(\gamma)$ is a monomorphism. Since the map $I \to J$ is between $K$-injectives we can model it as a chain map. So we can take $K$ to be its mapping cone and so $K$ is a bounded complex of finite direct sums of indecomposable injectives.
\end{proof}

\begin{lemma} \label{injfinlemma}
Suppose $X$ is $m$-$\inj$-pseudocoherent and $H^i(X) = 0$ for $i < m$ then the injective hull of $H^m(X)$ is a finite sum of indecomposable injectives.
\end{lemma}

\begin{proof}
There is an immaculate complex $I$ and a map $\alpha: X \to I$ be such that $H^i(\alpha)$ is an isomorphism for $i < m$ and $H^m(\alpha)$ is mono. Then $H^m(X) \subseteq H^m(I)$ and so it is enough to show that $H^m(I)$ embeds into a finite sum of indecomposable injectives. Note we have that $H^i(I) = 0$ for $i < m$. Suppose $n$ is the least such that $I^n \neq 0$. If $ n > m$ we have that $H^m(I) = 0$ and we are done. If $n < m$ then $0 = H^n(I) = \ker(d^n)$. It follows that $I^n \to d(I^n)$ is a null homotopic summand of $I$ and so we may replace $I$ with $I'= I/(I^n + d(I^n))$. Now the least $n'$ such that $(I')^{n'} \neq 0$ is greater than $n$. So by repeating this we may assume that eventually $n = m$. In this case we have that $H^m(I) = \ker(d^m) \subseteq I^m$ as required.  
\end{proof}

\begin{lemma} \label{injfinchar}
$X \in \mathcal{D}(\mathcal{A})$ is $\inj$-pseudocoherent if and only if it is $m$-$\inj$-pseudocoherent for every $m$.
\end{lemma}

\begin{proof}
If $X$ is $\inj$-pseudocoherent then there is a quasi-isomorphism $X \to I$ where $I$ is a bounded below complex of finite sums of indecomposable injectives. For any $m$, the map $X \to I \to \sigma_{\geq m}I$ is an isomorphism on cohomology in degrees $i < m$ and a monomorphism at $i =m$. 

Suppose $X$ is $m$-$\inj$-pseudocoherent finite for every $m$. Then $X$ has cohomology bounded below by some $N$. Suppose that for some $n \in \mathbb{Z}$, we have constructed a chain map $\alpha: X \to I$ where $I^i = 0$ for $i > n$ and $i < N$ and each $I^i$ is a finite direct sum of indecomposable injectives such that the map $H^i(\alpha)$ is an isomorphism for $i < n$ and that $H^n(\alpha)$ is a monomorphism. Our base case is $n = N-1$ where $I^i = 0$ for all $i$. Let $C$ denote the mapping cocone of $\alpha$ then the long exact sequence implies that $H^i(C) = 0$ for $i \leq n$. Since $I$ is immaculate it is $\inj$-pseudocoherent and so it is $m$-$\inj$-pseudocoherent for every $m$. By assumption $X$ is $m$-$\inj$-pseudocoherent for every $m$. So by Lemma \ref{triainjfin},  $C$ is $m$-$\inj$-pseudocoherent for every $m$. In particular it is $n+1$-$\inj$-pseudocoherent. So Lemma \ref{injfinlemma} implies that there is a map $H^{n+1}(C) \hookrightarrow I^{n+1}$ where $I^{n+1}$ is a finite sum of indecomposable injectives. By Lemma \ref{injlift} we can lift this to a morphism $C \to \Sigma^{-{n+1}}I^{n+1}$. Since $\Sigma^{-{n+1}}I^{n+1}$ is $K$-injective, we can assume it is a chain map and so the composite $C^n \xrightarrow{d^n_C} C^{n+1} \to I^{n+1}$ must vanish. Now $C^{n+1} = X^{n+1} \oplus I^n$ and the map $C^{n+1} \to I^{n+1}$ is of the form $(\alpha^{n+1},-d^n): X^{n+1} \oplus I^n \to I^{n+1}$ for some $\alpha^{n+1}: X^{n+1} \to I^{n+1}$ and $d^n: I^n \to I^{n+1}$. Since the composite mentioned above vanishes, $d^{n+1}d^n_C = 0$ and $\alpha^{n+1} d^n = d^{n+1} \alpha^n$ as required. It follows from a long exact sequence argument that $\alpha$ has the right properties on cohomology. 
\end{proof}

\begin{lemma} \label{summand}
If $X \oplus Y$ is $m$-$\inj$-pseudocoherent then so is $Y$.
\end{lemma}

\begin{proof} Since $X \oplus Y$ has bounded below cohomology so does $Y$. So there is some $n >> 0$ such that $\Sigma^{-n} Y$ is $m$-$\inj$-pseudocoherent. Indeed shift $Y$ to the right enough so that $H^i(Y) = 0$ for $i \leq m$ then the zero map ensures that $Y$ is $m$-$\inj$-pseudocoherent. Now there is a triangle $X \oplus Y \to X \oplus Y \to Y \oplus \Sigma Y$ so by Lemma \ref{triainjfin}, $Y \oplus \Sigma Y$ is $m-1$-$\inj$-pseudocoherent. Hence $\Sigma^{-1} Y \oplus Y$ is $m$-$\inj$-pseudocoherent. Repeating this argument to the summand $\Sigma^{-1} Y$ we deduce that $\Sigma^{-2} Y \oplus \Sigma^{-1} Y$ is $m$-$\inj$-pseudocoherent. Repeating this we see that $\Sigma^{-n} Y \oplus \Sigma^{-n+1 }Y$ is $m$-$\inj$-pseudocoherent. So the triangle 
\[
\Sigma^{-n+1} Y \to \Sigma^{-n} Y \oplus \Sigma^{-n+1} Y \to \Sigma^{-n} Y
\]
and (3) in Lemma \ref{triainjfin}, imply that $\Sigma^{-n+1} Y$ is $m$-$\inj$-pseudocoherent. Repeating this argument $n$ times we see that $Y$ is $m$-$\inj$-pseudocoherent.
\end{proof}

\begin{cor} \label{injfinthick}
The categories $\mathcal{D}_{\ipscoh}(\mathcal{A})$ and $\mathcal{D}^b_{\ipscoh}(\mathcal{A})$ are thick
\end{cor}

\begin{proof}
Lemmas \ref{summand} and \ref{injfinchar} imply that $\mathcal{D}_{\ipscoh}(\mathcal{A})$ is closed under summands. If $f$ is a morphism in $\mathcal{D}_{\ipscoh}(\mathcal{A})$ then since it is between $K$-injectives we may assume it is a chain map and model its cone as a mapping cone. Then it is clear that the cone is also a bounded below complex of finite sums of indecomposables. Clearly, $\mathcal{D}_{\ipscoh}(\mathcal{A})$ is closed under shifts and so we are done. The intersection of two thick subcategories is thick and so the second claim follows. 
\end{proof}

\begin{lemma} \label{immacchar}
A complex $N \in \mathcal{D}(\mathcal{A})$ is immaculate if and only if $N \in \mathcal{D}^b_{\ipscoh}(\mathcal{A})$ and there exists a $k \in \mathbb{Z}$ such that $\Hom_{\mathcal{D}(\mathcal{A})}(M,\Sigma^i N) = 0$ for $i >k$ for every $M \in \mathcal{A}$. 
\end{lemma}

\begin{proof}
If $N$ is immaculate then it is clearly contained in $\mathcal{D}^b_{\ipscoh}(\mathcal{A})$. Assume $N$ is a complex of injectives with $N^i = 0$ for $i > k$. Since $N$ is $K$-injective
\[
\Hom_{\mathcal{D}(\mathcal{A})}(M,\Sigma^i N) = \Hom_{K(\mathcal{A})}(M,\Sigma^i N)
\]
which vanishes for $i > k$ since $M$ is concentrated in degree zero.

Conversely, suppose $N$ is a bounded below complex of finite sums of indecomposable injectives with $H^i(N) = 0$ for $i > m$. Let $n > \max\{k,m\}$. Then there is a triangle $\sigma_{\geq n} N \to N \to \sigma_{< n} N$ and $\sigma_{\geq n} N \simeq \Sigma^{-n} H^n(\sigma_{\geq n} N)$. Then,
\[
\Hom_{\mathcal{D}(\mathcal{A})}(\sigma_{\geq n} N, N) = \Hom_{\mathcal{D}(\mathcal{A})}(H^n(\sigma_{\geq n} N),\Sigma^{n} N) = 0 
\]
Therefore the natural map $\sigma_{\leq n}N \to N$ is null homotopic. So there exists a map $s: N^n \to N^{n-1}$ such that $d^{n-1}s = \id_{\im(d^{n-1})}$. Therefore $\im(d^{n-1})$ is a summand of $N^n$ and so is a finite sum of indecomposable injectives. Therefore $N \simeq \tau_{\leq n} N$ is a bounded complex of finite sums of indecomposable injectives as required.
\end{proof}

\begin{cor}
The immaculate complexes form a thick subcategory. 
\end{cor}

\begin{proof}
The second condition in Lemma \ref{immacchar} is closed under summands, shifts and triangles since $\mathcal{D}^b_{\ipscoh}(\mathcal{A})$ is thick by Corollary \ref{injfinthick}.
\end{proof}

\begin{prop} \label{immacsarethick}
The smallest triangulated subcategory of $\mathcal{D}(\mathcal{A})$ containing the indecomposable injectives is $\mathcal{D}^{\immac}(\mathcal{A})$ and $\mathcal{D}^{\immac}(\mathcal{A}) = \thick(\inj)$. 
\end{prop}

\begin{proof}
Since $\inj \subseteq \mathcal{D}^{\immac}(\mathcal{A})$ and $\mathcal{D}^{\immac}(\mathcal{A})$ is triangulated, the smallest triangulated subcategory containing $\inj$ is contained in $\mathcal{D}^{\immac}(\mathcal{A})$. If $X$ is immaculate then by taking brutal truncations, one sees that $X$ is in the triangulated subcategory generated by $\inj$. Then the last claim follows since $\mathcal{D}^{\immac}(\mathcal{A})$ is thick.
\end{proof}

\bibliographystyle{alpha}
\bibliography{biblio}

\end{document}